\documentclass[10pt, dvipsnames]{article}

\usepackage{enumitem}
\usepackage[dvipsnames]{xcolor}
\usepackage{ucs}
\usepackage{longtable}
\usepackage[colorlinks=true,
  urlcolor=ForestGreen,
  linkcolor=Blue,
  citecolor=Violet,
  pagebackref=true,
  pdftitle={Titel},
  pdfsubject={},
  pdfauthor={},
  pdfkeywords={keywords}
]{hyperref}
\usepackage{graphicx}
\usepackage{subfigure}
\usepackage[subfigure]{tocloft}
\usepackage[font=footnotesize]{caption}
\usepackage{tabularx}
\usepackage{nameref}
\usepackage{tikz-cd}

\usepackage{amsmath} 	 %MATH1
\usepackage{amssymb}	 %MATH2
\usepackage{amsthm}	 %MATH3
\usepackage{mathtools}   %MATH4
\usepackage[normalem]{ulem} 
\usepackage{graphicx,wrapfig,lipsum}
\usepackage{float}
\usepackage{floatflt}

\usepackage[most]{tcolorbox}
\tcbuselibrary{breakable} %so that tcolorboxes can go over several pages
\usepackage{wrapfig}
\usepackage[toc,page]{appendix}

\usepackage[parfill]{parskip}

\newtcolorbox[auto counter, number within=section
                            ]{blueproposition}[1][]{enhanced jigsaw,
  colback=white, %white!60!blue
  coltext={black},
  coltitle={black},
  boxrule=0pt,
  arc=1mm,
  frame hidden,
  auto outer arc,
  boxsep=5pt,
  left=2pt,
  right=2pt,
  bottom=2pt,
  top=2pt,
  borderline west={1mm}{0mm}{white!60!blue},
  before skip=3mm,
  after skip=3mm,
  breakable,
  title={Proposition~\thetcbcounter.},
  attach title to upper=\quad,
  fonttitle=\bfseries,
  #1}

\newtcolorbox[auto counter, number within=section
  ]{bluecorollary}[1][]{enhanced jigsaw,
colback=white, %white!60!blue
coltext={black},
coltitle={black},
boxrule=0pt,
arc=1mm,
frame hidden,
auto outer arc,
boxsep=5pt,
left=2pt,
right=2pt,
bottom=2pt,
top=2pt,
borderline west={1mm}{0mm}{white!60!blue},
before skip=3mm,
after skip=3mm,
breakable,
title={Corollary~\thetcbcounter.},
attach title to upper=\quad,
fonttitle=\bfseries,
#1}

\newtcolorbox[auto counter, number within=section
                            ]{bluetheorem}[1][]{enhanced jigsaw,
  colback=white, %white!60!blue
  coltext={black},
  coltitle={black},
  boxrule=0pt,
  arc=1mm,
  frame hidden,
  auto outer arc,
  boxsep=5pt,
  left=2pt,
  right=2pt,
  bottom=2pt,
  top=2pt,
  borderline west={1mm}{0mm}{white!60!blue},
  before skip=3mm,
  after skip=3mm,
  breakable,
  title={Theorem~\thetcbcounter.},
  attach title to upper=\quad,
  fonttitle=\bfseries,
  #1}

\newtcolorbox[auto counter, number within=section
                            ]{redlemma}[1][]{enhanced jigsaw,
  colback=white, %white!60!blue
  coltext={black},
  coltitle={black},
  boxrule=0pt,
  arc=1mm,
  frame hidden,
  auto outer arc,
  boxsep=5pt,
  left=2pt,
  right=2pt,
  bottom=2pt,
  top=2pt,
  borderline west={1mm}{0mm}{white!60!red},
  before skip=3mm,
  after skip=3mm,
  breakable,
  title={Lemma~\thetcbcounter.},
  attach title to upper=\quad,
  fonttitle=\bfseries,
  #1}
  
\newtcolorbox[auto counter, number within=section
                          ]{debox}[1][]{enhanced jigsaw,
  colback=green!0!white,
  coltext={black},
  colframe={black},
  coltitle={black},
  boxrule=0pt,
  frame hidden,
  borderline west={1mm}{0mm}{white!50!green},
  arc=0mm,
  auto outer arc,
  boxsep=5pt,
  left=4pt,
  breakable,
  right=4pt,
  bottom=0pt,
  top=0pt,
  before skip=3mm,
  after skip=3mm,
  title={Definition~\thetcbcounter.},
  attach title to upper=\quad,
  fonttitle=\bfseries,
  #1}
\newtcolorbox[auto counter,number within=section
                          ]{yellowRemark}[1][breakable]{enhanced jigsaw,
  colback=green!0!white,
  coltext={black},
  colframe={black},
  coltitle={black},
  boxrule=0pt,
  frame hidden,
  breakable,
  borderline west={1mm}{0mm}{white!50!yellow},
  arc=0mm,
  auto outer arc,
  boxsep=5pt,
  left=4pt,
  right=4pt,
  bottom=0pt,
  top=0pt,
  before skip=3mm,
  after skip=3mm,
  title={Remark~\thetcbcounter.},
  attach title to upper=\quad,
  fonttitle=\bfseries,
  #1}
\newtcolorbox[auto counter,number within=section
                          ]{example}[1][breakable]{enhanced jigsaw,
  colback=green!0!white,
  coltext={black},
  colframe={black},
  coltitle={black},
  boxrule=0pt,
  frame hidden,
  breakable,
  borderline west={1mm}{0mm}{white!50!yellow},
  arc=0mm,
  auto outer arc,
  boxsep=5pt,
  left=4pt,
  right=4pt,
  bottom=0pt,
  top=0pt,
  before skip=3mm,
  after skip=3mm,
  title={Example~\thetcbcounter.},
  attach title to upper=\quad,
  fonttitle=\bfseries,
  #1}

  \newtcolorbox[auto counter,number within=section
                          ]{proofenv}[1][breakable]{enhanced jigsaw,
  colback=black!0!white,
  coltext={black},
  colframe={black},
  coltitle={black},
  boxrule=0pt,
  frame hidden,
  breakable,
  borderline west={1mm}{0mm}{white!50!black},
  arc=0mm,
  auto outer arc,
  boxsep=5pt,
  left=4pt,
  right=4pt,
  bottom=0pt,
  top=0pt,
  before skip=3mm,
  after skip=3mm,
  title={Proof.~},
  attach title to upper=\quad,
  fonttitle=\bfseries,
  #1}

  \newtcolorbox[auto counter,number within=section
                          ]{proofexp}[1][breakable]{enhanced jigsaw,
  colback=black!0!white,
  coltext={black},
  colframe={black},
  coltitle={black},
  boxrule=0pt,
  frame hidden,
  breakable,
  borderline west={1mm}{0mm}{white!50!black},
  arc=0mm,
  auto outer arc,
  boxsep=5pt,
  left=4pt,
  right=4pt,
  bottom=0pt,
  top=0pt,
  before skip=3mm,
  after skip=3mm,
  title={Proof and Explanation.~},
  attach title to upper=\quad,
  fonttitle=\bfseries,
  #1}
\usepackage{pictexwd,dcpic}

\newtcolorbox[% within=section
                          ]{shortsummary}[1][breakable]{enhanced jigsaw,
  colback=red!0!white,
  coltext={black},
  colframe={black},
  coltitle={black},
  boxrule=0pt,
  frame hidden,
  breakable,
  borderline west={1mm}{0mm}{white!50!red},
  arc=0mm,
  auto outer arc,
  boxsep=5pt,
  left=4pt,
  right=4pt,
  bottom=0pt,
  top=0pt,
  before skip=3mm,
  after skip=3mm,
  title={Short and important~\thetcbcounter:},
  attach title to upper=\quad,
  fonttitle=\bfseries,
  #1}

\usepackage{changepage}

\usepackage{geometry}
\def\wider{%
   \advance\leftskip -\parindent
   \advance\rightskip -\parindent}

\newcommand{\Ra}{\mathrel{ \Rightarrow }}

\newcommand{\ra}{\mathrel{ \rightarrow }}
\newcommand{\Lr}{\mathrel{ \Leftrightarrow }}

\newcommand{\te}[1]{\text{#1}}

\newcommand{\ov}[1]{\mathrel{ \overset{ \text{(#1)} }{=}  }}
\newcommand{\br}[1]{\mathrel{ \left(#1\right)  }}

\newcommand{\bp}[2]{\begin{blueproposition}[label=#1]#2\end{blueproposition}}
\newcommand{\bc}[2]{\begin{bluecorollary}[label=#1]#2\end{bluecorollary}}

\newcommand{\lm}[2]{\begin{redlemma}[label=#1]#2\end{redlemma}}
\newcommand{\gd}[2]{\begin{debox}[label=#1]#2\end{debox}}
\newcommand{\yr}[2]{\begin{yellowRemark}[label=#1]#2\end{yellowRemark}}
\newcommand{\yer}[2]{\begin{example}[label=#1]#2\end{example}}

\newcommand{\re}{\text{Re}}

\newcommand{\ovs}[1]{\mathrel{ \overset{ \text{(#1)} }{\simeq}  }}
\newcommand{\bo}[1]{\mathbf{#1}}
\newcommand{\fc}[1]{\widehat{\mathbf{#1}}}

\newcommand{\pr}[1]{\begin{proof}#1\end{proof}}
\newcommand{\prr}[1]{\pr{See Proof $~$\ref{proof:#1}.}}

\newcommand{\bt}[1][normal]{\begin{tikzcd}[ampersand replacement = \&, column sep=#1,row sep=#1]}
\newcommand{\et}{\end{tikzcd}}

  \tikzset{%this is for the adjunction symbol, 
    symbol/.style={%
        draw=none,
        every to/.append style={%
            edge node={node [sloped, allow upside down, auto=false]{$#1$}}}
    }
}
\makeatletter
\newcommand{\tpitchfork}{%
    \raise-0.1ex
    \vbox{
    \baselineskip\z@skip
    \lineskip-.52ex
    \lineskiplimit\maxdimen
    \m@th
    \ialign{##\crcr\hidewidth\smash{$-$}\hidewidth\crcr$\pitchfork$\crcr}
  }%
}
\makeatother

\newcommand{\mytoc}[1]{\hyperref[#1]{\ref{#1} $ \qquad $ \nameref{#1} \dotfill\pageref{#1}}}

\newcommand{\Hom}{\mathrm{Hom}}

\usepackage{tocloft}

\usepackage{newtxtext}
\usepackage{lipsum} 

\makeatletter
\def\moverlay{\mathpalette\mov@rlay}
\def\mov@rlay#1#2{\leavevmode\vtop{%
   \baselineskip\z@skip \lineskiplimit-\maxdimen
   \ialign{\hfil$\m@th#1##$\hfil\cr#2\crcr}}}
\newcommand{\charfusion}[3][\mathord]{
    #1{\ifx#1\mathop\vphantom{#2}\fi
        \mathpalette\mov@rlay{#2\cr#3}
      }
    \ifx#1\mathop\expandafter\displaylimits\fi}
\makeatother

\newcommand{\cupdot}{\charfusion[\mathbin]{\cup}{\cdot}}
\newcommand{\bigcupdot}{\charfusion[\mathop]{\bigcup}{\cdot}}

\begin{document}

\begin{center}
  \LARGE{Fuzzy simplicial sets and their application\\to geometric data analysis}\\~\\
  \large{Lukas Silvester Barth$^{1,*}$\qquad Hannaneh Fahimi$^{1,2}$\qquad Parvaneh Joharinad$^{1,2}$\\~\\ J\"urgen
  Jost$^{1,2,3,5}$\qquad Janis Keck$^{1,5,6}$\qquad Thomas Jan Mikhail$^{4,7}$}
  \\~\\ 
  \normalsize{$^1$Max Planck Institute for Mathematics in the Sciences, Leipzig, Germany\\
  $^2$Center for Scalable Data Analytics and Artificial Intelligence (ScaDS.AI), Dresden/Leipzig,
  Germany\\
  $^3$Santa Fe Institute for the Sciences of Complexity, New Mexico, USA\\
  $^4$Universitat Autònoma de Barcelona, Barcelona, Spain\\
  $^5$Max Planck Institute for Human Cognitive and Brain Sciences, Leipzig, Germany\\
  $^6$Max Planck School of Cognition\\
  $^7$University of Copenhagen, Copenhagen, Denmark\\
  $^*$Corresponding author\\~\\
  \href{mailto:lukas.barth@mis.mpg.de}{\color{Blue}{Lukas.Barth@mis.mpg.de}}, \href{mailto:fatemeh.fahimi@mis.mpg.de}{\color{Blue}{Fatemeh.Fahimi@mis.mpg.de}}, \href{mailto:Parvaneh.JohariNad@mis.mpg.de}{\color{Blue}{Parvaneh.JohariNad@mis.mpg.de}}, \href{mailto:jjost@mis.mpg.de}{\color{Blue}{Jost@mis.mpg.de}}, \href{mailto:Janis.Keck@maxplanckschools.de}{\color{Blue}{Janis.Keck@maxplanckschools.de}}, \href{mailto:tjm@math.ku.dk}{\color{Blue}{TJM@math.ku.dk}}
  }
  \end{center}

  \medskip 
  
\begin{adjustwidth}{10mm}{10mm}
  \small
  In this article, we expand upon the concepts introduced in \cite{Spivak09} about the relationship between the category $\bo{UM}$ of uber metric spaces and the category $\bo{sFuz}$ of fuzzy simplicial sets. We show that fuzzy simplicial sets can be regarded as natural combinatorial generalizations of metric relations. Furthermore, we take inspiration from UMAP (cf.~\cite{McInnes18}) to apply the theory to manifold learning, dimension reduction and data visualization, while refining some of their constructions. We generalize the adjunction between $\bo{UM}$ and $\bo{sFuz}$, derive an explicit description of colimits in $\bo{UM}$, and show that $\bo{UM}$ can be embedded into $\bo{sFuz}$. Furthermore, we prove analogous results for the category of extended-pseudo metric spaces $\bo{EPMet}$. We also provide rigorous definitions of functors that make it possible to recursively merge sets of fuzzy simplicial sets and provide a description of the adjunctions between the category of truncated fuzzy simplicial sets and $\bo{sFuz}$, 
  which we relate to persistent homology.
  Combining those constructions, we can show a surprising connection between the well-known dimension reduction methods UMAP and Isomap (cf.~\cite{Tenenbaum00}) and derive an alternative algorithm, which we call IsUMap, that combines some of the strengths of both methods.\footnote{The source code is available on github: \href{https://github.com/LUK4S-B/IsUMap}{https://github.com/LUK4S-B/IsUMap}} Additionally, we developed a new embedding method that allows to preserve clusters detected in the original metric space that we construct from the data. The visualization of the optimization process gives the user information both about the inner-cluster distributions in the original metric space and their inter-cluster relations. 
  We compare our new method with UMAP, Isomap and t-SNE on a series of low- and high-dimensional datasets and provide explanations for observed differences and improvements.

\medskip

\textbf{Keywords:} Fuzzy simplicial sets, metric spaces, manifold learning, dimension reduction\\
\textbf{MSC classes:} 18-02, 18-04, 51F04, 03E72
\normalsize 

\medskip

\textbf{Contents}

\mytoc{sec:Introduction}\\
\mytoc{sec:DiscreteApprox}\\
\mytoc{sec:FuzzySimplicialSets}\\
\mytoc{sec:nFuz}\\
\mytoc{sec:Adjunction}\\
\mytoc{sec:TheMergeFunctors}\\
\mytoc{sec:application}\\
% \mytoc{sec:declarations}\\
\mytoc{sec:appendixx}\\

\end{adjustwidth}

\section{Introduction}
\label{sec:Introduction}

The article presents two main achievements:
\begin{enumerate}[label=(A.\arabic*)]
  \item We show that fuzzy simplicial sets, as introduced by \cite{Spivak09}, can serve as natural combinatorial objects to represent data, and clarify and extend the corresponding theory.
  \label{A1}
  \item We use those theoretical advancements to put UMAP developed in \cite{McInnes18} on a more solid footing, show how it is related to Isomap (cf.~\cite{Tenenbaum00}) and introduce a natural combination of both methods that we call IsUMap. Furthermore, we enrich IsUMap with strong and interpretable clustering capabilities through a cluster separatation technique.
  \label{A2}
\end{enumerate}
In \cite{Spivak09}, an adjunction between the category of uber metric spaces $\bo{UM}$ and the category of fuzzy simplicial sets $\bo{sFuz}$ was sketched in analogy to the adjunction between topological spaces and (ordinary) simplicial sets. Building up on Spivak's work, \cite{McInnes18} introduced a similar but more discrete version of this adjunction and made use of it in the context of manifold learning and data visualization. 

In the first part of this work, we extensively study and develop the categorical theory behind these constructions.
Since the above mentioned works often skip details and sometimes contain unresolved problems, we begin with a thorough treatment of the underlying theory and re-exhibit both adjunctions as appropriate special cases of the nerve-realization theorem (cf.~Section \ref{sec:FuzzySimplicialSets} and Propositions \ref{prop:KanExtension}, \ref{prop:singFactorsThroughSfuz} and \ref{prop:skeletonSingFactorsThroughSfuz}).

To aid subsequent theory as well as computational applications, we provide a more explicit geometric description of the colimit in the category of uber metric spaces and extended pseudo-metric spaces in Proposition \ref{prop:colimitUM}, which in turn allows us to provide more explicit descriptions of the realization functors of the above mentioned adjunctions in Propositions \ref{prop:metricRealizationOfIdentity}, \ref{prop:faceMap} and \ref{prop:realizationOfFuzzySimpSet}. 

In Section \ref{sec:nFuz}, we generalize the skeleton and co-skeleton adjunctions known from the theory of (ordinary) simplicial sets to the setting of truncated fuzzy simplicial sets and also derive computationally explicit descriptions for these functors in Proposition \ref{prop:leftRightAdjointOfTruncation}.

Using these results, we can prove one of the main points of our article 
in Proposition \ref{prop:UMisInSfuz}, namely that $\bo{UM}$ and the category of extended pseudo-metric spaces $\bo{EPMet}$ can be embedded into $\bo{sFuz}$. This means that fuzzy simplicial sets can indeed be used to represent arbitrary metric data.
At the same time, we show 
in Proposition \ref{prop:SingSimpSetProperties} 
that only under specific conditions a fuzzy simplicial set corresponds to a uber metric space. This reveals that fuzzy simplicial sets have strictly stronger data representation capabilities than metric spaces. Below we explain how this can help in practice. 

Even though data sets are often finite subsets $X$ of some Euclidean space $\mathbb{R}^n$, and hence come equipped with a metric structure, the additional freedom can be used to encode relationships (for example in a social network graph) that do not satisfy all metric relations, like the triangle inequality, in a fuzzy simplicial set. Furthermore, the data set $X$ is often assumed to be sampled from a distribution, whose probability mass is centered around some lower-dimensional geometric object (like a manifold, or an object, that is a manifold almost everywhere, up to a subset of points with Lebesgue measure zero) and in that case the metric structure of $\mathbb{R}^n$ is not necessarily intrinsic to the data $X$. To extract the intrinsic structure, geodesic distances are often approximated with graph distances on local $k$-nearest-neighborhood graphs, and graphs are specific instances of (truncated) fuzzy simplicial sets. 

At the same time, fuzzy simplicial sets go beyond (possibly directed and weighted) graphs because they allow to encode higher-order relationships between more than $2$ points, which a metric cannot account for, while providing a rich combinatorial and category-theoretical structure that alternative structures, like for example hypergraphs, cannot offer. One important example of such higher relationships are those that arise in the area of persistent homology (see, for example \cite{Carlsson09}), where a simplicial complex or simplicial set can encode the relationships that arise due to a filtration. We show in 
Proposition \ref{prop:simplifySpivaksIdea} 
that the formalization of this process in terms of fuzzy simplicial sets, that was introduced in \cite{Spivak09}, can be carried out in a more direct way using the above mentioned skeleton and co-skeleton adjunctions we derived. 

Another example of a higher-order relationship, that one might like to encode in a fuzzy simplicial set, is the metric curvature between triples of points in a data set, which was introduced in \cite{Joharinad19}, and can be thought of as a metric generalization of the notion of Riemannian curvature.

Overall, our propositions thus allow to consider fuzzy simplicial sets as an appropriate combinatorial object for generalized metric and higher-order relationships, in particular those that arise in the area of data analysis as claimed in \ref{A1}. Furthermore, the described contributions can be used as a basis for further theoretical development.

\bigskip

The second contribution \ref{A2} is motivated by the importance of manifold learning, dimension reduction and data visualization methods in the area of data analysis. Such methods are a signifcant tool to get an overview of large high-dimensional data sets and can guide the application of further methods. Moreover, if the data can be embedded into a space with reduced dimension without too much distortion, then the embedding can be used as input for downstream tasks in order to increase the computational efficiency. For example, the data that is used for the training of a neural network could be preprocessed by a dimension reduction method to facilitate faster training. 

Prominent dimension reduction and data embedding methods include Principal Component Analysis (PCA, cf.~\cite{Pearson1901}), Laplacian Eigenmaps (LE, cf.~\cite{Belkin03}), (classical or (non-)metric) Multidimensional Scaling (MDS, cf.~\cite{Torgerson1952}, \cite{borg05}), Isomap (cf. \cite{Tenenbaum00}), t-SNE (cf.~\cite{vandermaaten08a}), 
LargeVis (cf.~\cite{tang2016largeVis}), 
Uniform Manifold Approximation (UMAP, cf.~\cite{McInnes18}) 
and PHATE (cf.~\cite{moon2017phate}). We also note that a general category-theoretical perspective on different manifold learning schemes and how they relate to clustering algorithms is given in the references \cite{shiebler2020} and \cite{shiebler2020clustering}. 

Of those methods, UMAP is interesting from a category-theoretical point of view because the authors pick up and modify the adjunction 
$\bt\te{Sing}:\bo{UM}\ar{r}\&\ar{l}\bo{sFuz}:\te{Re}\et$
introduced in \cite{Spivak09} to justify (a substantial part of) their procedure. It is useful for the purposes of this article to summarize UMAP concisely as follows:
\begin{enumerate}[label=(U.\arabic*)]
  \item Given a finite data set $X\subset \mathbb{R}^n$, and a hyperparameter $k$, the $k$ nearest neighbors of each point are computed using the Euclidean distance in $\mathbb{R}^n$. Then a set $\mathcal{X}:=\{(X,d_i)\}_{i \in \{1,\cdots,|X|\}}$ of metric spaces is defined, the underlying set of each being a copy of $X$ and the metric $d_i$ describing the distances between the $i$th point of $X$ and its $k$ nearest neighbors, while other distances are set to infinity.
  \label{U1}
  \item The distances $d_i$ are further modified in two ways, namely the nearest neighbor distance is subtracted
  and the distances are normalized. The extraction of $k$-neighborhoods and the normalization are justified with a geometric lemma, that we shall henceforth simply call ``Lemma 1''.
  \label{U2}
  \item Then a modification of the $\te{Sing}$-functor of the above mentioned adjunction is used to map $\mathcal{X}$ to a set of fuzzy simplicial sets $\mathcal{S}:=\{\te{Sing}(X,d_i)\}_{i}$.
  \label{U3}
  \item Once the metric spaces are converted to fuzzy simplicial sets, probabilistic operations, namely so-called t-conorms, can be used to merge them. 
  
  In UMAP, this merge is however only performed on the underlying directed weighted graphs $G_i$ obtained by truncating the fuzzy simplicial sets, i.e.~$G_i=\te{ctr}_1(\te{Sing}(X,d_i))$, and then the merge operation yields a single symmetric graph $G=\te{merge}(\{G_i\}_i)$.
  \label{U4}
  \item In order to visualize the symmetric graph $G$, an optimization procedure is carried out, in which the positions of a set of points $Y\subset \mathbb{R}^n$ (where $n$ is usually $2$ or $3$), with $|Y|=|X|$, is modified with gradient descent until $G$ and a graph $H(Y)$, obtained from $Y$, are close in the sense that their cross-entropy $\mathcal{L}(Y):=-\sum_{i,j}\{G_{ij}\log(H(Y)_{ij})+(1-G_{ij})\log(1-H(Y)_{ij})\}$ is low. Usually a spectral embedding of $G$ into $\mathbb{R}^m$ is used as initialization for the gradient-descent procedure.
  \label{U5}
\end{enumerate}
As one can see, steps (U.3) and (U.4) make use of category-theoretical constructions. The underlying idea is that the application of the Sing-functor brings the metric relations into a category, where merge operations are more easily described, and where the optimization procedure (U.5) takes a more flexible form.
While studying the UMAP algorithm, several questions arose:
\begin{enumerate}[label=(Q.\arabic*)]
  \item Lemma 1 appearing in \cite[section 2.1]{McInnes18} shall help ``to approximate the manifold we assume
  the data (approximately) lies on''. However, the proof requires the assumption that the metric $g$ of the manifold is a constant diagonal matrix in ambient coordinates. This implies that the derivative of the metric tensor and hence also all Christoffel symbols and the Riemann curvature tensor vanish locally. This means that Lemma 1 is only applicable to data sets that lie on locally flat manifolds, which is a very restrictive assumption. Is it possible to justify the approximation of geodesic distances by merging $k$-neighborhoods for more general manifolds?
  \label{q1}
  \item The proof provided in the appendix of \cite{McInnes18} about the existence of the adjunction between finite extended pseudo metric spaces \textbf{FinEPMet} and finite fuzzy simplicial sets \textbf{Fin-sFuz} is missing a proof of the existence of the colimits that appear in the Kan extension. The reason is that \textbf{FinEPMet}, while being finitely complete, does not admit all small colimits (e.g.~an infinite coproduct of finite metric spaces need not be finite). Even though one can perhaps show that colimits over the category of elements of finite fuzzy simplicial sets exist in \textbf{FinEPMet} if one defines finite fuzzy simplicial sets to be those fuzzy simplicial sets $S$, where $S(n,0)$ is a finite set for each $n$, we consider it more elegant to reconsider their adjunction in a bigger category, for example in the category $\bo{UM}$ or $\bo{EPMet}$, and want to explore if this is a theoretically consistent alternative.
  \label{q2}
  \item The merge operation or ``fuzzy set union'' as it is called in \cite[section 2.2, below Def. 9]{McInnes18} is not defined in categorical language. In fact, they do not define it at all in the category $\bo{sFuz}$ but instead only apply a t-conorm to an adjacency matrix that jointly represents all the classical fuzzy graphs obtained from the objects in $\bo{sFuz}$. Can one rigorously define the merge operations in $\bo{sFuz}$?
  \label{q3}
  \item As one can see in \ref{U1} to \ref{U5}, the full adjunction is never used, in the sense that, even though the Sing functor is employed to map both $X$ and $Y$ to fuzzy simplicial sets, the Re functor never maps them back to a metric space, and the final form of $Y$ is only obtained indirectly via gradient descent. Since, however, the data is visualized geometrically, the result of the process must be a metric space $(Y,d_{\mathbb{R}^m})$ again, and therefore a natural question is whether one could make use of the realization functor to create an alternative path to obtain the embedding. Diagrammatically, one can visualize the two pathways as follows (not all arrows are functors, and the diagram does not necessarily commute because the two embedding procedures are different):
  \begin{equation}
    \begin{split}
      \begin{tikzcd}[ampersand replacement=\&,column sep=large]
          \bo{Met} \ar{r}{\text{split}}
          \& 
          \bo{UM}^N_* \ar[',dashed]{d}{\top_{\bo{UM}}} \ar{r}{\text{Sing}^N}
          \& 
          \bo{sFuz}^N_* \ar{r}{\text{ctr}_1^N}
          \&
          \bo{c1Fuz}^N_* \ar{r}{\te{merge}_{\bo{c1Fuz}}} 
          \&
          \bo{c1Fuz} \ar[']{d}{\text{graph embedding}}\\
          \&
          \bo{UM}\ar[dashed]{rrr}{\te{metric embedding}}
          \& 
          \& 
          \& 
          \bo{Euc}
      \end{tikzcd}
    \end{split}
    \label{diag:umapAlternative}
  \end{equation}
  Here $\te{split}$ is the operation described in \ref{U1} and \ref{U2}, $\te{Sing}$ corresponds to step \ref{U3}, $\te{merge}_{\bo{c1Fuz}}\circ \te{ctr}_1^N$ corresponds to \ref{U4} and the embedding arrow on the right side corresponds to \ref{U5}. $\te{Euc}$ denotes the category of Euclidean spaces.\footnote{The $*$ next to a category (for example $\bo{UM}^N_*$ or $\bo{sFuz}^N_*$) marks that we actually have to consider appropriate pullback categories on which $\te{merge}$ can operate, that we describe in Section \ref{sec:TheMergeFunctors}.}

  As one can see, even though the Sing functor is employed, the Re functor never maps the merged fuzzy simplicial set back to a metric space because the final form in $\bo{Euc}$ is obtained indirectly via gradient descent.\footnote{Since the form of the Sing-functor is agnostic to the choice of the metric category, UMAP in principle works both in $\mathbf{UM}$ and $\mathbf{EPMet}$.}
  Therefore a natural question is whether one could make use of the Re functor to walk along the alternative path that is visualized with dashed lines in Diagram \eqref{diag:umapAlternative}, where $\top_{\bo{UM}}:=\te{Re}\circ\te{merge}_{\bo{sFuz}}\circ\te{Sing}^N$ denotes the T-conorm combination (or T-combination, for short) of (uber) metric spaces, and use a metric embedding method to obtain an alternative embedding in $\bo{Euc}$. (We choose the name T-combination because t-conorms are ultimately defining $\te{merge}_{\bo{sFuz}}$ as we explain in Section \ref{sec:TheMergeFunctors}, which in turn define $\top_{\bo{UM}}$.)
  The question arising from this analysis is whether one can derive explicit formulas for the computations involved in this alternative path and implement an algorithm that makes use of them?
  \label{q4}
\end{enumerate}
To resolve \ref{q1}, we show in 
Section \ref{sec:DiscreteApprox}
that radial distances in Riemann normal coordinates are well approximated by Euclidean distances in small neighborhoods and combine this idea with convergence proofs provided in \cite{Bernstein00}. An important difference to Lemma 1 is, however, that this justification does not require the uniformity assumption or any kind of normalization. 

\ref{q2} is answered affirmatively by reconsidering the adjunction introduced in \cite{McInnes18} as one of all possible nerve-realization adjunctions between $\bo{UM}$ and $\bo{sFuz}$ or $\bo{EPMet}$ and $\bo{sFuz}$, the existence of which is shown in
Propositions \ref{prop:KanExtension} and \ref{prop:umapReIsFunWithCodUM}.
This reconsidered adjunction is the coskeleton adjunction that we already mentioned above in connection to persistent homology and 
Proposition \ref{prop:simplifySpivaksIdea}. 

We also provide a positive answer to \ref{q3} and describe general merge operations in $\bo{sFuz}$ in Section \ref{sec:TheMergeFunctors}. This allows to apply a t-conorm recursively across an entire family of fuzzy simplicial sets, generalizing the simpler graph case in a clean category-theoretical setting.

With those general merge operations at hand, and the derivation of an explicit formula for the Re functor of the skeleton adjunction provided in 
Proposition \ref{prop:realizationOfFuzzySimpSet}
(which in turn required a more explicit description of the colimit in $\bo{UM}$ provided in Proposition \ref{prop:colimitUM}),
we can finally provide an answer to \ref{q4}. In particular, Proposition \ref{prop:mergeUM} describes how to transfer the merge operation to the T-combination operation over the category $\bo{UM}$ and Section \ref{sec:application} provides the description of a concrete algorithm for computing the combination of metric spaces and the embedding into some Euclidean space for dimension reduction or data visualization purposes. It turns out that upon choosing a particularly simple t-conorm and a particularly simple, unnormalized metric for the local neighborhoods, the alternative pathway recovers the Isomap algorithm proposed in \cite{Tenenbaum00}. At the same time, our derivation shows how to naturally incorporate the local distances used in UMAP into the alternative pathway, which then yields a natural combination of the two methods.
Since our algorithm represents a combination of Isomap and UMAP, and takes place entirely in the category $\bo{UM}$, we decided to call it \textbf{IsUMap}. 

As an additional contribution, we enrich IsUMap with clustering capabilities vital to a wide range of applications by developing a new \textbf{Cluster Separation Optimization} method. This method combines our initial IsUMap embedding with a subsequent optimization procedure that uses high-dimensional information to pull clusters in the low-dimensional embedding apart. This procedure resolves the crowding problem (cf.~\cite{vandermaaten08a} and \cite{olszewski2025data}) that causes clusters to overlap in the embedding space of distance-preserving dimension reduction methods as we explain in Section \ref{sec:embedding}.
More specifically, after the merged simplicial sets have been geometrically realized in a metric space (which is not yet embedded into any ambient space), clusters can be determined based on the geodesic distances in that space before any information is destroyed through an embedding into a lower-dimensional space. Then, after embedding via metric multidimensional scaling, the overlapping clusters in the low-dimensional embedding are pulled apart via a stochastic gradient descent procedure that minimizes the overlap of the convex hulls of the points that were determined to be in different clusters before the embedding. This results in a clean separation of concerns (distance-preserving embedding and clustering occur separately) that yields very interpretable results. Moreover, our method allows to track and visualize the optimization path during the separation of the convex hulls in the low-dimensional space (as we show for example in Figures \ref{fig:dynamicCRCCbefore} to \ref{fig:dynamicBMMCLeiden}), which gives the user additional information about the original high-dimensional configuration.
The method can be freely combined with any clustering method applicable to geodesic distance matrices. We employ the Linkage 
\cite{Linkage2001,mullner2011linkage} and Leiden clustering \cite{Traag2019Leiden} in our simulations but others would work as well.

In Section \ref{sec:visualization}, we compare our new approach with UMAP, Isomap and t-SNE on a series of low- and high-dimensional datasets, ranging from low-dimensional manifolds over higher-dimensional datasets like Fashion-MNIST and newsgroup documents to real-world RNA sequencing datasets of bone marrow mononuclear cells and KRAS tumor organoid cells.
We demonstrate that our method often outperforms the others in terms of metric faithfulness, visual clarity and clustering capability.

At the end, in Section \ref{sec:effectsUMAP}, we provide a detailed analysis of the implicit effects taking place during the embedding performed in UMAP, which allows to better understand the differences between the visualizations of the various methods. While these implicit effects often force UMAP to display structures that are not inherent to the data but merely artifacts of the method, IsUMap adheres very closely to its theoretical description, which helps to interpret its visual results. The reader interested in further empirical experiments with additional datasets, can have a look at our related article \cite{VRpaper} and our book \cite{IsumapBook}.
\section{Geometric considerations and colimits in UM and EPMet}
\label{sec:DiscreteApprox}

\gd{def:uberMetricSpace}{
An uber-metric space $(X,d)$ is a set $X$ equipped with  a map $d:X\times X\to \mathbb{R}_{\ge 0}\cup \{\infty\}$ such that
\begin{enumerate}
\item $d(x,x) = 0$;
\item $d(x,y)=d(y,x)$;
\item $d(x,z)\le d(x,y)+d(y,z)$.
\end{enumerate}
The category of uber-metric spaces $\mathbf{UM}$ has as
objects uber-metric spaces and as morphisms non-expansive maps. 
}
A closely related concept, important in the context of \cite{McInnes18}, is a so-called \textbf{extended pseudo-metric space}, the definition of which is equivalent to an uber metric space, except that the 3rd condition is amended to ``3. Either $d(x,z)=\infty$ or $d(x,z)\le d(x,y)+d(y,z)$'', i.e.~the triangle inequality does not have to hold for infinite distances in an extended pseudo-metric space. We denote the category of extended pseudo-metric spaces and non-expansive maps by $\mathbf{EPMet}$.

As explained in \ref{q1}, we would like to replace \cite[Lemma 1]{McInnes18} with another justification for the approximation of geodesic distances of the manifold, that is assumed to approximately underly the data, even when the manifold is not assumed to be a locally flat submanifold of some Euclidean space. This involves aspects of intrinsic and extrinsic geometry. For the intrinsic geometry, locally (that is, within the injectivity radius) the distances from a point $p$ equal the Euclidean distances in a Riemann normal coordinate chart centered at $p$. The distances between other points are different, unless the metric is locally flat, but in fact deviate from the Euclidean ones in the normal coordinates only to second order w.r.t. the distance from $p$, with the approximation error controlled by the Riemann curvature tensor. But since the construction to be developed will only depend on the distances from $p$ of its $k$ nearest neighbors, this is not of much concern. When a smooth manifold $M$ is embedded into Euclidean space $\mathbb{R}^n$, it inherits a Riemannian metric as a submanifold of $\mathbb{R}^n$. The  distances w.r.t. that metric will, however,  because constrained to $M$, be different from and in fact larger than those of $\mathbb{R}^n$ unless $M$ is a linear subspace. The difference is controlled by the second fundamental tensor of $M$. More precisely, a geodesic on $M$ from  $p$ to some other point $q$ in the vicinity of $p$ is approximated to second order by a circle in $\mathbb{R}^n$ of some radius $r$ that is controlled by the principal curvatures of $M$ in $\mathbb{R}^n$. If the absolute value of those curvatures is bounded by $\kappa$, then $\frac{1}{r^2}\le \kappa$. The intrinsic distance $d(p,q)$ from $p$ to $q$ thus is approximated  by the arclength on the circle, that is, by $r\phi$ for some angle $\phi \in [0,\pi]$  whereas the Euclidean distance is the length of the straight line from $p$ to $q$, that is, in our approximation, the length $2r\sin \frac{\phi}{2}$ of the chord. By Taylor expansion, the difference between the arc and the chord length is of order $r\phi^3$, that is controlled by $\kappa d(p,q)^3$.

This shows that the radial distances in $k$-neighborhoods can indeed be approximated by the Euclidean distances of the embedding space $\mathbb{R}^n$, as long as the points in the $k$-neighborhoods are sufficiently close relative to the extrinsic curvature of the manifold from which they are assumed to be approximately sampled. However, it does not yet show that merging all those fuzzy simplicial sets, or neighborhood graphs, indeed results in a graph whose graph distances approximate geodesic distances when the density of sample points goes to infinity. This, in turn, was shown in \cite[Main Theorem C, Section 5]{Bernstein00}. Hence, we conclude that the procedure is justified for sufficient sampling density, even for non-locally flat underlying manifolds.
As mentioned in Section \ref{sec:Introduction}, 
this justification does not require the uniformity assumption used in \cite{McInnes18} and it leaves further room for refinement by the notion of curvature.

Note that one has some freedom in the choice of the local metrics $d_i$ in step \ref{U1}. In \cite{McInnes18}, they are
\begin{equation}
  \begin{split}
    d_i(x_i,x_{i_j})&=d_i(x_{i_j},x_i)=(d(x_i,x_{i_j})-\rho_i)/\sigma_i\quad
                                 \text{ for }j=1,\dots k\\
    d_i(x,x)&= 0\quad \text{ for all }x\in X_i,\qquad\te{ and }\qquad
    d_i(x_j,x)=\infty \quad \text{ in all other cases.}
\end{split}
\label{eq:localDists}
\end{equation}
where $x_{i_j}$ is the $j$th neighbor of $x_i$ (ordered by distance), $\rho_i:=d(x_i,x_{i_1})$  and $\sigma_i$ is some normalization factor. In UMAP, this normalization factor is inspired by t-SNE (cf.~\cite{vandermaaten08a}), namely given by the solution to the equation $\sum_{j=1}^k \exp\br{-d_i(x_i,x_{i_j})/\sigma_i}=\log_2(k)$ but one can also more simply choose $\sigma_i=d(x_i,x_{i_k})$ to be the distance to the $k$th nearest neighbor, to obtain an approximately uniform data distribution. Another reasonable choice is $\rho_i=0$ and $\sigma_i=1$, which gives the local distances used in Isomap (cf.~\cite{Tenenbaum00}).

In $\bo{UM}$, \eqref{eq:localDists} is not an appropriate metric because the triangle inequality is violated for some distances. Thus, we propose the following new family of local metrics in $\bo{UM}$.
\begin{equation}
  \begin{split}
    d_i(x_i,x_{i_j}) &= d_i(x_{i_j},x_i) = (d(x_i,x_{i_j})-\rho_i)/\sigma_i \quad \text{ for }j=1,\dots k\\
    d_i(x_{j},x_{\ell}) &= f(d_i(x_{i_j},x_i),d_i(x_i,x_{i_\ell}),x_i,x_{i_j},x_{i_\ell}) \quad \text{ for }j,\ell=1,\dots k\\
    d_i(x,x) &= 0 \quad \text{ for all }x\in X_i,\qquad\te{ and }\qquad
    d_i(x_j,x)=\infty \quad \text{ in all other cases.}
\end{split}
\label{eq:localDistsUM}
\end{equation}
Here $f$ can be any function that respects the triangle inequality. For example, canonical choices include $f(d_1,d_2,x,x_1,x_2)=d_1+d_2$ or (as suggested by Michael Freedman) $f(d_1,d_2,x,x_1,x_2)=(d_1+d_2)/\sqrt{2}$. The latter choice takes into account that most directions emanating from a point in a high-dimensional space are approximately orthogonal. If the coordinates $x$ are themselves recorded in a high-dimensional Euclidean space, then $f$ could of course also be the Euclidean distance, appropriately rescaled to account for the effect of $\rho_i$ and $\sigma_i$, i.e.
\begin{equation}
  \begin{split}
    f(d_1,d_2,x,x_1,x_2) = \sqrt{\left(\frac{(x_1-x)}{||x_1-x||} d_1 - \frac{(x_2-x)}{||x_2-x||}d_2\right)^2}.
  \end{split}
\end{equation}
In our algorithmic implementation (cf.~Section \ref{sec:application}), we provide all $3$ options.

Whether one chooses to work in $\bo{EPMet}$ or $\bo{UM}$ depends on what one would like to achieve in a given application. From a geometric point of view, we consider $\bo{UM}$ to be somewhat more natural as uber metric spaces always respect the triangle inequality.

Our next task is to connect geometric considerations to category-theoretical ones by proving an explicit formula for the computation of a colimit in the category of uber metric spaces.

\bp{prop:colimitUM}{
    Let $\bo{M}$ be either $\bo{UM}$ or $\bo{EPMet}$. The colimit of a small diagram $D:\bo{I}\to\bo{M}$ is given by 
    \begin{equation}
        \begin{split}
            \te{colim}(D)=\left(\te{colim}(FD),~d_{\te{colim}}\right)
        \end{split}
    \end{equation}
    where $F$ is the forgetful functor $F:\bo{M}\to\bo{Sets}$, and $\te{colim}(FD)$ is given by
    \begin{equation}
        \begin{split}
            &\te{colim}(FD)\simeq X/\sim,\quad\te{ where }X:=\bigcupdot_{I\in\bo{I}}~FD(I)\\
            &\te{ and }\sim\text{ is generated by }~(~x\sim x' 
            \quad\text{iff}\quad x'=FDu(x)~),
        \end{split}
        \label{eq:colimSetsFD}
    \end{equation}
    (where $u$ is a morphism in the indexing category $\bo{I}$ of the diagram $D$). 

    For the case $\bo{M}=\bo{UM}$,
    $d_{\te{colim}}$ is defined by $d_{\sim}$, defined by 
    \begin{equation}
        \begin{split}
            d_{\sim} ([x], [x']) = \inf(d_X (p_1 , q_1 ) + \cdots + d_X (p_n , q_n )),
        \end{split}
        \label{eq:simMetric2}
    \end{equation}
    where the infimum is taken over all pairs of sequences $(p_1 , \cdots , p_n ),~ (q_1 , \cdots , q_n )$ of
    elements of $X$, such that
    \begin{equation}
        \begin{split}
            p_1 \sim x,\quad q_n \sim x',\quad \text{ and }\quad p_{i+1} \sim q_i ~\te{ for all } ~1 \le i \le n - 1,
        \end{split}
    \end{equation}
    whereas, for the case $\bo{M}=\bo{EPMet}$, $d_{\te{colim}}$ is defined by $d_{\sim}^{\te{EPMet}}$, given by
    \begin{equation}
        \begin{split}
            d_{\sim}^{\te{EPMet}} ([x], [x']):= 
            \begin{cases}
              \infty,&\te{if }\inf_{y\in[x],y'\in [x']}=\infty,\\
              d_{\sim}\qquad\te{(as in \eqref{eq:simMetric2})},&\te{else.}
            \end{cases}
        \end{split}
        \label{eq:simMetricEPMet2}
    \end{equation} 
    For both cases ($\bo{UM}$ or $\bo{EPMet}$), $d_X$ is always defined by
    \begin{equation}
        \begin{split}
            d_{X}(p_i,q_i) := \begin{cases}
                d_{J}(p_i,q_i),&\text{if }p_i,q_i\in FD(J)\\
                \infty,&\text{else.}
            \end{cases}
        \end{split}
        \label{eq:dXmetric}
    \end{equation}
}
\prr{prop:colimitUM}
The colimit in $\bo{UM}$ turns out to be exactly equal to the natural gluing operation that was known in metric geometry for quite some time, cf.~\cite[Chapter 3]{Burago01}, which manifests an example of the idea that category theory can accurately guide the extraction of meaningful formulas and generalizations from abstract nonsense.
Note that the above proposition is what allows us to derive an explicit description of the realization functor in Proposition \ref{prop:realizationOfFuzzySimpSet}. 
Intuitively, this colimit glues metric spaces together in a very natural way:
\begin{enumerate}
    \item One first takes the disjoint union of all spaces in a diagram and on each disjoint space one imposes the original distance,
    \item then one identifies all of the points in the disjoint union that are mapped to each other in the diagram,
    \item and finally imposes the infimum path distance on the resulting space, which means that the new distance between any two points is just the distance along the shortest path that one can take after the identification in step 2.
\end{enumerate}
\yer{example:colimitUM}{
    Suppose our domain category $\bo{I}$ consists of two objects $O_1$ and $O_2$ and there are two morphisms from $O_1$ to $O_2$. Assume further that $D:\bo{I}\to\bo{UM}$ maps $O_1$ to the one-point-space $(\{\star\},d_1)$ (where $d_1$ is defined by $d_1(\star,\star)=0$) and maps $O_2$ to the line interval $([0,1],d_E)$ (where $d_E$ is the Euclidean metric defined by $d_E(x,y)=|x-y|$), and maps the two morphisms in $\bo{I}$ to the two morphisms in $\bo{UM}$ that map $\star$ to $0$ and to $1$ respectively. The colimit of the diagram 
    is then isomorphic to the circle $S^1$ (because $0\in [0,1]$ and $1\in [0,1]$ are now identified) and distances are also measured along the circle (because distances are now distances along shortest paths in the glued space).
}
For more geometric details, we refer the reader to our related book \cite{IsumapBook}.

\section{Fuzzy simplicial sets}
\label{sec:FuzzySimplicialSets}

Our approach to fuzzy sets and fuzzy simplicial sets is based on \cite{barr1986fuzzy} and \cite{Spivak09}.

Let $\bo{I} := [0,1]$ be the closed interval of real numbers, interpreted as a poset category. The unique morphism corresponding to $a \le b$ is denoted by $i_{ab} : a \to b$. Whenever we refer to a subset $B \subset \bo{I}$, we implicitly think of $B$ as a full subcategory.

\gd{def:fuzzySet}{
  A \textit{fuzzy set} $S$ is a presheaf on $\bo{I}$ for which all restriction maps $S(i_{ab}:a\to b):S(b)\to S(a)$ are injections and which satisfies the gluing condition $\cap_{b\in B}S(b)\simeq S(\sup B)$ for all non-empty $B\subset \bo{I}$.
  Their category is denoted by $\bo{Fuz}$.
}
Note that $\bo{I}$ can equivalently be viewed as the category of open sets of a topological space in which the open sets are the intervals $[0,a)$ (for all $a\in [0,1]$) and $[0,1]$, ordered by inclusion.

If we had permitted $B$ to be empty in the above definition, then the gluing condition could be interpreted as a sheaf condition on this topological space.
However, if S were a sheaf, then $S(0)$ would have to be the
singleton set. 
Together with the injectivity condition, this would exclude most fuzzy sets.  
To overcome this, \cite{barr1986fuzzy} augments $\bo{I}$ with an additional bottom element 
$\bot$, and imposes injectivity on all restriction maps except those of the form $S(a)\to S(\bot)$, while \cite{wegmann2026theoryumap} sets $S(0)$ equal to the singleton set and does not require injectivity of the maps $S(a)\to S(0)$.

Our Definition \ref{def:fuzzySet} of fuzzy sets is equivalent to the original definition of~\cite{zadeh1965fuzzy} (see Def.~\ref{def:classicalFuzzySets} and Prop.~\ref{prop:isoClassicalAndSimplicialFuzzySets} below) that we refer to as classical fuzzy sets. Definition~\ref{def:classicalFuzzySets} distinguishes classical fuzzy sets even when they only differ up to elements of strength $0$.
As a result, our category of classical fuzzy sets is equivalent to the category $\bo{Fuz}(L)$ in \cite{barr1986fuzzy} when $L=[0,1]$ but not equivalent to that defined in \cite{Spivak09} and \cite{wegmann2026theoryumap} where $(0,1]$-valued sets are used.

We find our approach convenient because it avoids case distinctions arising in applications of the sheaf approach, where the set $S(\bot)$ or $S(0)$ often has to be treated as special case. (Also \cite{wegmann2026theoryumap} attempts to reduce case distinctions with the notion of bottomless sheaf condition.)
Furthermore, in our application to UMAP and IsUMap, edges with strength $0$ in our fuzzy simplicial sets can directly encode infinite distances. Among other things, this simplifies the definition of recursive T-conorm merge operations (described in Section \ref{sec:TheMergeFunctors}) because the underlying edge sets of the fuzzy simplicial sets involved in our construction then coincide.
Finally, our definition connects with the original work of \cite{zadeh1965fuzzy} and the work that builds up on it.

The following defines a category whose objects correspond to fuzzy sets as originally defined in \cite{zadeh1965fuzzy}.
\gd{def:classicalFuzzySets}{
  A \textit{classical fuzzy set} is a pair $(X,\eta)$ where $X$ is a set and $\eta:X\to [0,1]$ is a function, often called \textit{membership or strength function}. The morphisms  $f:(X,\eta)\to (Y,\xi)$ in the \textit{category of classical fuzzy sets} $\bo{cFuz}$ are functions $f:X\to Y$ that are increasing in the sense that $\xi(f(x))\ge \eta(x)~\forall x\in X$.
}
The next proposition shows the equivalence to our definition of fuzzy sets.
\bp{prop:isoClassicalAndSimplicialFuzzySets}{
  There is an equivalence of categories $C:\bo{Fuz} \to \bo{cFuz}$ that assigns to a fuzzy set $S$ a tuple $C(S):=(S(0),\xi)$, where 
  \begin{equation}
      \begin{split}
          \xi(s\in S(0)):=\te{sup}\{~a~|~s\in \te{im}(S(i_{0a}))~\}.
      \end{split}
      \label{eq:maxStrength}
  \end{equation}
  The inverse $C^{-1}$ assigns, to each pair $(A,\xi)$ the fuzzy set 
  \begin{equation}
      \begin{split}
          S:\bo{I}^{\te{op}}\to\bo{Sets},\quad S(a):=\xi^{-1}([a,1]).
      \end{split}
      \label{eq:SaProperty}
  \end{equation}
}
\pr{See Proof \ref{proof:prop:isoClassicalAndSimplicialFuzzySets}. Note that an analogous Proposition was already given in \cite{Spivak09}.}

Given a category of fuzzy sets, one can define fuzzy \textit{simplicial} sets
in complete analogy to simplicial sets
(i.e.~functors from $\Delta^{\text{op}}$ to $\bo{Set}$),
as a functor from $\Delta^{\text{op}}$ to $\bo{Fuz}$.
\gd{def:simplicialIndexingCategory}{
$\Delta$ denotes the \textit{simplicial indexing category}. Its objects are given
by finite totally ordered sets $[n] := \{0,1,\dots, n\}$ with exactly $n+1$ elements (we follow the standard convention in this context to start counting at $0$) and its morphisms
are order preserving maps ($f:[n]\to [m]$ s.t. $f(a)\ge f(b)$ if $a\ge b$).\footnote{If one thinks of $[n]$ as the poset category $0\to 1\to 2 \to \cdots \to n$, then
order preserving maps become functors and $\Delta$ becomes a subcategory of the category of categories $\bo{Cat}$.}
}
Note that, when no confusion with natural numbers may arise, we write $n$ instead of $[n]$.
It is well-known that each such order preserving $f$ can be written as a composition of face maps $ \delta_i $ and degeneracy maps $ \sigma_i $ defined as follows.
\gd{def:simplicialIdentities}{
    The \textit{face and degeneracy maps} $\delta_i^n:(n-1)\to n$ (for $n\in \{1,2,\cdots\}$) and $\sigma_i^n:(n+1) \to n$ (for $n\in \{0,1,\cdots\}$) are defined,  for $0\le i \le n$, by
    \begin{equation}
        \begin{split}
          \delta_i^n(j)&:= \begin{cases}
                j,&\te{if }j<i,\\
                j+1,&\te{if }j\ge i,
            \end{cases},\qquad 
          \sigma_i^n(j):= \begin{cases}
                j,&\te{if }j\le i,\\
                j-1,&\te{if }j > i.
            \end{cases}
        \end{split}
    \end{equation}
    (Intuitively, $\delta_i^n$ omits the $i$-th index and $\sigma_i^n$ repeats the $i$-th index.)
}
The following lemma is taken from \cite[Chapter VII.5]{MacLane78}.
\bp{prop:AllMorphismsInSimplexCategoryAreCompositionsOfFaceAndDegeneracyMaps}{
    In $\Delta$, any arrow $f:n\to n'$ has a unique representation
    \begin{equation}
        \begin{split}
            f=\delta_{i_1}\circ \cdots \circ \delta_{i_k}\circ \sigma_{j_1}\circ\cdots \circ\sigma_{j_h},
        \end{split}
        \label{eq:canonicalFormInSimplexCategory}
    \end{equation}
    where the ordinal numbers $h$ and $k$ satisfy $n-h+k= n'$, while the strings of the subscripts $i$ and $j$ satisfy $n' >i_1>\cdots>i_k\geq 0$, and $0\le j_1 <\cdots < n-1$.
}
\pr{See \cite[Lemma on p. 177 in Ch. VII.5]{MacLane78}}
In particular, the composite of any two $\delta$'s and $\sigma$'s may be put into the canonical form \eqref{eq:canonicalFormInSimplexCategory}, which implies that there must exist three kinds of binary identities \eqref{eq:identities}, and one can show that those freely generate $\Delta$:
\bp{prop:simplicialIdentitiesDefineSimplicialSet}{
    The face and degeneracy maps fulfill the following \textit{simplicial relations} (whenever the compositions appearing in the relations are well-defined):
    \begin{equation}
        \begin{split}
            \delta_j\circ \delta_i &= \delta_{i}\circ \delta_{j-1},\quad\te{if }i<j,\\
            \sigma_j\circ \sigma_i&=\sigma_{i}\circ \sigma_{j+1},\quad\te{if }i\le j,\\
            \sigma_j\circ \delta_i &= \begin{cases}
                \delta_i\circ \sigma_{j-1},&\te{if }i<j,\\
                \te{id},&\te{if }i\in\{j,j+1\},\\
                \delta_{i-1}\circ \sigma_j,&\te{if }i>j+1.
            \end{cases}
        \end{split}
        \label{eq:identities}
    \end{equation} 
    Furthermore, the simplex category $\Delta$ is freely generated by the arrows $\delta_i^n$ and $\sigma_j^m$ subject to those relations.
}
\pr{See \cite[Proposition 2 on p. 178 in Ch. VII.5]{MacLane78}.\footnote{In \cite{MacLane78}, the first identity in \eqref{eq:identities} is actually written $\delta_i\circ \delta_j = \delta_{j+1}\circ \delta_i$ if $i\le j$, while in some other sources, for example in \cite[p. 22, before Lemma 2.2.]{Gabriel67}, the form $\delta_j\circ \delta_i  = \delta_{i}\circ \delta_{j-1}$ if $i<j$ is used. Those are equivalent. To see this, exchange the left-hand side and the right-hand side of the equation $\delta_i\circ \delta_j = \delta_{j+1}\circ \delta_i$, send $j$ to $j-1$ and use $i\le j-1~\Leftrightarrow~i<j$. 
We use the notation from \cite{Gabriel67} because we consider it somewhat more consistent to have all left-hand side indices of the identities \eqref{eq:identities} in the order $j$ followed by $i$. Furthermore, when considering the corresponding simplicial identities for simplicial sets, \cite{MacLane78} turns around the right-hand side and left-hand side of the first two equations in \eqref{eq:identities} but not of the third. In \cite{Gabriel67}, the third equation in the simplicial identity contains a typo because they write $i=j+1$ instead of $i>j+1$. The identities seem to be correct in \cite{Goerss09} but then they were copied incorrectly from there to the nlab article about simplicial identities, cf.~\cite{nlab:simplicialIdentities}.}}

Equipped with the above definitions, we are now ready to define fuzzy simplicial sets.
\gd{def:fuzzySimSets}{
A \textit{fuzzy simplicial set} is a functor $\Delta^{\text{op}}\to \bo{Fuz}$. Their category, in which the morphisms are the natural transformations, is denoted by $\bo{sFuz}$.}

\bp{prop:simplexFaceStrength}{The strength of a simplex is at most equal to the minimum of the strengths of its faces. All degeneracies of a simplex have the same strength as the simplex.}
\pr{This Proposition was asserted in \cite{Spivak09} but without proof and we provide a proof in \ref{proof:prop:simplexFaceStrength}.}
\gd{def:classicalFuzzySimplicialSets}{
    The category of \textit{classical fuzzy simplicial sets} $\bo{csFuz}$ is the category where objects are functors $S:\Delta^{\te{op}}\to\bo{cFuz}$ and morphisms are natural transformations.
}
\gd{def:csFuz}{
    Define $C_s:\bo{sFuz}\to\bo{csFuz}$ to be the functor given by postcomposing with $ C $. 
}

In more detail we define $ C_s $ by assigning to each fuzzy simplicial set $S\in [\Delta^{\te{op}},\bo{Fuz}]$ a functor $C_s(S)\in [\Delta^{\te{op}},\bo{cFuz}]$, by defining, for each $n$, $C_s(S)(n):=C(S(n))=(S(n)(0),\xi_n)$ according to the definition in Proposition \ref{prop:isoClassicalAndSimplicialFuzzySets} and by assigning, to each morphism $ f :m \to n $ in $\Delta$, the morphism $C_s(S)(f):C(S(n))=(S(n)(0),\xi_n) \to C(S(m))=(S(m)(0),\xi_m)$ that consists of the unique map $S(f)(i_{00}):S(n)(0)\to S(m)(0)$ between the underlying sets of $C(S(n))$ and $C(S(m))$.

\bp{prop:Fsiso}{
    $C_s:\bo{sFuz}\to\bo{csFuz}$ is an equivalence.
}
\prr{prop:Fsiso}

\bp{prop:cFuzCompat}{
    Classical fuzzy simplicial sets $S$ fulfill the property that
    \begin{equation}
        \begin{split}
            \xi_n(s)\le \min\{~\xi_{n-1}(f)~|~f=S(\delta_i)(s),~\delta_i\te{ a face map in }\Delta~\}
        \end{split}
        \label{eq:compatCondCsFuz}
    \end{equation}
    for all $n\in\mathbb{N}$ and $s\in S(n)(0)$ (where, as before, $S(n)=(S(n)(0),\xi_n)$ and $\xi_n$ is defined as in \eqref{eq:maxStrength}).
}
\prr{prop:cFuzCompat}

\section{Truncated fuzzy simplicial sets}
\label{sec:nFuz}

In this section, we provide definitions and propositions in relation to truncated fuzzy simplicial sets that are not only necessary for later Propositions but also interesting in their own right.
\gd{def:truncation}{
  For each $n\in\mathbb{N}$, let $\bo{n}$ denote the full subcategory of $\Delta$ which only contains the objects $0$ to $n$. The corresponding embedding (fully faithful functor) is called $i^n:\bo{n}\to\Delta$ and we also write $i_n:\bo{n}^{\te{op}}\to\Delta^{\te{op}}$.
  We also define $I^n:=i^n\times \te{id}:\bo{n}\times \bo{I}\to\Delta\times \bo{I}$ and $I_n:(\bo{n}\times \bo{I})^{\te{op}}\to(\Delta\times\bo{I})^{\te{op}}$.
}
Using $\bo{n}$ and $I_n$, we define the category $ \bo{nFuz}$ and a truncation functor:
\gd{def:truncation2}{
  Let $ \bo{nFuz} := [\bo{n}^\te{op}, \bo{Fuz}] $ be the \textit{category of ($n$-)truncated fuzzy simplicial sets}.
  The \textit{truncation functor} $\te{tr}_n:\bo{sFuz}\to\bo{nFuz}$ is defined by precomposition:
  \begin{equation}
    \begin{split}
      \te{tr}_n:=-\circ I_n,\quad\te{ i.e.~}\quad
      \te{tr}_n(S)(m,a):=S(I_n(m,a)).
    \end{split}
    \label{eq:truncationFunctorDefinition}
  \end{equation}
  where $ S \in \bo{sFuz} $ is some fuzzy simplicial set and $ (m,a) \in \bo{n}\times \bo{I} $.
}
\bp{prop:leftRightAdjointOfTruncation}{
  Kan extensions equip the truncation functor $\te{tr}_n:\bo{sFuz}\to\bo{nFuz}$ with both a left and a right adjoint, called the skeleton and coskeleton respectively.
  \begin{equation}
    \label{diag:adjunctionTruncationSkeletonCoskeleton}
    \begin{split}
      \bt
      \bo{nFuz}\ar[bend left=60,""{name=A, above}]{rr}{\te{sk}_n}\ar[bend right=60,""{name=C,above}]{rr}{\te{cosk}_n}\&\&\ar[',""{name=B,above}]{ll}{\te{tr}_n}\ar[from=A, to=B, symbol=\dashv]
      \ar[from=B, to=C, symbol=\dashv]\bo{sFuz}
      \et
    \end{split}
  \end{equation}
  The skeleton functor (left Kan extension) is naturally isomorphic to the following coproduct,
  \begin{equation}
    \te{sk}_n(N)(m,a) := \te{Lan}_{I_n}(N)(m,a) \simeq \smashoperator[r]{\coprod_{\substack{f:m \twoheadrightarrow k\\ k\le n}}} \ N^{\te{nd}}(k,a),
    \label{eq:skFormula}
  \end{equation}
  where $N^{\te{nd}}(k,a)$ is the subset of non-degenerate simplices of $N(k,a)$, while the action on morphisms is explained around Diagram and \eqref{diag:nondegenerate_naturality} and \eqref{eq:ActionOnMorphismsOfCoproduct}.

  The coskeleton functor (right Kan extension) is naturally isomorphic to the following Hom-functor,
  \begin{equation}
    \begin{split}
      \te{cosk}_n(N)(m,a) &:= \te{Ran}_{I_n}(N) \simeq 
      \te{Hom}_{\bo{nFuz}}(\te{tr}_n(\Delta^m_{a}),N),
    \end{split}
    \label{eq:coskFormula}
  \end{equation}
  where $\Delta^m_a := \bo{y}(m,a)$ and $\bo{y}:\Delta\times\bo{I}\to \bo{PSh}(\Delta\times\bo{I})$ is the Yoneda embedding.
}
\prr{prop:leftRightAdjointOfTruncation}

We remark that our proof of the above proposition in fact shows that the adjunctions displayed in Diagram \eqref{diag:adjunctionTruncationSkeletonCoskeleton} (as well as the explicit descriptions in \eqref{eq:skFormula} and \eqref{eq:coskFormula}) exist more generally between the presheaf categories $[(\bo{n}\times \bo{I})^{\te{op}},\bo{Set}]$ and $[(\bo{Delta}\times \bo{I})^{\te{op}},\bo{Set}]$, while additionally restricting appropriately to adjunctions between $\bo{nFuz}$ and $\bo{sFuz}$.

By Prop.~\ref{prop:adjunctionsPreserveLimits}, the adjunctions in Prop.~\ref{prop:leftRightAdjointOfTruncation} give rise to the following corollary:
\bc{cor:trPreservesLimsColims}{
  $\te{tr}_n $ preserves both limits and colimits, i.e.~for any small diagram $D:\bo{A}\to\bo{sFuz}$, we have 
  \begin{equation}
    \begin{split}
      \te{tr}_n(\te{lim}(D))\simeq\te{lim}(\te{tr}_n(D))\quad\te{ and }\quad \te{tr}_n(\te{colim}(D))\simeq\te{colim}(\te{tr}_n(D)).
    \end{split}
  \end{equation}
}
Furthermore, for $m\le n$, $\te{sk}_n$ and $\te{cosk}_n$ act like the identity:
\bp{cor:trcoskId}{
  The skeleton and coskeleton functor are fully faithful, or equivalently
  \begin{equation}
    \begin{split}
        \te{id} \cong \te{tr}_n\circ \te{sk}_n \qquad\qquad \te{tr}_n\circ \te{cosk}_n \cong \te{id}
    \end{split}
    \label{eq:idTr}
  \end{equation}
  where the isomorphism are exhibited by the unit and the counit respectively. In particular
 \begin{equation}
      \begin{split}
          \te{tr}_m(\te{sk}_n(N))(k,a)\cong N(k,a)\cong \te{tr}_m(\te{cosk}_n(N))(k,a)
      \end{split}
      \label{eq:idTrmn}
  \end{equation}
  for all $ m \le n $.
}
\prr{cor:trcoskId}

Similarly to the functor $C_s:\bo{sFuz}\to\bo{csFuz}$, there is also, for each $n\in\mathbb{N}$, a category $\bo{cnFuz}$ and a functor $C_n:\bo{nFuz}\to\bo{cnFuz}$:
\gd{def:cnFuz}{
  Define
  \begin{equation}
      \begin{split}
          \bo{cnFuz}:=[\bo{n}^{\text{op}},\bo{cFuz}].
      \end{split}
  \end{equation}
  Furthermore, let $C_n:\bo{nFuz}\to\bo{cnFuz}$ the functor given by postcomposing with $ C $ from proposition \eqref{prop:isoClassicalAndSimplicialFuzzySets}. 
}

In more detail $C_n:\bo{nFuz}\to\bo{cnFuz}$ is defined by assigning to each $S\in \bo{nFuz}$ a functor $C_s(S)\in \bo{cnFuz}$, which in turn is defined, for each $m$, by $C_n(S)(m):=C(S(m))=(S(m)(0),\xi_m)$ according to the definition in Proposition \ref{prop:isoClassicalAndSimplicialFuzzySets} and by assigning, to each morphism $ f :l \to m$ in $\bo{n}$, the morphism $C_n(S)(f):C(S(m))=(S(m)(0),\xi_m) \to C(S(l))=(S(l)(0),\xi_l)$ that consists of the map $S(f)(i_{00}):S(n)(0)\to S(m)(0)$ between the underlying sets of $C(S(n))$ and $C(S(m))$.

\bp{prop:Cniso}{
  $ C_n : \bo{nFuz} \to \bo{cnFuz} $ is an equivalence of categories.
}
\pr{
  Same argument as in the proof of Proposition \ref{prop:Fsiso}.
}
\yr{rem:weightedgraphs}{
  Note that $\bo{c1Fuz}$ is the category of directed weighted graphs with weights in $[0,1]$, in which vertices also carry weights and edge weights are at most equal to the minimum of the weights of their vertices. 
  However, one can single out those $S\in\bo{c1Fuz}$, where vertices all have weight $1$, and then only the weights on the edges matter, giving rise to ordinary directed weighted graphs. 
}
Using that $C_n$ and $C_s$ are equivalences, we can define 
\begin{align}
        \te{csk}_n: & =C_s\circ \te{sk}_n\circ C_n^{-1},\\
        \te{ctr}_n: & =C_n\circ \te{tr}_n\circ C_s^{-1},\\
        \te{ccosk}_n: & =C_s\circ \te{cosk}_n\circ C_n^{-1}
    \label{eq:ctrnEtc}
\end{align}
We finish the section with the following corollary.
\bc{prop:classicalVersionOfAdjunctions}{
  We have $\te{csk}_n \dashv \te{ctr}_n \dashv \te{ccosk}_n$, i.e.~$\te{csk}_n$ is left adjoint to $\te{ctr}_n$ and $\te{ctr}_n$ is left adjoint to $\te{ccosk}_n$.
}

\section{Metric realizations}
\label{sec:Adjunction}

As mentioned in Section \ref{sec:Introduction}, \cite{Spivak09} introduced an adjunction 
between $\bo{UM}$ and $\bo{sFuz}$. This adjunction arises as the restriction to suitable subcategories of an abstract construction known as the nerve-realization adjunction. In \cite{McInnes18}, a similar adjunction was provided between $\bo{EPMet}$ and $\bo{sFuz}$ that emphasizes spaces where the underlying sets are finite.
In the following, we provide all relevant proofs, filling in some previously missing details, provide a slight generalization of Spivak's construction and use our Proposition \ref{prop:colimitUM} to derive more explicit formulas for the metric realizations. Furthermore, we prove that (a non-finite version of) the adjunction defined in \cite{McInnes18} can actually be used to prove that $\bo{UM}$ and $\bo{EPMet}$ can be embedded into $\bo{sFuz}$. 
More precisely, $ \bo{UM} $ and $\bo{EPMet}$ turn out to be reflective subcategories of $ \bo{sFuz} $ and  the adjunction restricts to an equivalence on the essential image of the nerve functor ($\text{Sing}$). This result is one of the reasons why we believe that fuzzy simplicial sets are the right combinatorial structure for encoding metric relations as remarked in \ref{A1}: They can represent any metric space exactly but there are also fuzzy simplicial sets that are less constrained. We also characterize the essential image of the nerve functor $\text{Sing}$ in Prop.~\ref{prop:SingSimpSetProperties} and finally outline some relationships to persistent homology in Sec.~\ref{sec:relationToPersHomology}.

Let $\bo{y}(\bo{C})$ be the subcategory of the presheaf category $\bo{PSh}(\bo{C})=[\bo{C}^{\te{op}},\bo{Set}]$ obtained by the Yoneda embedding.
To begin, we recall the following formulation of the well-known nerve-realization theorem:
\bp{prop:KanExtension}{
  For any category $\bo{C}$, suppose that $\te{Re}:\bo{y}(\bo{C})\to\bo{M}$ is any functor and that $\bo{M}$ has small colimits. Then $\te{Re}$ can be extended to the following functor:
    \begin{equation}
        \begin{split}
            \te{Re}:&~\bo{PSh}(\bo{C})\to\bo{M},\\
            \te{Re}(S):&=\te{colim}(D_S)\\\quad\te{where}\quad D_S&=\te{Re} \circ \bo{y}\circ P_S :\bo{El}(S)\to\bo{M},
        \end{split}
        \label{eq:ReExtended}
    \end{equation}
    where the category of elements $\bo{El}(S)$ consists of objects which are
    pairs $(A,x)$ with $A\in \bo{C}$ and $x\in S(A)$ and $P_S:\bo{El}(S)\to \bo{C}$ is the projection functor from the category of elements to $\bo{C}$.

    Furthermore, the following functor is right-adjoint to $\te{Re}$:
    \begin{equation}
        \begin{split}
            \te{Sing}:\bo{M}&\to\bo{PSh}(\bo{C}),\\
                      Y&\mapsto \te{Sing}(Y):\bo{C}^{\te{op}}\to\bo{Sets},\\
                      &\qquad\qquad\quad\quad~C\mapsto \te{Hom}_{\bo{M}}(\te{Re}(\bo{y}(C)),~Y).
        \end{split}
    \end{equation}
}
\pr{
  We use Prop.~\ref{prop:kanExtension} about Kan extensions to provide a proof in App.~\ref{proof:prop:KanExtension}. For an alternative proof, see also \cite[Proposition 9.16.]{awodey10}.
}

Now consider $\bo{C}=\Delta\times \bo{I}$ and $\bo{M}=\bo{UM}$ or $\bo{M}=\bo{EPMet}$.
As we have shown in Proposition \ref{prop:colimitUM}, both $\bo{UM}$ and $\bo{EPMet}$ have small colimits.
Therefore, whenever we can define a functor $\te{Re}$ from $\bo{y}(\Delta\times \bo{I})$ to $\bo{UM}$ or $\bo{EPMet}$, then the above proposition gives us an adjunction between $\bo{PSh}(\Delta\times \bo{I})$ and the corresponding metric category. For a given $(n,a)\in \Delta\times\bo{I}$, we subsequently use the notation $\Delta^n_a:=\bo{y}(n,a)$.

For every natural number $p \ge 1$ we now define a functor $\text{Re}^p : \bo{y}(\Delta \times \bo{I}) \to \bo{UM} $. Let $\Vert - \Vert_p$ denote the $\ell^p$-norm, that is $\parallel x \parallel_p := (\sum_{i=1}^{n+1}x_i^p)^{1/p}$. Given a pair $(n,a) \in \Delta \times \bo{I} $ we define the functors on objects by
\begin{equation}
  \label{eq:spivakReObjects}
    \text{Re}^p(\Delta_{a}^n) := \left\{ x \in [0,1]^{n+1}~\bigg|~
    \parallel x \parallel_p =
    1 \right\}
\end{equation}
as sets. We endow these sets with the structure of a metric space $(\text{Re}^p(\Delta_{a}^n), d_a)$: for $a \not= 0$ we take $d_a(x, y) = -\log a\parallel x - y \parallel_p $, while for $a = 0 $ we take $d_0$ to be the metric defined by $d_0(x,y)=\infty$ whenever $x \not = y$ (and otherwise $d_0(x,y)=0$).

To define the action on morphisms, we introduce the notation $\parallel x_I \parallel_p ~:= (\sum_{i\in I}x_i^p)^{1/p}$, where $I\subset \{0,\cdots,n\}$.
Then, given morphisms $\sigma : n \to m$ in $\Delta$ and $i_{ab} : a \to b$ in $\bo{I}$, we set
\begin{equation} 
  \label{eq:spivakRe-a->b}
  \text{Re}^p(\bo{y}(\sigma,i_{ab}))(x) := 
  \left( \parallel x_{\sigma^{-1}(0)} \parallel_p ,~\cdots~,\parallel x_{\sigma^{-1}(m)} \parallel_p\right).
\end{equation}
For $p=1$ this reduces to Spivak's definition. Note that the case $p=1$ only makes sense when the $\ell^1$-norm is used, as was pointed out by \cite{wegmann2026theoryumap}. However, it is possible to generalize to higher $p$ if one imposes the $\ell^p$-norm in \eqref{eq:spivakReObjects}, as we show below.

\bp{prop:spivakReNonExpansive}{
  For all $p\in\mathbb{N}$, $\text{Re}^p$ defined in \eqref{eq:spivakReObjects} and \eqref{eq:spivakRe-a->b} is indeed a functor. (In particular, \eqref{eq:spivakRe-a->b} defines a non-expansive map.)
}
\prr{prop:spivakReNonExpansive}

\yr{rem:ReDifferentMetricsSemiSemplicial}{
  What prevents us from equipping the standard simplex $\text{Re}^1(\Delta_{a}^n)$ with the standard Euclidean $\ell^2$ metric is that the action of $\te{Re}$ defined in \eqref{eq:spivakRe-a->b} on degeneracy maps then becomes non-expansive. However, instead of imposing the $\ell^1$ metric on $\text{Re}^1(\Delta_{a}^n)$ to save non-expansiveness, one could alternatively explore the possibility to use a semi-semplicial indexing category $\Delta_+$ (in which one only allows for face maps and their compositions) to build up a theory of semi-simplicial fuzzy sets, where \eqref{eq:spivakRe-a->b} remains non-expansive even when combining $\text{Re}^1(\Delta_{a}^n)$ with $\ell^2$.
}

Given any natural number $p \ge 1$ and $\text{Re}^p$ as defined above, we can use Prop.~\ref{prop:KanExtension} to extend it to a functor $\text{Re}^p \in \bo{PSh}(\Delta\times \bo{I})$, and to obtain an adjunction
\begin{equation}
  \label{eq:presheafAdjunction}
    \begin{tikzcd}
    \bo{PSh}(\Delta\times \bo{I}) \arrow[r, bend left=25, "\text{Re}^p", ""{name=A, below}] & 
    \bo{UM} \arrow[l, bend left=25, "\text{Sing}^p", ""{name=B, above}]
    \arrow[phantom, from=A, to=B, "\dashv" rotate=-90],
    \end{tikzcd}
\end{equation}
where $\text{Sing}^p(Y)(n,a) = \te{Hom}_{\bo{M}}(\te{Re}^p(\Delta^n_a),~Y)$. However, to prove that this adjunction factors through $\bo{sFuz}$ is non-trivial because it requires to show that, for every $Y \in\bo{UM}$ and every $n\in \Delta$, the presheaf $\text{Sing}^p(Y)(n,-):\bo{I}^{\te{op}}\to\bo{Set}$ is indeed a fuzzy set, that is it fulfills the injectivity and gluing conditions of a fuzzy set. This is what the next proposition is about.
\bp{prop:singFactorsThroughSfuz}{
   For every $p\in\mathbb{N}$ and $Y\in\bo{UM}$, $\text{Sing}^p(Y)$ is a fuzzy simplicial set, implying that the adjunction \eqref{eq:presheafAdjunction} can be restricted to an adjunction between $\bo{sFuz}$ and $\bo{UM}$.
}
\prr{prop:singFactorsThroughSfuz}

Note that $\te{Re}^p$ is still a well-defined functor if the codomain category is $\bo{EPMet}$ instead of $\bo{UM}$. Therefore, there is also a corresponding adjunction between $\bo{sFuz}$ and $\bo{EPMet}$.

Next, we answer \ref{q2}, by showing that the functor $\te{FinReal}: \bo{y}(\Delta\times \bo{I})\to\bo{FinEPMet}$ defined in \cite{McInnes18} is still a functor when reconsidering it with codomain $\bo{UM}$. We call the resulting functor $\te{Re}^{\te{skeleton}}:\bo{y}(\Delta\times \bo{I})\to\bo{UM}$ because the Sing-functor adjoint to $\te{Re}^{\te{skeleton}}$ turns out to be the functor that maps subsets of points of a metric space to their simplicial $1$-(co)skeleton, as we show later in Proposition \ref{prop:simplifySpivaksIdea}.
\bp{prop:umapReIsFunWithCodUM}{
    $\te{Re}^{\te{skeleton}}:~\bo{y}(\Delta\times \bo{I})\to \bo{UM}$ defined by 
    \begin{equation}
        \begin{split}
            \te{Re}^{\te{skeleton}}:~&\bo{y}(\Delta\times \bo{I})\to \bo{UM},
            \quad \Delta_{a}^n \mapsto (\{x_0,\cdots,x_n\},d_a),\\
            &\te{ where }\quad d_a(x_i,x_j) := \begin{cases}
                -\log(a),&\te{if }i\ne j,\\
                0,&\te{else.}
            \end{cases}
        \end{split}
        \label{eq:reDeltaUMAP}
    \end{equation}
    and
    \begin{equation}
        \begin{split}
            \te{Re}^{\te{skeleton}}(\sigma,i_{ab}):~&\te{Re}(\Delta_{a}^n)\to \te{Re}(\Delta_{b}^m),\quad x_i\mapsto x_{\sigma(i)}
        \end{split}
        \label{eq:reDeltaMorphisms}
    \end{equation}
    is a functor. (In particular, morphisms defined by \eqref{eq:reDeltaUMAP} are non-expansive.)
    Similarly, $\te{Re}^{\te{skeleton}}:~\bo{y}(\Delta\times \bo{I})\to \bo{EPMet}$, defined in the same way, is also a functor.
}
\prr{prop:umapReIsFunWithCodUM}

In analogy to $\te{Re}^p$, we can now again use Prop.~\ref{prop:KanExtension} to extend $\te{Re}^{\te{skeleton}}$ to a functor from $\bo{PSh}(\Delta\times\bo{I})$ to $\bo{UM}$ or $\bo{EPMet}$. As before, we then obtain an adjunction
\begin{equation}
  \label{eq:presheafAdjunctionSkeleton}
    \begin{tikzcd}
    \bo{PSh}(\Delta\times \bo{I}) \arrow[r, bend left=25, "\text{Re}^{\te{skeleton}}", ""{name=A, below}] & 
    \bo{UM} \arrow[l, bend left=25, "\text{Sing}^{\te{skeleton}}", ""{name=B, above}]
    \arrow[phantom, from=A, to=B, "\dashv" rotate=-90],
    \end{tikzcd}
\end{equation}
The following proposition establishes that also this adjunction factors through $\bo{sFuz}$.

\bp{prop:skeletonSingFactorsThroughSfuz}{
    For every $Y\in\bo{M}$, $\text{Sing}^{\te{skeleton}}(Y)$ is a fuzzy simplicial set, implying that the adjunction \eqref{eq:presheafAdjunctionSkeleton} can be restricted to an adjunction between $\bo{sFuz}$ and $\bo{M}$, where $\bo{M}=\bo{UM}$ or $\bo{EPMet}$.
}
\prr{prop:skeletonSingFactorsThroughSfuz}

Now that we derived general formulas for the Re and Sing functors, that specialize to those provided in \cite{Spivak09} and \cite{McInnes18}, our next task is to attempt to derive more explicit formulas for their computation. Generally, once the $\te{Re}$ functor is given, the Sing functor already has a comparatively explicit description in terms of a Hom-set but the Re functor is defined by a colimit that can be non-trivial to compute. First of all, it requires the derivation of a more explicit formula for the colimit in $\bo{UM}$ and $\bo{EPMet}$. Fortunately, we already derived such a formula in Proposition \ref{prop:colimitUM}.
Hence, the only thing that remains to be done is to analyze how the structure of the fuzzy simplicial set $S$ is related to the diagram $D_S$ and its colimit expressed in \eqref{eq:ReExtended}.

To analyze the diagram $D_S$, we make use of the fact that the only relevant maps in $\Delta$, that we need to look at, are the identity, face and degeneracy maps (cf.~Proposition \ref{prop:AllMorphismsInSimplexCategoryAreCompositionsOfFaceAndDegeneracyMaps}). Hence the relevant maps in $\Delta\times \bo{I}$ are of the form $(\te{id},i_{ab})$, $(\sigma,i_{ab})$ and $(\delta,i_{ab})$. 

We next look at what happens when considering such maps in the category $\bo{El}(S)$ and \textit{by abuse of notation simply refer to them in the same way}. We start with $(\te{id},i_{ab})$. Iterating maps of this form yields diagrams of the form 
\begin{equation}
    \begin{split}
      \te{Identity maps id}\qquad\bt\cdots 
        \ar[yshift=0.7em]{r}{(\te{id},i_{ca})}\& \ar[yshift=-0.85em]{l}{S(\te{id},i_{ca})} 
        \begin{pmatrix}
            m,a\\
            s
        \end{pmatrix}
          \ar[yshift=0.7em]{r}{(\te{id},i_{ab})}\& \ar[yshift=-0.85em]{l}{S(\te{id},i_{ab})} 
          \begin{pmatrix}
              m,b\\
              s
          \end{pmatrix} 
        \et
    \end{split}
    \label{eq:CollapseOfDiagramsSequence}
\end{equation}
whose metric realization often turns out to have an especially simple form.
Furthermore, diagrams such as \eqref{eq:CollapseOfDiagramsSequence} appear as building blocks for diagrams of degeneracy and face maps. Next, Proposition \ref{prop:simplexFaceStrength} implies that degeneracy and face maps give rise to diagrams that look as follows in the category of elements:
\begin{equation}
    \begin{split}
      \te{Degeneracy maps $\sigma$}\qquad
        \bt 
        \cdots 
        \ar[yshift=0.7em]{r}{(\te{id},i_{ca})}\& \ar[yshift=-0.85em]{l}{S(\te{id},i_{ca})} 
        \begin{pmatrix}
            m+1,a\\
            s
        \end{pmatrix}
        \ar[xshift=-0.5em,']{d}{(\sigma,i_{aa})}
          \ar[yshift=0.7em]{r}{(\te{id},i_{ab})}\& \ar[yshift=-0.85em]{l}{S(\te{id},i_{ab})} 
          \begin{pmatrix}
              m+1,b\\
              s
          \end{pmatrix}
          \ar[xshift=-0.5em,']{d}{(\sigma,i_{bb})} \\
        \cdots 
        \ar[yshift=0.7em]{r}{(\te{id},i_{ca})}\& \ar[yshift=-0.85em]{l}{S(\te{id},i_{ca})} 
        \begin{pmatrix}
            m,a\\
            s'
        \end{pmatrix}
        \ar[xshift=0.5em,']{u}{S(\sigma,i_{aa})}
          \ar[yshift=0.7em]{r}{(\te{id},i_{ab})}\& \ar[yshift=-0.85em]{l}{S(\te{id},i_{ab})} 
          \begin{pmatrix}
              m,b\\
              s'
          \end{pmatrix} 
          \ar[xshift=0.5em,']{u}{S(\sigma,i_{bb})}
        \et
    \end{split}
    \label{eq:degeneracyDiagram}
\end{equation}
Here, one can observe that the maximum strength of both, the degeneracy $s$ and the simplex $s'$ is equal (Proposition \ref{prop:simplexFaceStrength}).
\begin{equation}
    \begin{split}
      \te{Face maps $\delta$}\qquad
        \bt 
        \cdots 
        \ar[yshift=0.7em]{r}{(\te{id},i_{ca})}\& \ar[yshift=-0.85em]{l}{S(\te{id},i_{ca})} 
        \begin{pmatrix}
            m,a\\
            s
        \end{pmatrix}
        \ar[xshift=-0.5em,']{d}{(\delta,i_{aa})}
          \ar[yshift=0.7em]{r}{(\te{id},i_{ab})}\& \ar[yshift=-0.85em]{l}{S(\te{id},i_{ab})} 
          \begin{pmatrix}
              m,b\\
              s
          \end{pmatrix}
          \ar[xshift=-0.5em,']{d}{(\delta,i_{bb})}
          \ar[yshift=0.7em]{r}{(\te{id},i_{bd})}\& \ar[yshift=-0.85em]{l}{S(\te{id},i_{bd})} 
          \begin{pmatrix}
              m,d\\
              s
          \end{pmatrix} \\
        \cdots 
        \ar[yshift=0.7em]{r}{(\te{id},i_{ca})}\& \ar[yshift=-0.85em]{l}{S(\te{id},i_{ca})} 
        \begin{pmatrix}
            m+1,a\\
            s'
        \end{pmatrix}
        \ar[xshift=0.5em,']{u}{S(\delta,i_{aa})}
          \ar[yshift=0.7em]{r}{(\te{id},i_{ab})}\& \ar[yshift=-0.85em]{l}{S(\te{id},i_{ab})} 
          \begin{pmatrix}
              m+1,b\\
              s'
          \end{pmatrix} 
          \ar[xshift=0.5em,']{u}{S(\delta,i_{bb})} \&
        \et
    \end{split}
    \label{eq:faceMapDiagram}
\end{equation}
Here, the upper row might go on to the right because the face might have a higher maximum strength than the simplex (Proposition \ref{prop:simplexFaceStrength}).

For both adjunctions, determined by \eqref{eq:spivakRe-a->b} and \eqref{eq:reDeltaUMAP}, respectively, it is relatively easy to derive explicit descriptions of $(\te{id},i_{ab})$ and $(\sigma,i_{ab})$. Namely, we can derive the following result.

\bp{prop:metricRealizationOfIdentity}{
    For both, the adjunction determined by \eqref{eq:spivakRe-a->b} and the adjunction determined by \eqref{eq:reDeltaUMAP}, the metric realization of diagrams of the form \eqref{eq:CollapseOfDiagramsSequence}
    equals $\te{Re}(\Delta^m_{b})$. In other words, a given fuzzy $m$-simplex is always realized with its maximally available strength. Furthermore, the metric realization of diagrams of the form \eqref{eq:degeneracyDiagram} is also isomorphic to $\te{Re}(\Delta^m_{b})$.
}
\prr{prop:metricRealizationOfIdentity}

However, for Spivaks adjunction \eqref{eq:spivakRe-a->b}, an explicit formula for face maps, that goes beyond a characterization in terms of equivalence classes, is hard to derive. For that reason, we restrict ourselves in the following to a derivation of an explicit formula for the skeleton adjunction, i.e.~for the rest of this section, we make the assumption 
\begin{equation}
  \begin{split}
      \te{Re} = \te{Re}^{\te{skeleton}} : \bo{y}(\Delta\times \bo{I})\to \bo{M},\te{ where }\bo{M}\te{ is either }\bo{UM}\te{ or }\bo{EPMet},
  \end{split}
  \label{eq:ReIsReSkeleton}
\end{equation}
where $\te{Re}^{\te{skeleton}}$ is defined in Proposition \ref{prop:umapReIsFunWithCodUM}. 
In contrast to other possible choices of $\te{Re}$, the colimit of face map diagrams simplifies when assuming \eqref{eq:ReIsReSkeleton} and we obtain the following result.

\bp{prop:faceMap}{
  Assuming \eqref{eq:ReIsReSkeleton}, the metric realization of a face map diagram \eqref{eq:faceMapDiagram}, with face map $\delta:m\to m+1$, results in a metric space of the following form:
  \begin{equation}
    \begin{split}
      (X,d_\delta)&=(\{x_0,\cdots,x_{m+1}\},d_\delta),~\te{ with}\\
      d_\delta(x_i,x_j)&= \begin{cases}
        0,&\te{if }i=j,\\
        -\log(d),&\te{if }i,j\in\te{im}(\delta)\te{ but }i\ne j,\\
        -\log(b),&\te{if }i\notin \te{im}(\delta)\te{ or }j\notin \te{im}(\delta)\te{ but }i\ne j. 
      \end{cases}
    \end{split}
    \label{eq:facemaprealization}
  \end{equation}
}
\prr{prop:faceMap}

This nice result is what allows us to give an explicit characterization of the metric realization functor in the next definition and proposition.

\gd{def:Rec2}{
  Suppose we are given an object $S$ in $\bo{c1Fuz}$ (defined in section \ref{sec:nFuz}). Recall that $S$ is described by $S(0)=(S_0,\xi_0)$ and $S(1)=(S_1,\xi_1)$ and the two face maps $\bt[large] S_0\&\ar[']{l}{S(\delta_0)}S_1\ar{r}{S(\delta_1)} \& S_0\et$, where $\xi_0:S_0\to [0,1]$ and $\xi_1:S_1\to [0,1]$ are some functions such that condition \eqref{eq:compatCondCsFuz} is fulfilled (for $n\in\{0,1\}$). (There is also a trivial degeneracy map $S(\sigma):S_0\to S_1$.) 

  Define the realization functor $\te{Re}_{c1}:\bo{c1Fuz}\to\bo{UM}$ by $\te{Re}_{c1}(S):=(S_0,d)$, where 
  \begin{equation}
      \begin{split}
          d(x,y)&:=\inf_{x=x_1,\cdots,x_n=y} \sum_{i=1}^{n-1} d_{\min}(x_i,x_{i+1}),
      \end{split}
  \label{eq:realizationExplicit}
  \end{equation}
  where 
  \begin{equation}
      \begin{split}
          d_{\te{min}}(x_1,x_2)&:=\min\{-\log(\xi_1(s))~|~\exists~s\in S_1~:~c_1(s)\te{ or }c_2(s)~\},\\
          \te{where }&\qquad c_1(s):=(x_1=S(\delta_0)(s)\te{ and }x_2=S(\delta_1)(s)),\\
          \te{and }&\qquad c_2(s):=(x_1=S(\delta_1)(s)\te{ and }x_2=S(\delta_0)(s)).
      \end{split}
      \label{eq:nSimplexRealizationExplicit}
  \end{equation}
  (If no $s$ exists such that $c_1(s)$ or $c_2(s)$, then $d_{\te{min}}(x_0,x_1):=\infty$.) 

  Furthermore, define the realization functor $\te{Re}_{c1}:\bo{c1Fuz}\to\bo{EPMet}$ by $\te{Re}_{c1}(S):=(S_0,d_{\te{EPMet}})$, where 
  \begin{equation}
      \begin{split}
          d_{\te{EPMet}}(x,y)&:=
          \begin{cases}
            \infty,&\te{if }\xi_1(s)=0~\forall s\in S_1~:~c_1(s)\te{ or }c_2(s),\\
            d(x,y)\qquad\te{(as in \eqref{eq:realizationExplicit})},&\te{else.}
          \end{cases}
      \end{split}
  \label{eq:realizationExplicitEPMet}
  \end{equation}
}
\bp{prop:realizationOfFuzzySimpSet}{
  Assuming \eqref{eq:ReIsReSkeleton}, we have
  \begin{equation}
      \begin{split}
          \te{Re}(S)\simeq \te{Re}_{c1}(\te{ctr}_1(C_s(S)))\simeq \te{Re}_{c1}(C_{1}(\te{tr}_1(S))),
      \end{split}
  \end{equation}
  i.e.~the following diagram commutes:
  \begin{equation}
    \begin{split}
        \bt 
        \bo{csFuz} \ar{d}{\te{ctr}_1} \& \bo{sFuz} \ar[']{l}{C_s}\ar{r}{\te{tr}_1} \ar{d}{\te{Re}} \& \bo{1Fuz} \ar{d}{C_1} \\
        \bo{c1Fuz} \ar{r}{\te{Re}_{c1}} \&\bo{M} \& \bo{c1Fuz} \ar[']{l}{\te{Re}_{c1}}
        \et
    \end{split}
\end{equation}
}
\prr{prop:realizationOfFuzzySimpSet}

Note that an explicit description of $\te{ctr}_1(C_s(S))$ is as follows. According to Definition \ref{def:csFuz} and \eqref{eq:ctrnEtc}, $\te{ctr}_1(C_s(S))$ is the classical fuzzy graph $\bt[huge] S(0,0)\&\ar[']{l}{S(\delta_0,i_{00})}S(1,0)\ar{r}{S(\delta_1,i_{00})} \& S(0,0)\et$ of $S$ with weight functions $\xi_0:S(0,0)\to [0,1]$ and $\xi_1:S(1,0)\to [0,1]$, defined as in \eqref{eq:maxStrength}.
Note further that $\te{ctr}_1\circ C_s \simeq C_1 \circ \te{tr}_1$ simply follows from propositions \ref{prop:Fsiso} and \ref{prop:Cniso} but the main point of the proposition is that $\te{Re}$ factors through $\bo{c1Fuz}$, i.e.~the only information needed for the metric realization of the simple skeleton adjunction is contained in the fuzzy 1-simplices. 

Intuitively, one can think of the metric realization as a space of points, where the distance between any two points is either simply the minimum distance among all the face realizations in which those two points are contained, or, if any of their common faces is longer than the length of a path along several shorter faces, then one obtains their distance as the geodesic distance by hopping along the shortest path connecting them. Since the simplices are not filled as in the case of Spivak's realization functor, one does not have to worry about any other points that are not themselves the vertex points of the simplices. One can also think of those vertex points as the zero-skeleton of a simplicial complex.

Now that we have completed an explicit description of $\te{Re}^{\te{skeleton}}$, we would like to do the same for $\te{Sing}^{\te{skeleton}}$. Fortunately, this is a much simpler task, achieved by the following proposition.

\bp{prop:SingFun}{
  Assuming \eqref{eq:ReIsReSkeleton}, $\te{Sing}(X,d)(n,a)$ is in bijective correspondence with
  \begin{align}
    [X,d]^n_a := \{ (r_0, \dots, r_n) \in X^{n+1} \, | \, d(r_i, r_j) \le - \log(a) \ \ \forall i \neq j \}.
  \end{align}
  where the latter sets can in fact be assembled into a functor isomorphic to $ \te{Sing} $. The strength of a simplex $ s \in \te{Sing}(X,d)(n,a) $, seen as a tuple $ (r_0, \dots, r_n) \in X^{n+1} $, is
  \begin{align}
    \te{str}(s) = \min_{0 \le i,j \le n}\exp(-d(r_{i_j},r_{i_k})) \ge a
  \end{align}
  In particular $(r_0, \dots, r_n) $ is in $ \te{Sing}(X,d)(n,a) $ iff all pairs $ (r_i,r_j) $ are in $ \te{Sing}(X,d)(1,a) $.
}
\prr{prop:SingFun}

We now continue our analysis of the relation between uber metric spaces and fuzzy simplicial sets. The next proposition says that $ \te{Sing} $ exhibits $\bo{UM}$ and $\bo{EPMet}$ as full subcategories of $\bo{sFuz}$.

\bp{prop:UMisInSFuz2}{
Assuming \eqref{eq:ReIsReSkeleton}, $ \te{Sing} : \bo{M} \to \bo{sFuz} $ is fully faithful.
}
\prr{prop:UMisInSFuz2}

\bc{prop:UMisInSfuz}{
  Assuming \eqref{eq:ReIsReSkeleton}, the counit $\te{Re}\circ\te{Sing} \stackrel{\cong}{\longrightarrow} \te{id}_{\bo{M}}$ is an isomorphism.
}
\prr{prop:UMisInSfuz}

To show how fuzzy simplicial sets generalize metric spaces, we next characterize the essential image of Sing inside of $\bo{sFuz}$.
For this consider the diagram category $ \mathcal D_n $ where $ n \ge 2 $, defined as follows:
\begin{itemize}
  \item $ \te{Ob}(\mathcal D_n) := \{i\}_{0 \le i \le n} \sqcup \{ (i,j) \}_{0 \le i < j \le n} $
  \item The only non-trivial morphisms are of the form $ k \to (i,j) $ if either $ k = i $, or $ k = j $.
\end{itemize}
Define a functor $ F_n : \mathcal D_n \to \Delta $, which on objects is given by $ i \mapsto [0] $ and $ (i,j) \mapsto [1] $ while the morphisms $ i \to (i,j) $ and $ i \to (k,i) $ are mapped to $ ,\delta_1 $ and $ \delta_0 $ respectively. Precomposing this with a fuzzy simplicial set $ S$ at a fixed $ a \in [0,1] $ can then be represented pictorially as
\begin{center}
  \begin{tikzcd}[ampersand replacement=\&]
    {\overset{(0,1)}{S(1,a)}} \arrow[d] \arrow[rrd] \& \& {\overset{(0,2)}{S(1,a)}} \arrow[lld] \arrow[rrd] \& \& {\overset{(1,2)}{S(1,a)}} \arrow[lld] \arrow[d] \& \dots \\
    {\underset{0}{S(0,a)}} \& \& {\underset{1}{S(0,a)}} \& \& {\underset{2}{S(0,a)}} \& \dots
    \end{tikzcd}
\end{center}
where we have decorated the sets with the corresponding objects in $ \mathcal D_n $. For each $ 0 \le i < j \le n $ let $ f_{ij}^k : [1] \to n $ be the map with image $ i < j $ in $ n $ and for each $ 0 \le i \le n $ let $ g_i : [0] \to n $ be the map with image $ i $ in $ n $. With these maps $ n $ forms a cocone over $ F_n $ in $ \Delta $. This in turn, gives rise to a cone $ S(n,a) $ over $ S(-,a) \circ F_n $. In particular, by the universal property there is a canonical map
\begin{align}
  u_{n,a} : S(n,a) \to \lim_{\mathcal D_n} S(F_n(-),a).
\end{align}
This map is natural in $ (n,a) $. To see this, notice that any morphism $ f : n \to m $ induces a functor $ \mathcal D_f: \mathcal D_n \to \mathcal D_m $ defined by $ i \mapsto f(i) $ and $ (i,j) \mapsto (f(i), f(j)) $ satisfying $ F_m \circ \mathcal D_f = F_n $. This gives rise to a natural transformation $ \mathcal D_f^* : S(F_m(-), a) \to S(F_n(-), a) $ all components of which are the identity. At the same time, given a morphism $ i_{ab} : a \to b $ we also get a natural transformation $ S(F_m(-), b) \to S(F_m(-),a) $. Both natural transformations induce morphisms on the respective limits, which are compatible with the morphisms $ u_{n,a} $. All of this is summarized by the following diagram
\begin{center}
  \begin{tikzcd}
    {S(m,b)} \arrow[r] \arrow[d, "{u_{m,b}}"'] & {S(m,a)} \arrow[r] \arrow[d, "{u_{m,a}}"]  & {S(n,a)} \arrow[d, "{u_{n,a}}"]  \\
    {\lim_{\mathcal D_m}S(F_m(-),b)} \arrow[r] & {\lim_{\mathcal D_m}S(F_m(-),a)} \arrow[r] & {\lim_{\mathcal D_n}S(F_n(-),a)}
    \end{tikzcd}
\end{center}
which commutes by the uniqueness of the universal property. The outer diagram encodes the naturality of $ u_{n,a} $.

With this at hand, we can now spell out the characterization.

\bp{prop:SingSimpSetProperties}{
  Assuming \eqref{eq:ReIsReSkeleton}, a fuzzy simplicial set $ S \in \bo{sFuz} $ is in the essential image of $ \te{Sing} $ if and only if:
  \begin{enumerate}
    \item $ S(\te{id}_0, i_{ab}) : S(0,b) \to S(0,a) $ is an isomorphism for all $ a,b \in I $.
    
    \item $ \bigl( S(\delta_1,i_{00}), S(\delta_0, i_{00}) \bigr) : S(1,0) \to S(0,0)\times S(0,0) $ is an isomorphism; we denote the inverse by $ \langle -, - \rangle $.
    
    \item The map $ p : S(0,0)\times S(0,0) \to \bo{I} $ defined by
    \begin{align*}
      p(r,r') := \sup \big\{a \in \bo{I} \ \big| \ \langle r, r' \rangle \in S(1,a) \big\}
    \end{align*}
    satisfies the properties
    \begin{itemize}
      \setlength\itemsep{0.4em}
      \item $ p(r,r) = 1 $
      \item $ p(r,r') = p(r',r) $
      \item $\begin{cases}p(r,r')p(r',r'') \le p(r,r''),&\te{(for $\bo{M}=\bo{UM}$)},\\
      p(r,r')p(r',r'') \le p(r,r'')\te{ or }p(r,r'')=0,&\te{(for $\bo{M}=\bo{EPMet}$)}.
      \end{cases}$\\
    \end{itemize}
    for all $ r,r',r'' \in S(0,0) $.

    \item The canonical map $ u_{n,a} : S(n,a) \to \lim_{\mathcal D_n} \bigl(S(F_n(-),a)\bigr) $ is an isomorphism.
  \end{enumerate}
}
\prr{prop:SingSimpSetProperties}

The three equations in the third property correspond to the three axioms of uber metric spaces in the same order as in Definition \ref{def:uberMetricSpace} and they make sure that the axioms are satisfies by the metric we construct. Regarding the last condition, it is interesting to note its resemblance to the \textbf{Segal condition}, which characterizes those simplicial sets which are in the essential image of the nerve of a category.

This proposition demonstrates the additional flexibility alluded to in \ref{A1} that fuzzy simplicial sets provide in contrast to metric spaces for encoding relationships of data. In particular, non-symmetric and higher-order relationships can be encoded, as they arise, for example, in social networks or biological complex systems.

\subsection{Relation of the skeleton adjunction with persistent homology}
\label{sec:relationToPersHomology}

In \cite[Section 4: Persistent Homology]{Spivak09}, David Spivak suggests a category-theoretical formulation for persistent homology using fuzzy simplicial sets. Let us develop this here, although our conclusions will differ somewhat from Spivak's. To  a metric space $(M,d)$, he assigns a fuzzy graph via 
\begin{equation}
    \begin{split}
        G(M)(0,a):&=M\te{ for all }a\in \bo{I},\\
        G(M)(1,a):&=\{(m_1 , m_2 ) \in M \times M ~|~ d(m_1 , m_2 ) \le - \log(a)\}\ .
    \end{split}
  \end{equation}
Now, there is a limit-preserving right adjoint (and not a left adjoint, as Spivak suggested) $F$ to the forgetful functor from $\bo{sFuz}$ to $\bo{1Fuz}$, given by the coskeleton functor, 
\begin{equation}
    \begin{split}
        F:=\te{cosk}_1&:~\bo{1Fuz}\to\bo{sFuz} \\
        \te{cosk}_1(N \in \bo{1Fuz})(m\in\Delta^{\te{op}},a) &= 
        \te{Hom}_{\bo{1Fuz}}(\te{tr}_1(\Delta^m_{a}),N),
    \end{split}
    \label{eq:coskFormulaRepeat}
  \end{equation}
  where the truncation functor $\te{tr}_1$ is defined in Definition \ref{def:truncation}. This is a special case of formula \eqref{eq:coskFormula}, which defines the coskeleton functor.
We also recall that, by Yoneda,  $\te{Hom}_{\bo{1Fuz}}(\te{tr}_1(\Delta^m_{a}),N)\simeq N(m,a)$ for $m\le 1$. Thus,
\begin{equation}
  \begin{split}
      F (G)(n,a) = \te{Hom}_{\bo{1Fuz}}(\te{tr}_1(\Delta^n_{a}), G).
  \end{split}
\end{equation}
Thus, as $G(M)$ is an object of $\bo{1Fuz}$, Spivak suggests to compute the persistent homology of $F(G(M))$. The basic idea thus is to assign to a fuzzy simplicial set $N$ for each $n$ the free $\mathbb{Z}$-module generated by the $n$-simplices of strength at least $a$ in $N$ and to compute the boundary maps in the usual way to obtain homology groups.
With this right adjoint  $F$, and using the results in Propositions \ref{prop:leftRightAdjointOfTruncation}, \ref{prop:realizationOfFuzzySimpSet} and \ref{prop:SingFun}, we can briefly show that, when we assume \eqref{eq:ReIsReSkeleton}, $F(G(M))$
is actually simply equal to $\te{Sing}(M)$:
\bp{prop:simplifySpivaksIdea}{
    Assuming \eqref{eq:ReIsReSkeleton}, i.e.~$\te{Sing}=\te{Sing}^{\te{skeleton}}$, we have
    \begin{equation}
        \begin{split}
            F(G(M)) 
            \simeq 
            \te{Sing}(M).
        \end{split}
    \end{equation}
}
\prr{prop:simplifySpivaksIdea}
This is an interesting result because it shows that the adjunction provided in
\cite{Spivak09}, obtained by choosing $\te{Re}$ to be defined by
\eqref{eq:spivakRe-a->b} was
actually not the most natural choice for defining persistent homology using
fuzzy simplicial sets. Instead, the skeleton adjunction, i.e.~the choice
\eqref{eq:ReIsReSkeleton}, that is, where to $\Delta^n_{a}$ we instead
assign the uber metric space with $d_a(x_i,x_j)=-\log a$ for $i\neq j$, is more natural for that purpose. 

However, we also remark that the adjunction provided by Spivak could be more useful for the purpose of interpolating between data points. Indeed, as shown in Proposition \ref{prop:UMisInSfuz}, $\te{Re}^{\te{skeleton}}\circ\te{Sing}^{\te{skeleton}}=\te{id}_{\bo{UM}}$, and therefore no new points are obtained by going through the adjunction, while Spivaks definition fills the faces between the vertices. Unfortunately, it is more difficult to obtain an explicit description of the realization functor when choosing to define $\te{Re}$ via \eqref{eq:spivakRe-a->b} because Proposition \ref{prop:faceMap} does not have an easily computable counterpart for that choice. This impedes a direct practical application. However, a modification of the definition might result in something fruitful from the practical perspective.
\section{The merge functors}
\label{sec:TheMergeFunctors}

\subsection{Merging fuzzy simplicial sets}

In this section we follow up on question \ref{q3}, i.e.~we rigorously define a recursive merge operation in full generality in the category $\bo{sFuz}$. The fundamental probabilistic operations, that can ultimately be recursively applied to pairs of fuzzy membership strengths of simplices, are the so called t-conorms.

\gd{def:tConorm}{
  A \textit{t-conorm}\index{t-conorm} is a function $T:[0,1]\times [0,1]\to[0,1]$ that fulfills the following axioms:
  \begin{enumerate}
      \item $T(a,b)=T(b,a)$ (Symmetry),
      \item $T(a,T(b,c))=T(T(a,b),c)$ (Associativity),
      \item $T(a,0)=a$ (Bounded from below),
      \item $a\le a'$ and $b\le b'$ implies $T(a,b)\le T(a',b')$ (Functoriality).
  \end{enumerate}
}
\bp{prop:tConormAnnihilatingElement}{
  For every t-conorm, $1$ is an ``annihilating element'', meaning that 
  \begin{equation}
      \begin{split}
          T(a,1)=1=T(1,a)~\forall a.
      \end{split}
  \end{equation}
}
\prr{prop:tConormAnnihilatingElement}

We first define the t-conorm merging procedure for two classical fuzzy sets, that are not simplicial. As a preparation, consider the following pullback in the category of categories,
\begin{equation}
\begin{split}
    \bt
    \bo{cFuz}\times_{\bo{Sets}} \bo{cFuz} \ar{r}\ar{d}\&  \bo{cFuz} \ar{d}{F} \\
    \bo{cFuz} \ar{r}{F} \& \bo{Sets},
    \et
\end{split}
\end{equation}
where $F:\bo{cFuz}\to \bo{Sets}$ is simply the forgetful functor that is defined by $F(A,\xi):=A$. This pullback will be the domain of our merge functor.
\gd{def:classicalTconormMerging}{
    Given a t-conorm $T$ and two classical fuzzy sets $(A,\xi_1)$ and $(A,\xi_2)$, with the same underlying set $A$, define $\te{merge}_{\bo{cFuz}}:\bo{cFuz}\times_{\bo{Sets}} \bo{cFuz}\to\bo{cFuz}$ by
    \begin{equation}
        \begin{split}
            \te{merge}_{\bo{cFuz}}((A,\xi_1),(A,\xi_2))&:=(A,\xi),\\
            \te{where}\quad\xi(a)&:= T(\xi_1(a),\xi_2(a)).
        \end{split}
        \label{eq:classicalTconormMerging}
    \end{equation}
    On morphisms, $\text{merge}_{\bo{cFuz}}$ acts as the identity: $\text{merge}_{\bo{cFuz}}(f,f):=f$.
}
\bp{prop:mergeIsAfunctor}{
  $\te{merge}_{\bo{cFuz}}$ is a functor.
}
\prr{prop:mergeIsAfunctor}
We now proceed to use the equivalence $C:\bo{Fuz}\to\bo{cFuz}$, described in Proposition \ref{prop:isoClassicalAndSimplicialFuzzySets}, to define a merging procedure in $\bo{Fuz}$.

\gd{def:mergeFuz}{
  Given fuzzy sets $S_1, S_2\in \bo{Fuz}$ that fulfil the condition $S_1(0)=S_2(0)$, and a $t$-conorm $T$, define the functor $\te{merge}_{\bo{Fuz}}:\bo{Fuz}\times_{\bo{Sets}} \bo{Fuz}\to \bo{Fuz}$ by
  \begin{equation}
      \begin{split}
          \te{merge}_{\bo{Fuz}}(S_1,S_2):=C^{-1}(\te{merge}_{\bo{cFuz}}(C(S_1),C(S_2))).
      \end{split}
      \label{eq:mergeFuz}
  \end{equation}
}
Again, the domain of $\te{merge}_{\bo{Fuz}}$ is the pullback $\bo{Fuz}\times_{\bo{Sets}} \bo{Fuz}$, where the forgetful functor $F:\bo{Fuz}\to\bo{Sets}$, along which we take the pullback, is defined by $F(S):=S(0)$. Functoriality follows from functoriality of $\te{merge}_{\bo{cFuz}}$ and functoriality of $C$.

\bp{prop:associativity}{
    $\te{merge}_{\bo{Fuz}}$ is
    symmetric in the sense that 
    \begin{equation}
        \begin{split}
            \te{merge}_{\bo{Fuz}}(S_1,S_2)=\te{merge}_{\bo{Fuz}}(S_2,S_1)
        \end{split}
    \end{equation}
    and 
    is associative in the sense that 
    \begin{equation}
        \begin{split}
            \te{merge}_{\bo{Fuz}}(S_1,\te{merge}_{\bo{Fuz}}(S_2,S_3))=\te{merge}_{\bo{Fuz}}(\te{merge}_{\bo{Fuz}}(S_1,S_2),S_3).
        \end{split}
    \end{equation}
}
\prr{prop:associativity}

As a result, the following definition makes sense:
\gd{def:fuzzyPowerSetMerge}{
  Given a finite set $\mathcal{S}$ of fuzzy sets that fulfill the property $s(0)=s'(0)$ for all $s,s'\in\mathcal{S}$, we can recursively define 
  \begin{equation}
      \begin{split}
          \te{merge}_{\bo{Fuz}}(\mathcal{S}):=
          \begin{cases}
              \te{merge}_{\bo{Fuz}}(s\in\mathcal{S},\text{merge}_{\bo{Fuz}}(\mathcal{S}\backslash \{s\})),&\te{for }|\mathcal{S}|>2,\\
              \te{merge}_{\bo{Fuz}}(s,s'),&\te{for }\mathcal{S}=\{s,s'\}.
          \end{cases}
      \end{split}
  \end{equation}
}
Due to symmetry and associativity, we do not even need to impose an order on the set $\mathcal{S}$ and can randomly pick elements for merging until only two elements are left. Furthermore, the following proposition shows that we can first map all $s\in\mathcal{S}$ to classical fuzzy sets, perform the merging of all of them in $\bo{cFuz}$ and then just apply the inverse once:
\bp{prop:mergeAllInCfuzAndApplyInverseOnce}{
  We have 
  \begin{equation}
      \begin{split}
          \te{merge}_{\bo{Fuz}}(\mathcal{S})=C^{-1}(\te{merge}_{\bo{cFuz}}(\{C(s)\}_{s\in\mathcal{S}})),
      \end{split}
  \end{equation}
  where $\te{merge}_{\bo{cFuz}}$ is defined on sets (of classical fuzzy sets) in complete analogy to Definition \ref{def:fuzzyPowerSetMerge} with $\bo{Fuz}$ replaced by $\bo{cFuz}$.    
}
\prr{prop:mergeAllInCfuzAndApplyInverseOnce}

The above preparations now allow us to define a merge operation for two fuzzy simplicial sets.

\gd{def:mergeSfuz}{
    Given $S_1,S_2\in \bo{sFuz}$, which fulfill the condition that  $S_1(n,0) \simeq S_2(n,0)$ for all $n$, define $\te{merge}_{\bo{sFuz}}:\bo{sFuz}\times_{\bo{sSet}}\bo{sFuz}\to\bo{sFuz}$ by
    \begin{equation}
        \begin{split}
            \te{merge}_{\bo{sFuz}}(S_1,S_2)(n,a):=\te{merge}_{\bo{Fuz}}(S_1(n,-),S_2(n,-))(a).
        \end{split}
        \label{eq:fuzzySimpSetDef}
    \end{equation}
    On morphism $\te{merge}_{\bo{sFuz}}$ acts again as the identity: $\te{merge}_{\bo{sFuz}}(m,m):=m$.
}
Here, the proper domain of $\te{merge}_{\bo{sFuz}}$ is the pullback in the diagram 
\begin{equation}
    \begin{split}
        \bt
        \bo{sFuz}\times_{\bo{sSet}} \bo{sFuz} \ar{r}\ar{d}\&  \bo{sFuz} \ar{d}{F} \\
        \bo{sFuz} \ar{r}{F} \& \bo{sSets},
        \et
    \end{split}
\end{equation}
Functoriality follows from functoriality of $\te{merge}_{\bo{Fuz}}$ and the next proposition.

\bp{prop:fuzzySimSetMergeWelldefined}{
    The operation in equation \eqref{eq:fuzzySimpSetDef} is well-defined in the sense that $\te{merge}_{\bo{sFuz}}(S_1,S_2)$ is indeed again a fuzzy simplicial set.
}
\prr{prop:fuzzySimSetMergeWelldefined}

As in \ref{def:fuzzyPowerSetMerge}, we can use this pairwise merge to define a merge operation for finite sets of objects in $\bo{sFuz}$.

\gd{def:simplicialfuzzyPowerSetMerge}{
  Given a finite set $\mathcal{S}$ of fuzzy simplicial sets, which fulfills the condition that
  \begin{equation}
      \begin{split}
          s(n,0)=s'(n,0) \quad \forall~s,s'\in\mathcal{S},~n\in\mathbb{N},
      \end{split}
      \label{eq:conditionOnMergeSet}
  \end{equation}
  we can define 
  \begin{equation}
      \begin{split}
          \te{merge}_{\bo{sFuz}}(\mathcal{S}):=
          \begin{cases}
              \te{merge}_{\bo{sFuz}}(s\in\mathcal{S},\text{merge}_{\bo{sFuz}}(\mathcal{S}\backslash \{s\})),&\te{for }|\mathcal{S}|>2,\\
              \te{merge}_{\bo{sFuz}}(s,s'),&\te{for }\mathcal{S}=\{s,s'\}.
          \end{cases}
      \end{split}
      \label{eq:mergeSetFuzzySim}
  \end{equation}
}
Its domain is $\bo{sFuz}\times_{\bo{sSets}}\cdots\times_{\bo{sSets}}\bo{sFuz}$. The above definition is our final answer to \ref{q3}.

\subsection{Combining metric spaces}

We next investigate the special case, that occurs in step \ref{U1} to \ref{U3} of UMAP, where the set of fuzzy simplicial sets comes from a set of metric spaces $\{(X,d_i)\}_{i}$, where the underlying set is always $X$. 

\bp{prop:canMerge}{
  Any finite set of fuzzy simplicial sets of the form $\{S_i:=\te{Sing}(X,d_i)\}_{i\in I}$ fulfills condition \eqref{eq:conditionOnMergeSet} and can thus be merged with $\te{merge}_{\bo{sFuz}}$ described in \eqref{eq:mergeSetFuzzySim}.
}
\prr{prop:canMerge}

Combining this last result with Definition \ref{def:simplicialfuzzyPowerSetMerge}, that describes the merge operation in $\bo{sFuz}$ and the explicit descriptions of $\te{Re}^{\te{skeleton}}$ and $\te{Sing}^{\te{skeleton}}$ provided in Propositions \ref{prop:realizationOfFuzzySimpSet} and \ref{prop:SingFun}, we can finally derive an explicit description of what we call the \textit{T-combinations} $\top_{\bo{UM}}$ and $\top_{\bo{EPMet}}$ to address the first part of \ref{q4}. Note that we chose $\top$ to be reminiscent of the letter T as in t-conorm. Furthermore, $\top$ can be used in infix notation, as in $(X,d_1)\top (X,d_2)$.

\bp{prop:mergeUM}{
  Assuming \eqref{eq:ReIsReSkeleton}, the functor
  \begin{equation}
      \begin{split}
          \top_{\bo{UM}}:=\te{Re}\circ \te{merge}_{\bo{sFuz}}\circ \te{Sing}^N:\bo{UM}\times_{\bo{Sets}}\cdots\times_{\bo{Sets}}\bo{UM}\to\bo{UM}.
      \end{split}
      \label{eq:mergeFunDefUM}
  \end{equation}
  can be given the following explicit description:
  \begin{equation}
      \begin{split}
          \top_{\bo{UM}}(&(X,d_1),\cdots ,(X,d_N)) = (X,d),\te{ where}\\
          d(x,y) :&= \inf_{x=x_1,\cdots,x_n=y}\sum_{i=1}^{n-1}(-\log(T_1(x_i,x_{i+1}))),
      \end{split}
      \label{eq:mergemetricUM}
  \end{equation}
  where $T_1$ is defined recursively:
  \begin{equation}
      \begin{split}
          T_{k\in \{1,\cdots,N-1\}}(x,y):&=T(e^{-d_k(x,y)},T_{k+1}(x,y)),\te{ and}\\
          T_N(x,y):&=e^{-d_N(x,y)}.
      \end{split}
      \label{eq:T1}
  \end{equation}
  Similarly, the functor 
  \begin{equation}
      \begin{split}
          \top_{\bo{EPMet}}:=\te{Re}\circ \te{merge}_{\bo{sFuz}}\circ \te{Sing}^N:\bo{EPMet}\times_{\bo{Sets}}\cdots\times_{\bo{Sets}}\bo{EPMet}\to\bo{EPMet}.
      \end{split}
  \end{equation}
  can be given the description:
  \begin{equation}
      \begin{split}
          \top_{\bo{EPMet}}(&(X,d_1),\cdots ,(X,d_N)) = (X,d_{\te{EPMet}}),\te{ where}\\
          d_{\te{EPMet}}(x,y) :&=
          \begin{cases}
            \infty,&\te{if }T_1(x,y)=0,\\
            d(x,y)\quad\te{(as in \eqref{eq:mergemetricUM})},&\te{else.}
          \end{cases}
      \end{split}
      \label{eq:mergeEPMetinfExpression}
  \end{equation}
  On morphisms, $\top_{\bo{UM}}$ and $\top_{\bo{EPMet}}$ act as the identity on the functions (that are all assumed to be identical) on the underlying sets (for example, $\top_{\bo{UM}}(f,\cdots,f)=f$).
}
\prr{prop:mergeUM}
\yr{rem:associativityOfTconormResult}{By the symmetry and associativity of the t-conorm\index{t-conorm}, the result of \eqref{eq:T1} does not depend on the order of the metric spaces $(X,d_j)$.}
The recursive definition of $T_1$ in \eqref{eq:T1} simplifies for the max-t-conorm.
\bc{cor:specialCase}{
    For the special case in which $T(a,b)=\max(a,b)$,\index{t-conorm} we obtain
    \begin{equation}
        \begin{split}
            -\log(T_1(x,y))=\min(d_1(x,y),\cdots,d_N(x,y)).
        \end{split}
    \end{equation}
    $~$
}
\prr{cor:specialCase}
\yr{rem:naturalityOfMaxTconorm}{
  Note that in light of the fact that equation \eqref{eq:nSimplexRealizationExplicit}  also contains this minimum operation 
  due to the form of the colimit described in Proposition \ref{prop:colimitUM}
  one might consider this maximum-t-conorm a rather natural choice for $T$.
}
Another useful simplification arises when the metric spaces are $k$-neighborhood spaces. 
\bc{cor:knnMetrics}{
  Suppose that $\{(X,d_i)\}_i$ is a finite set of metric spaces where $d_i$ is given by \eqref{eq:localDists}. Then 
  \begin{equation}
    \begin{split}
      T_1(x_i,x_j)=T(e^{-d_i(x_i,x_j)},e^{-d_j(x_j,x_i)}).
    \end{split}
  \end{equation}
  Therefore, equation \eqref{eq:mergemetricUM} simplifies to 
  \begin{equation}
    \begin{split}
        d(x,y) :&= \inf_{x=x_1,\cdots,x_n=y}\sum_{i=1}^{n-1}(-\log(T(e^{-d_i(x_i,x_{i+1})},e^{-d_{i+1}(x_{i+1},x_{i})}))).
    \end{split}
    \label{eq:knnMerge}
  \end{equation}
}
\prr{cor:knnMetrics}

Now one can recognize that choosing $d_i$ as in equation \eqref{eq:localDists} and $T(a,b)=\max(a,b)$ results, by Corollary \ref{cor:specialCase} and \eqref{cor:knnMetrics}, in the distance 
\begin{equation}
  \begin{split}
    d(x,y)=\inf_{x=x_1,\cdots,x_n=y}\sum_{i=1}^{n-1}\min(d_i(x_i,x_{i+1}),d_{i+1}(x_{i+1},x_{i})).
  \end{split}
  \label{eq:IsomapDistance}
\end{equation}
If one additionally chooses $\rho_i=0$, $\sigma_i=1$ in the definition of $d_i$ in equation \eqref{eq:localDists}, then the expression \eqref{eq:IsomapDistance} is exactly the distance used in Isomap, \cite{Tenenbaum00}. This demonstrates the surprising connection between Isomap and UMAP mentioned in Section \ref{sec:Introduction}. At the same time, Proposition \ref{prop:mergeUM} with arbitrary t-conorm and the additional choices in \eqref{eq:localDists} and \eqref{eq:localDistsUM} naturally generalize Isomap and give rise to a natural combination of the ideas of both UMAP and Isomap in a unified algorithm.

The fact that $\te{Re}\circ \te{Sing} \simeq \te{id}$, as proven in Proposition \ref{prop:UMisInSfuz}, furthermore demonstrates the important point that performing the combination of metric spaces entirely inside of $\bo{UM}$ or $\bo{EPMet}$ does not come at a loss of information. Nevertheless, we show in the next section that the practical implementation details of the embedding, the final step in the paths of Diagram \eqref{diag:umapAlternative} in \ref{q4}, can strongly influence the final outcome.

Finally, note that there is a crucial difference between $\top_{\bo{UM}}$ and $\top_{\bo{EPMet}}$, that can ultimately be traced back to the difference of the metric of the co-equalizers in Proposition \ref{prop:coequalizerUber}. This difference has practical consequences: If we employ $\bo{EPMet}$, then the unification of the local metric spaces does not require to make global distances consistent, whereas $\bo{UM}$ does enforce global consistency with the triangle inequality. From a geometric point of view, we consider $\bo{UM}$ the more natural choice, which is why we employ it subsequently in our algorithm \textbf{IsUMap}.
\section{Application to data visualization: IsUMap}
\label{sec:application}

In this final section, we describe a practical algorithm that applies the theory of previous sections to the problem of data analysis and dimension reduction.
This algorithm completes the description of the alternative path in Diagram \eqref{diag:umapAlternative} in \ref{q4}.

\subsection{Embedding method}
\label{sec:embedding}

While Proposition \ref{prop:mergeUM} gives a solid basis to determine the metric space whose distances are the geodesic distances of the manifold that we assume to underlie the (possibly uniformized) data, it leaves open the question on how to embed this space in the best possible way. While Whitney's and Nash's embedding theorems \cite{whitney1936differentiable,nash1956imbedding}, as well as the Johnson-Lindenstrauss Lemma \cite{johnson1984extensions}, guarantee that good embeddings in Euclidean space are possible for sufficiently high embedding dimensions, challenges arise in particular for low-dimensional embeddings, where the so-called ``crowding problem'' \cite{vandermaaten08a,olszewski2025data} causes considerable difficulties. Intuitively, the crowding problem arises whenever the low dimensional embedding space does not have enough degrees of freedom to accommodate the degrees of freedom of the data manifold. Many more points can be approximately equally far away from each other in a high-dimensional space than in a low-dimensional space and the number of points that can be placed in equally spaced distances on a sphere of constant radius around the center of a given point increases with the number of dimensions. If the embedding method aims at preserving distances, then, the more faithfully shorter distances are preserved, the more distorted medium distances must become and vice versa. In particular the preservation of medium and long distances then leads to a situation where all points crowd together in the center. Clusters that are distinguishable in the metric space then often overlap in the low-dimensional embedding space as a consequence of these effects.

Methods like t-SNE \cite{vandermaaten08a} and UMAP \cite{McInnes18} present strategies to alleviate the crowding problem. t-SNE does so very explicitly by employing the student's t-distribution, while UMAP implicitly employs a heavy-tailed distribution in their smooth approximation and a negative sampling strategy, both of which we explain in more detail in Section \ref{sec:effectsUMAP}, when discussing the results.

However, these solutions exhibit the problem that the distortion from an embedding, that attempts perfect preservation of distances, is not clear anymore. Especially in the case of UMAP, this can lead to the formation of strong artifacts, as one can observe in Table \ref{fig:comparisonUMAPandIsumap}. Furthermore, it is unclear if and how the full information about the clusters in the original metric space is used or not. We thus developed a new embedding method that is supposed to improve upon these shortcomings.

Our embedding method combines a metric embedding with a subsequent optimization procedure that pulls clusters of points apart that a clustering method distinguished in the original geodesically complete metric space (which we in turn construct according to Proposition \ref{prop:mergeUM}).

For the metric embedding, natural choices are classical, metric or non-metric multidimensional Scaling (MDS, cf.~\cite{Torgerson1952}, \cite{borg05}).

\begin{figure}[H]
  \centering
  \includegraphics[width=0.9\textwidth]{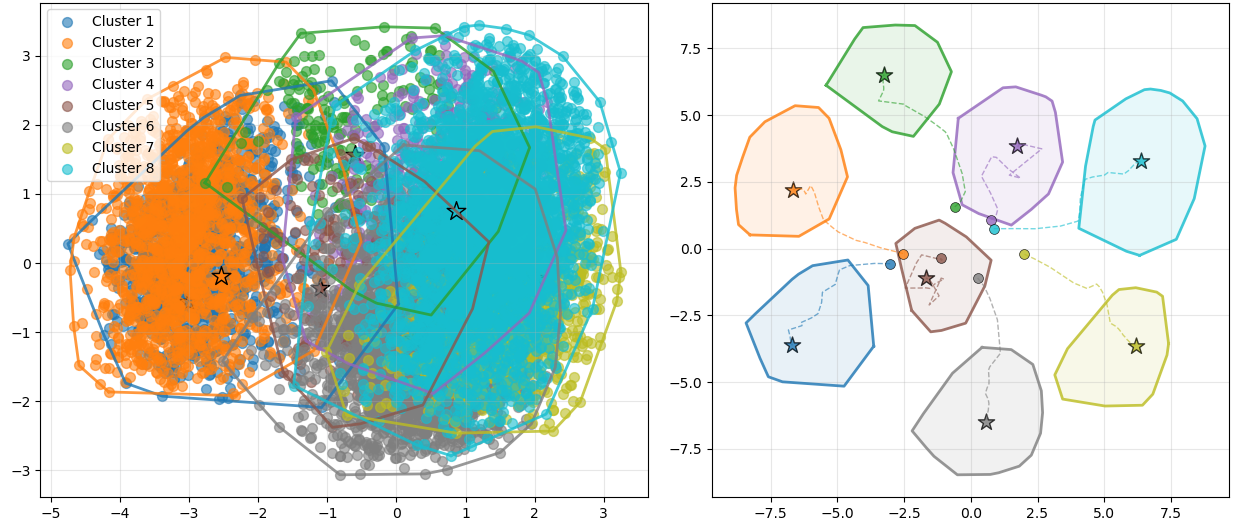}
  \caption{Left: The crowding problem is illustrated by overlapping convex hulls of clusters. Right: The Cluster Separation Optimization pulls the clusters apart and traces out paths in the embedding space.}
	\label{fig:dynamicCRCCbefore}
\end{figure}

To determine the clusters, any clustering method can be applied to the geodesic distance matrix and is compatible with our pipeline. For the simulations that we present in subsequent sections, we have chosen to illustrate the idea with the Linkage 
\cite{Linkage2001,mullner2011linkage} and Leiden clustering \cite{Traag2019Leiden} 
methods but others would work too. The clusters are pulled apart via a stochastic gradient descent procedure that minimizes the overlap of the convex hulls of the clusters. Furthermore, the inter-cluster distance is optimized as a monotone function of the distance of the medoids of the clusters in the original metric space, which allows to retain information about the distances in the original space in the embedding. One very nice feature of the procedure is that one can visualize not only the end result but also the paths that the (medoids of the) clusters took during the optimization that separated the clusters. Before introducing more formal arguments, we want to give the reader an intuition with the plot depicted in Figure \ref{fig:dynamicCRCCbefore}.

%\graphicspath{{figures/}}

The figure shows an embedding of single cell RNA sequencing data, that we will describe in more detail in Section \ref{sec:visualization}. For now, the important point is that the data is high-dimensional and the embedding on the left side of the figure thus exhibits the crowding problem. As a consequence, clusters, marked by their surrounding convex hulls, that can be detected in the original metric space,  overlap in the embedding space. On the right side of the figure, one can see how those clusters are then pulled apart by an optimization procedure that aims at minimizing the overlap of the convex hulls. The colored round points close to the center are the starting points of the medoids of the clusters before the cluster separation optimization procedure and the colored stars are the medoid positions in the embedding after separation. Between them, we can see a path that is traced out by the procedure that helps to understand the relationship between the metric embedding and the separated configuration. It can help to understand how entangled the clusters really are in the high-dimensional space. Furthermore, one can then also re-insert the data points into the separated convex hulls and color them with label information, as the reader can see on the right side of the corresponding Figure \ref{fig:BMMC}.

We believe this novel visualization to be useful because the user immediately obtains some intuition about the distortion that was introduced into the visualization in comparison to the local optimum of a metric preserving embedding, which is something that is not easily available in other methods.
Another advantage is that users can employ cluster algorithms that are already prevalently used in their field, such as Leiden clustering in single cell RNA sequencing analysis frameworks \cite{satija2015Seurat,stuart2019Seurat,wolf2018scanpy}. This allows to obtain a visualization that is both, aligned with the clusters encoding relevant information in the original space, and rich in information about the metric structure of the embedding.
Below, we provide an algorithmic summary.

\subsection{Algorithmic description}
\label{sec:algorithmicDescription}

Input: $X\subset \mathbb{R}^n,~|X|<\infty,~k\in\mathbb{N},~m\le n$.
\begin{enumerate}
  \item As in \ref{U1}, we split $X\subset \mathbb{R}^n$ into $N:=|X|$ metric spaces $(X,d_i)$, where $d_i$ is either defined by \eqref{eq:localDists} or \eqref{eq:localDistsUM}, as explained in Section \ref{sec:DiscreteApprox}. $\rho_i$, $\sigma_i$ and $f$ can be chosen according to the needs of the application. For example, if uniformization is desired, a simple and effective choice is $\rho_i=0$ and $\sigma_i=d(x_i,x_{i_k})$.
  \item Then we apply the T-combination functor defined in Proposition \ref{prop:mergeUM}, or more precisely, equation \eqref{eq:mergemetricUM} or \eqref{eq:mergeEPMetinfExpression}, to obtain the metric space $(X,d)$. Note that one can choose an arbitrary t-conorm. A computationally effective choice is the max-t-conorm because of Corollary \ref{cor:specialCase}. However, for high-dimensional datasets the probabilistic sum often turns out to be advantageous.
  \item Next, we initialize an embedding of $(X,d)$ into $\mathbb{R}^m$ with classical multidimensional scaling (cMDS).
  \item Additionally, we apply metric multidimensional scaling if desired.
  \item Finally, we may apply a clustering method (e.g.~Leiden or Linkage clustering) to $(X,d)$ to obtain an index set of clusters and pull the convex hulls of the corresponding clusters in the embedding space apart via stochastic gradient descent, while ensuring that the inter-cluster distance is at least as big as a monotonic function of the distance of the medoids of the original clusters in $(X,d)$.
\end{enumerate}
Output: $Y\subset \mathbb{R}^m,~|Y|=|X|$, usually $m=2$.

We next describe the concrete algorithms that we use to perform the mathematical steps described in the list above. We represent the metric space by a sparse tensor $R$, whose $i, j, k$-th component denotes the distance between points $j$ and $k$ in the $i$-th local space. All operations above are then performed on this tensor:
\begin{enumerate}
  \item 
  The $k$ nearest neighbors can be found with a time complexity of $\mathcal{O}(N(k+\log N))$ using $kd$-Trees (cf.~\cite{Bentley1975}). However, the runtime does not scale well for high $d=n$ but one can achieve even better results using approximate gradient descent algorithms that can be parallelized (cf., for example, \cite{douze2024faiss}). One then assigns non-infinite distances to the sparse tensor $R$ according to \eqref{eq:localDists} or \eqref{eq:localDistsUM}, while infinite distances are not explicitly assigned (but respected during subsequent computations) to preserve sparsity.
  \item To apply the combination operation, \eqref{eq:knnMerge}, we 
  only need to loop over the non-zero entries of the sparse tensor for computing the recursive t-conorm because $T(a,0)=a$ by Def.~\ref{def:tConorm}. The result is again assigned to a sparse matrix. This makes the algorithm quite efficient.
  
  Thereafter, we apply the log-operation and then execute the Dijkstra routine to compute the infimum in \eqref{eq:mergemetricUM}.
  Dijkstra has a time complexity of $N(N\log(N)+E)$ for computing all shortest paths between any two points. Since the number of edges $E$ is at most equal to $kN$ (because we consider only the $k$-nearest neighbor graph), and $k$ usually depends on $N$ at most logarithmically, we have $E\propto N\log(N)$ as well and hence obtain a complexity of roughly $N^2\log (N)$, which bounds the overall complexity of our algorithm. We implemented a parallel version of Dijkstra, that can roughly reduce the computation time by a factor proportional to the total number $T$ of CPU threads, i.e.~the sequential complexity becomes proportional to $N^2\log(N) / T$, resulting in a significant speedup compared to previous implementations of similar algorithms like Isomap.
  \item For the metric embedding, classical MDS (cMDS) requires to determine the $m$ biggest eigenvectors of the matrix $D$. This can be done using the Lanczos algorithm (cf.~\cite{lanczos1950}, and later improvements), which has a high efficiency due to stepwise determination of the eigenvectors, of which often only $m=2$ or $m=3$ have to be computed for data visualization. 
  \item Optionally, metric or non-metric MDS is thereafter initialized with cMDS and performed using stochastic gradient descent (SGD), which we implemented ourselves in pytorch \cite{pytorch2019}.
  Note that there is also the SMACOF algorithm \cite{SMACOF} to solve the metric MDS problem. It has a higher chance to evade running into a local minimum but also a higher time complexity of $\mathcal{O}(N^2)$, though there are parallel versions, \cite{Orts2019}. With small batch sizes, SGD can achieve a complexity significantly lower than $\mathcal{O}(N^2)$.
  \item A chosen clustering algorithm (like Leiden or Linkage clustering) is applied to the original fuzzy graph or geodesic distance matrix and returns cluster indices that we use to obtain clusters in the embedding space. We also determine the medoids and their distances of the original distance matrix, which correspond to special points in the embedding space. Thereafter we compute the convex hulls of the embedded clusters with scipy \cite{2020SciPy-NMeth} and finally use stochastic gradient descent (again implemented in pytorch) to pull the convex hulls apart. Specifically, we minimize the following objective via gradient descent
  \begin{equation}
    \begin{split}
      L(\theta) = \alpha \sum_{i \ne j=1}^c \sum_{k=1}^{n_j} \tanh(f_{ij}^k(d(T_{\theta_j}x_{j,k}, T^H_{\theta_i}H_i)) + \beta ||\theta||_2^2,
    \end{split}
  \end{equation}
  where $c$ is the number of clusters, $n_j$ is the number of points in cluster $j$, and $x_{j,k}$ is the $k$-th point in cluster $j$, and $T_{\theta_j}$ transforms all points of cluster $j$ in the same way through some rotation and/or translation, parameterized by $\theta_j$, and $H_i$ is the convex hull of cluster $i$ (described by the set of vertex points of the hull), and the action of $T^H_{\theta_i}$ on $H_i$ is induced by the action of $T_i$ on all points in the convex hull, and $d(a,B)$ measures the minimal distance of the point $a$ to the boundary of the convex hull $B$, and $f_{ij}^k$ is a function whose effect depends on whether $T_{\theta_j}x_{j,k}$ is inside or outside of $T^H_{\theta_i}H_i$. If $T_{\theta_j}x_{j,k}$ is outside, then $f_{ij}^k(d)$ defaults to $D_{ij}+d$ and otherwise to $\max(0,D_{ij}-d)$, where $D_{ij}$ is the pre-computed inter-cluster distance of the original metric space, measured by the distance of the medoids of the clusters in the original metric space $(X,d)$. In this way, the intercluster distance is a monotone function ($\tanh$) of relevant distances of the original space. A monotone function  is used to prevent the clusters from becoming too distant for a nice visualization. Other monotonic functions could be used as well. This formulation might seem complicated but it really just is one formalization of the idea of imposing a repelling force on overlapping convex hulls, while retaining some information about intercluster distances. The drift term $\beta ||\theta||_2^2$ ensures that clusters do not drift apart too much, though we found that the setting is robust enough to allow for $\beta=0$ in all of our experiments, especially when only employing translations.
  
  In practice a more efficient samplig strategy is implemented that does not require to sum over all points. Pytorch then backpropagates through the entire stack of computations involved in finding the minimum distance to the convex hull for each given point and cluster. We precompute the hulls to increase efficiency. In principle, more efficient implementations should be possible involving only vertex points of hulls. 
  However, since there are only $3c$ parameters, where $c$ is the number of clusters ($2c$ for translation and $c$ for rotation, if used), the procedure is comparatively efficient, though we expect further improvements to be possible. Since we use pytorch, the implementation can run on the GPU.
\end{enumerate}
This concludes the description of our algorithm both from a theoretical and from a practical point of view, including rough bounds on its time complexity, and hence fully answers \ref{q4}.

In addition to the algorithm itself, we also provide elaborate functions that can visualize the optimization path during the cluster optimization phase, as illustrated in Figure \ref{fig:dynamicNewsgroup}.

Note that the algorithm presents an alternative path to UMAP as presented in diagram \ref{diag:umapAlternative} but its metric embedding procedure (Steps 1 to 3) is also similar to Isomap \cite{Tenenbaum00} (even though our method allows for more general T-combination procedures and implements a more efficient Dijkstra routine). This makes a previously unsuspected connection between those methods apparent, and shows how to combine them in a meaningful way. Since the algorithm represents this combination of Isomap and UMAP, and makes strong use of the category $\bo{UM}$, we decided to call it \textbf{IsUMap}.

All source code is freely available at: \href{https://github.com/LUK4S-B/IsUMap}{https://github.com/LUK4S-B/IsUMap}.

\subsection{Visualizations}
\label{sec:visualization}

In this subsection, we present some empirical simulation results. Source code for reproducing them is provided in the Simulations folder of the above referenced github repository.

Table \ref{fig:comparisonUMAPandIsumap} is dedicated to low-dimensional manifolds. We show some visualizations created with IsUMap (where we use \eqref{eq:localDists}, set $\rho_i=0$, $\sigma_i=d(x_i,x_{i_k})$ and use classical MDS without additional clustering steps) and compare them with those produced by UMAP and Isomap.
\\
One can observe that IsUMap imposes only a small distortion on the embedding, similar to Isomap, while UMAP significantly distorts the embedding. Particularly interesting is the first row, where we used as input data a hemisphere with a non-uniform data distribution (the top of the hemisphere has a higher density of points). One can see that IsUMap can succesfully uniformize the data distribution of the embedding (especially notable when looking at how the area of the green points expanded in the embedding), simply by changing $\sigma_i=1$ (the default choice of Isomap) to $\sigma_i=d(x_i,x_{i_k})$. While UMAP also shows some uniformization effect, the distortion imposed on the embedding makes it somewhat less pronounced. Of course, whether one wants an embedding with uniform distribution depends on the application but in any case IsUMap can deliver whatever option is desired.

\begin{figure}
  \centering
  \includegraphics[width=0.8\textwidth]{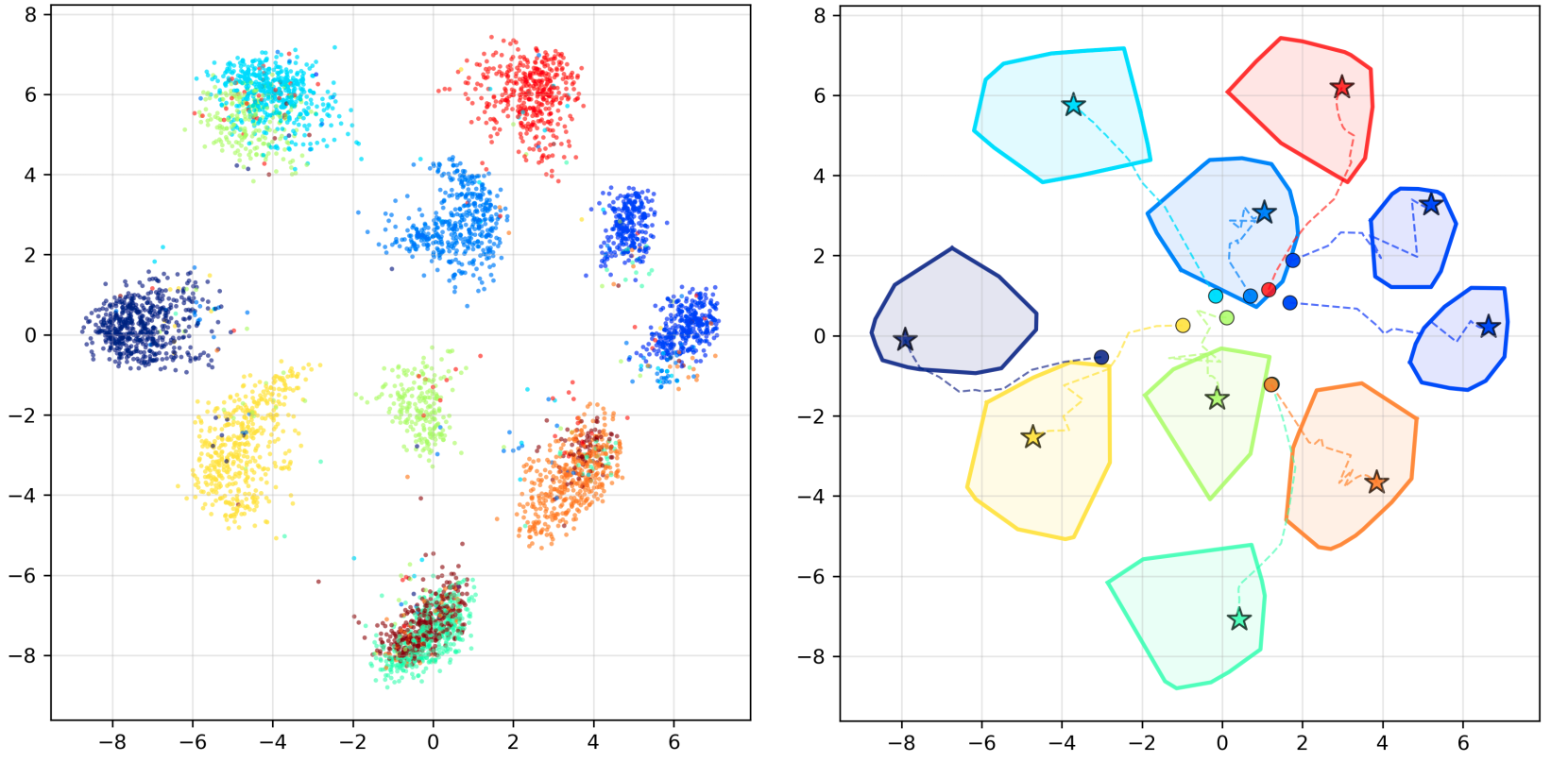}
	\caption{Custer Separation Optimization paths on MNIST dataset embeddings.}
	\label{fig:dynamicMNIST}
\end{figure}

\begin{figure}
  \centering
  \includegraphics[width=0.8\textwidth]{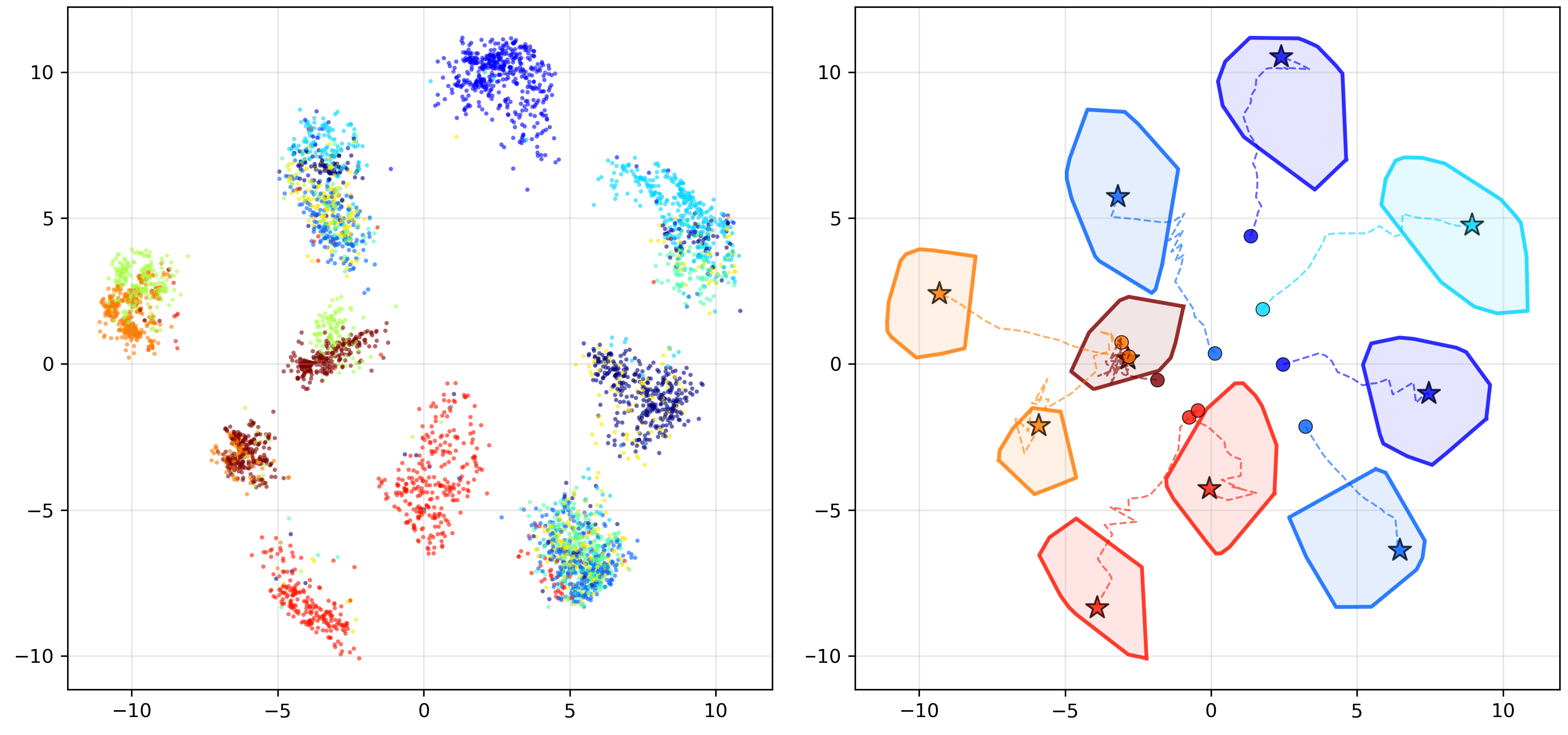}
	\caption{Custer Separation Optimization paths on FashionMNIST dataset embeddings.}
	\label{fig:dynamicFashionMNIST}
\end{figure}

Table \ref{fig:MNISTandFashionMNIST} presents embedding results of the well-known high-dimensional datasets MNIST and FashionMNIST. We use similar parameters for IsUMap as above but with additional clustering steps, employing linkage clustering in our cluster separation pipeline. One can observe that the additional separation clearly improves upon Isomap's result (in the sense that labels and visualized clusters are much more aligned). Furthermore, we also show the optimization paths in 
Figure \ref{fig:dynamicMNIST} and  \ref{fig:dynamicFashionMNIST}. 
Arguably, this plot provides more information than the corresponding UMAP plot because the initial round points show how a metric preserving embedding would look like, i.e. how close the medoids of the clusters actually are in the original metric space.

\begin{figure}
  \centering
  \includegraphics[width=0.9\textwidth]{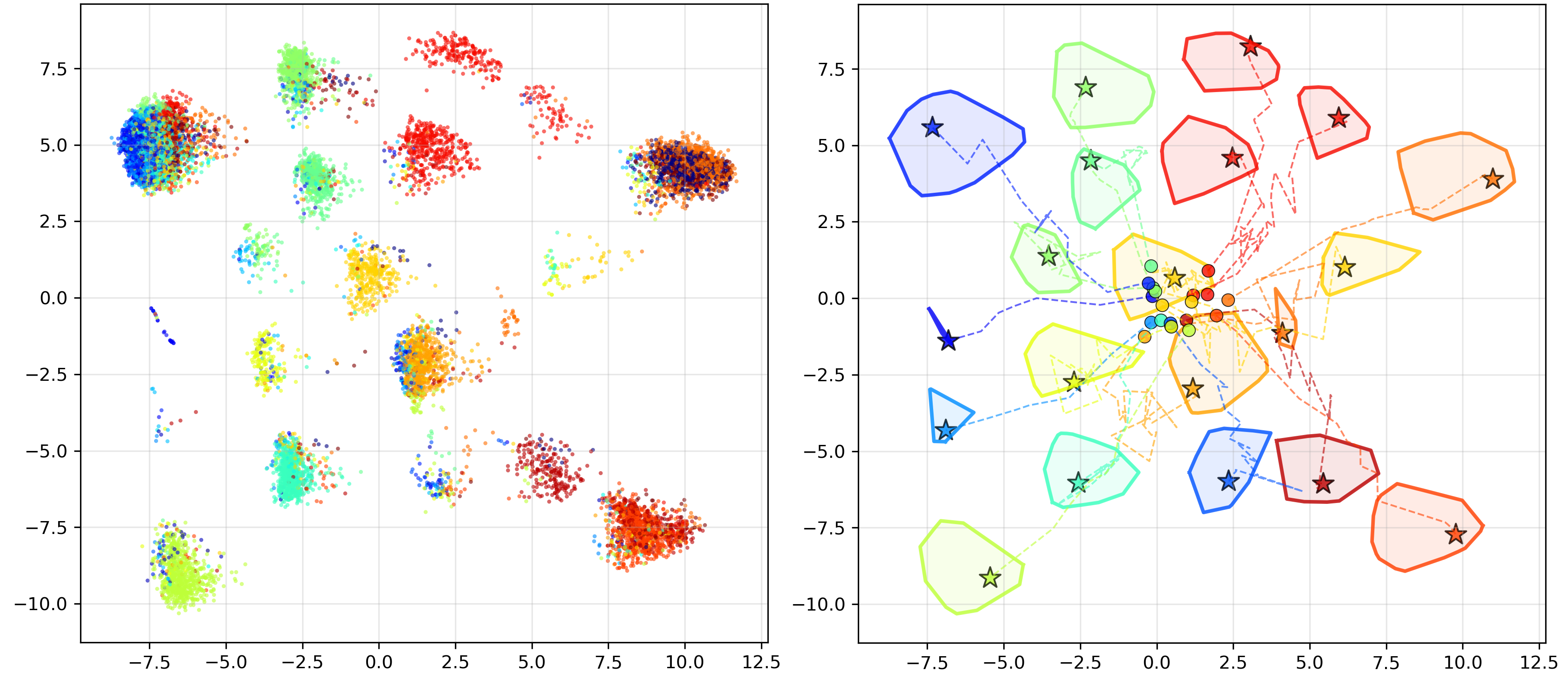}
  \caption{The Cluster Separation Optimization paths of IsUMap during the embedding of the 20-newsgroups dataset.}
	\label{fig:dynamicNewsgroup}
\end{figure}

Next, we present results for high-dimensional text data to showcase the algorithm in the domain of document embedding and natural language processing. The documents of the 20-newsgroups dataset (retrieved through scikit-learn \cite{scikit-learn}) are vectorized with a bag-of-words approach to obtain a matrix that essentially contains token counts for each document. Thereafter the counts are normalized to obtain probability distributions over tokens (for each document) and a distance matrix is finally computed between the distributions using the Hellinger distance. This distance matrix is then embedded using UMAP, IsUMap and t-SNE and the results are compared in Table \ref{fig:newsgroups}. The colors show how well the document categories fit the clusters. Furthermore, the IsUMap cluster separation optimization process is visualized in Figure \ref{fig:dynamicNewsgroup}. 
\\
Here, the paths in Figure \ref{fig:dynamicNewsgroup} show that the crowding problem described in Section \ref{sec:embedding} is very severe for this dataset and makes an embedding without cluster separation strategy almost useless.
Arguably, the IsUMap clusters are the most clear ones because the different clusters do not stick together as tightly as the ones in the UMAP and t-SNE plots, while Figure \ref{fig:dynamicNewsgroup} additionally contains information about how they relate to each other. This is especially helpful to distinguish classes if no label information is provided. However, we admit that there is always a subjective element to the judgment of such visualizations.

\begin{figure}
  \centering
  \includegraphics[width=0.49\textwidth]{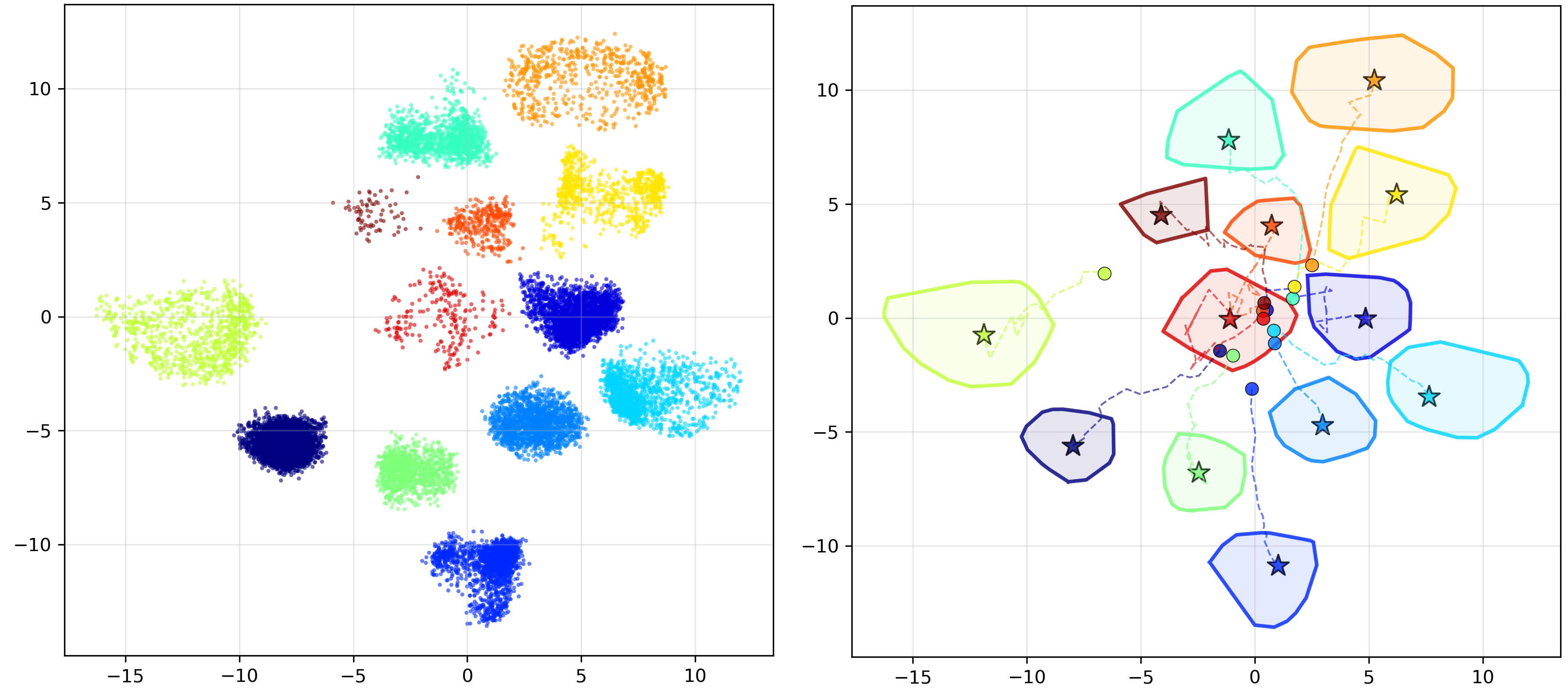}
  \includegraphics[width=0.49\textwidth]{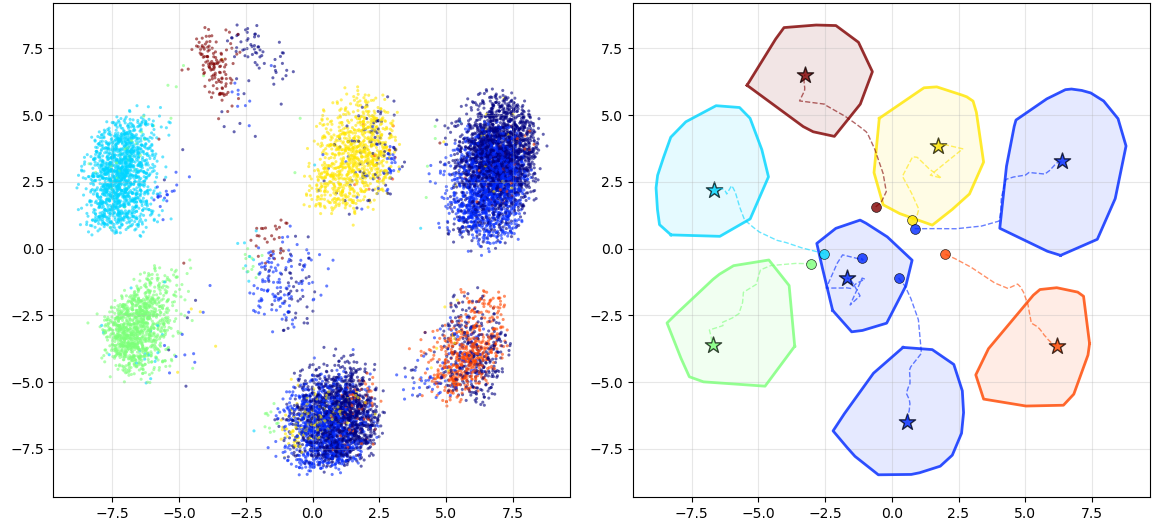}
  \caption{Left: IsUMap with cluster separation, employing Leiden clustering, applied to the BMMC dataset (cf.~Table \ref{fig:BMMC}), colored with Leiden cluster labels. Right: Cluster Separation Optimization paths for CRCC (Leiden) (cf.~Table \ref{fig:CRCC}).}
	\label{fig:dynamicBMMCLeiden}
\end{figure}

Table \ref{fig:BMMC} shows results on single cell RNA sequencing data.
The cells are bone marrow mononuclear cells (BMMCs) of healthy human donors. The dataset was one of NeurIPS 2021 open problem benchmarking datasets \cite{luecken2021sandbox} and we follow the scanpy tutorial for preprocessing it \cite{scanpyTutorial}, including a PCA projection to 50 dimensions, which is conventional in single cell sequencing workflows.
Here, we test how well the different embedding methods align with the Leiden clustering method that is commonly employed in single cell sequencing workflows.  \\
IsUMap's cluster separation optimization would in principle allow to separate Leiden cluster labels perfectly as shown in Figure \ref{fig:dynamicBMMCLeiden} because it can determine the convex hulls of the separated clusters in the embedding based on the Leiden algorithm itself. Even though this would render the correlation between cluster shapes and colors meaningless, the inner-cluster distributions as well as the inter-cluster distances and the cluster optimization paths would still contain meaningful information because they are based on the geodesic distances of the original metric space. We thus believe this might be of some interest if a meaningful visualization is desired that strictly separates certain cluster groups.\\
However, more often the aim might be to cross-check if the cluster shapes obtained by the embedding method correlate with the label colors obtained from an independent clustering method. For example, one might want to cross-check if UMAP's cluster shapes correlate with the coloring of Leiden clustering. IsUMap's cluster separation optimizer facilitates to perform such cross-checks very systematically because it allows to employ any clustering algorithm in its pipeline. For example, in Table \ref{fig:BMMC} (2nd column), we let IsUMap optimize cluster separation based on Linkage clustering but use Leiden clustering to color the cells. The combination of this plot with the paths on the right side of Figure \ref{fig:dynamicBMMCLeiden} can then provide more information about the metric relations of the datapoints than a single embedding alone. Furthermore, for each new clustering method, we could use our method to produce several more of such cross-checks.

Finally, we look at a dataset from \cite{Krauss2025}, containg samples from KRAS tumor organoid cells, sequenced in a study about colorectal cancer (CRC). We thus refer to them as colorectal cancer cells (CRCCs). The data was preprocessed using Seurat \cite{stuart2019Seurat}. In addition, PCA was employed to reduce the dimension to 50 before applying any other dimension reduction method. Embeddings of UMAP, IsUMap and t-SNE are shown in Table \ref{fig:CRCC}. The top row shows labels obtained from Leiden clustering (while IsUMap employed Linkage for clustering) and the bottom row shows the log-scaled expression count of gene SMOC2. The latter is an important marker because it turned out to be highly expressed in AKPE organoids, while being absent in AKP organoids. We find that IsUMap produces more clearly separated clusters in this case.

%\graphicspath{{figures/}}
\begin{table}[H]
	\centering
	\tiny
	\begin{tabular}{>{\centering\arraybackslash}m{1.4cm}|*{4}{>{\centering\arraybackslash}m{2.6cm}|}}
	 & \textbf{Input data} & \textbf{UMAP} & \textbf{IsUMap}  & \textbf{Isomap} \\
		\hline 
		\textbf{(1) Hemisphere} & \includegraphics[width=0.15\textwidth]{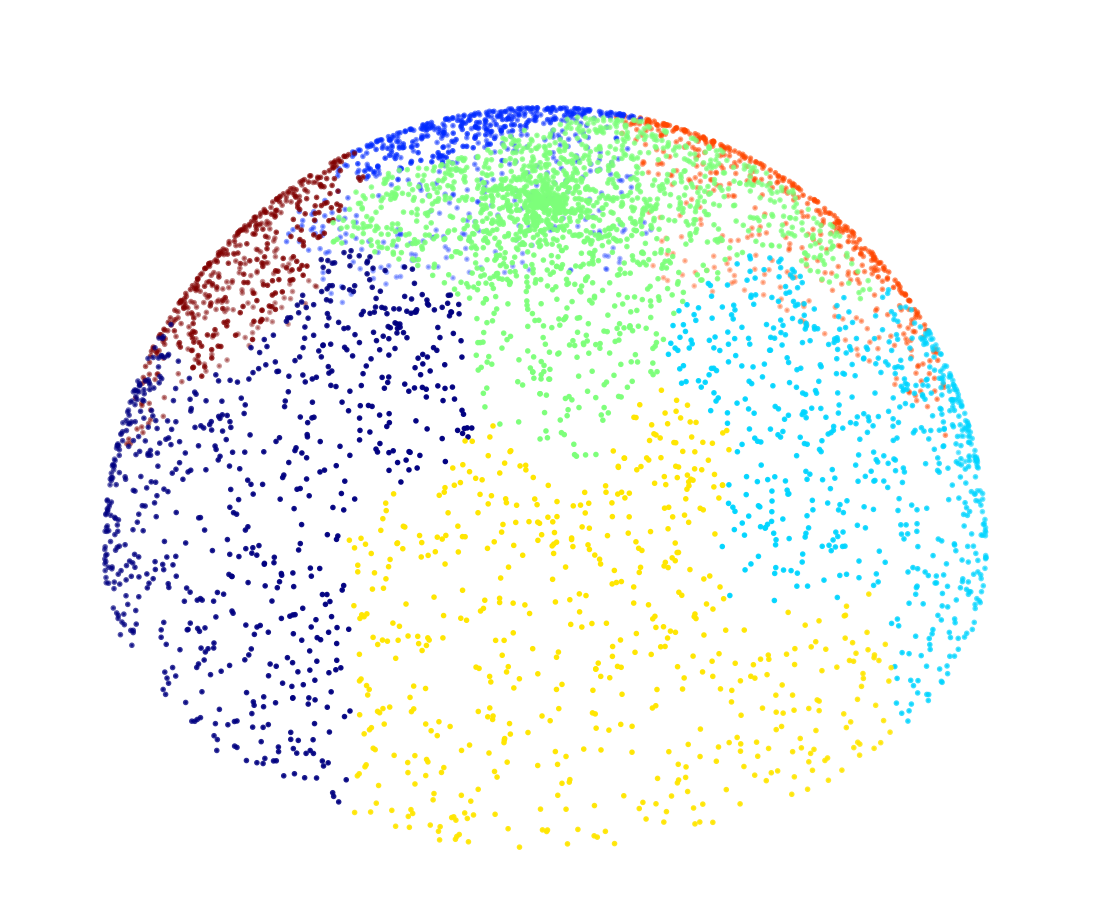} & \includegraphics[width=0.15\textwidth]{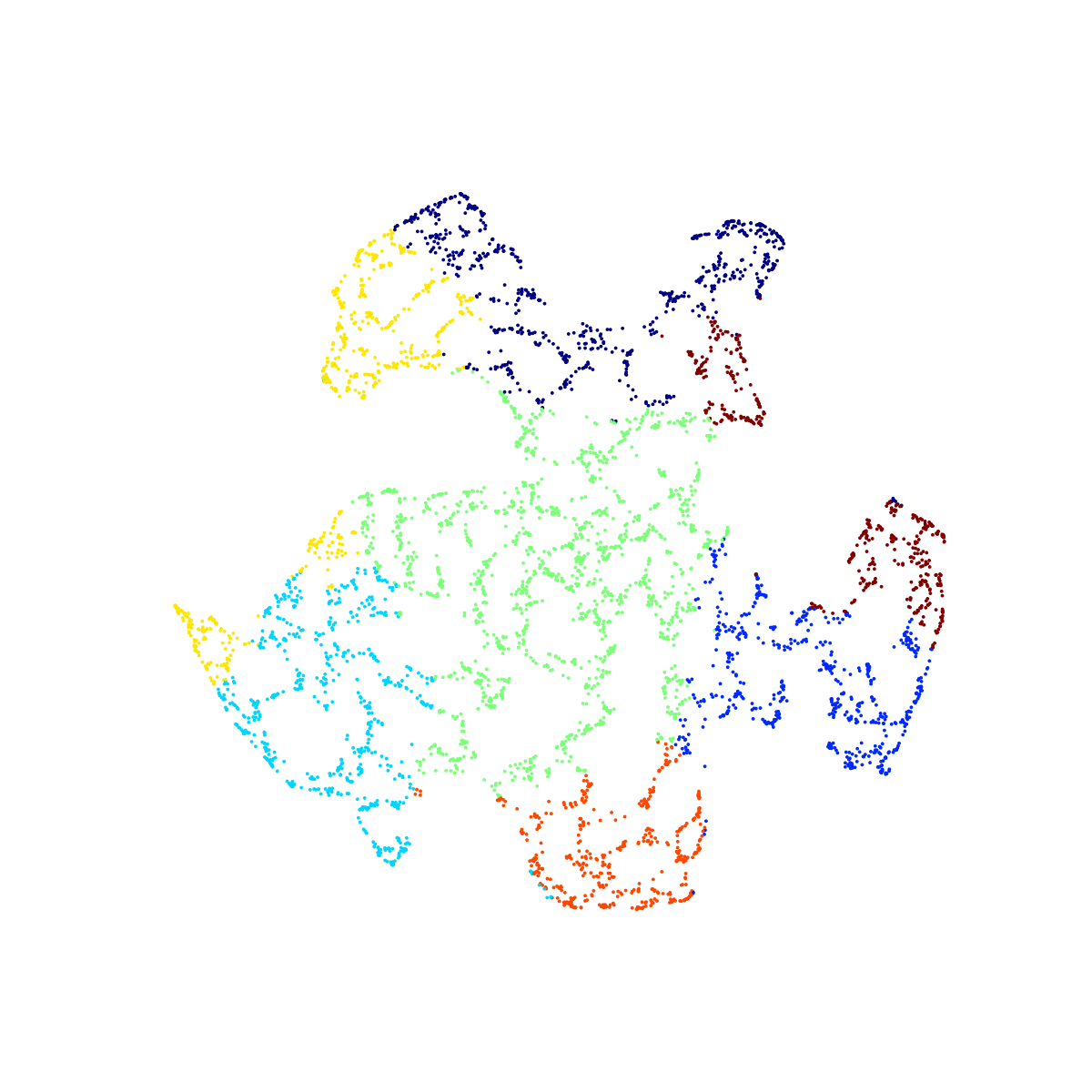} & \includegraphics[width=0.15\textwidth]{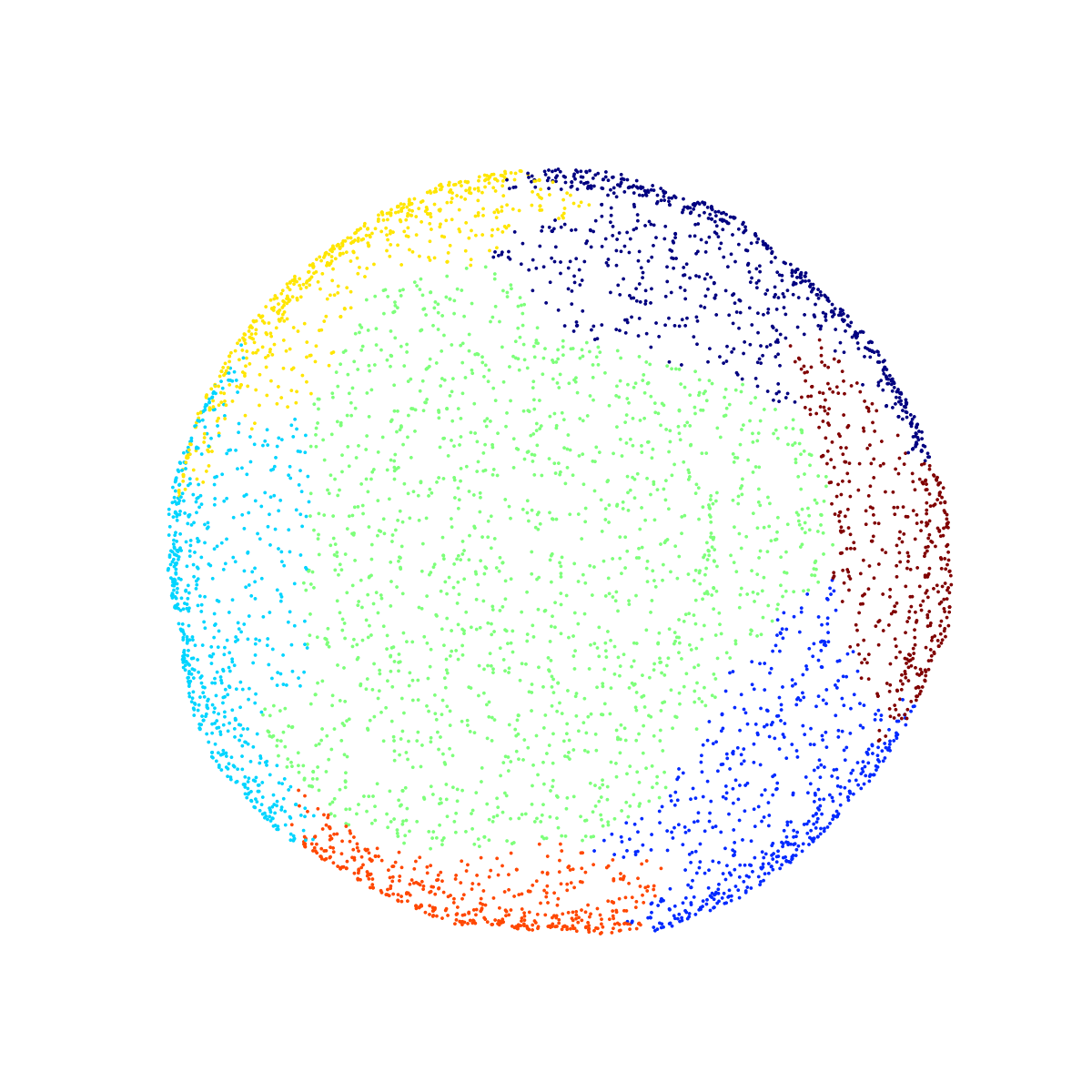}& \includegraphics[width=0.15\textwidth]{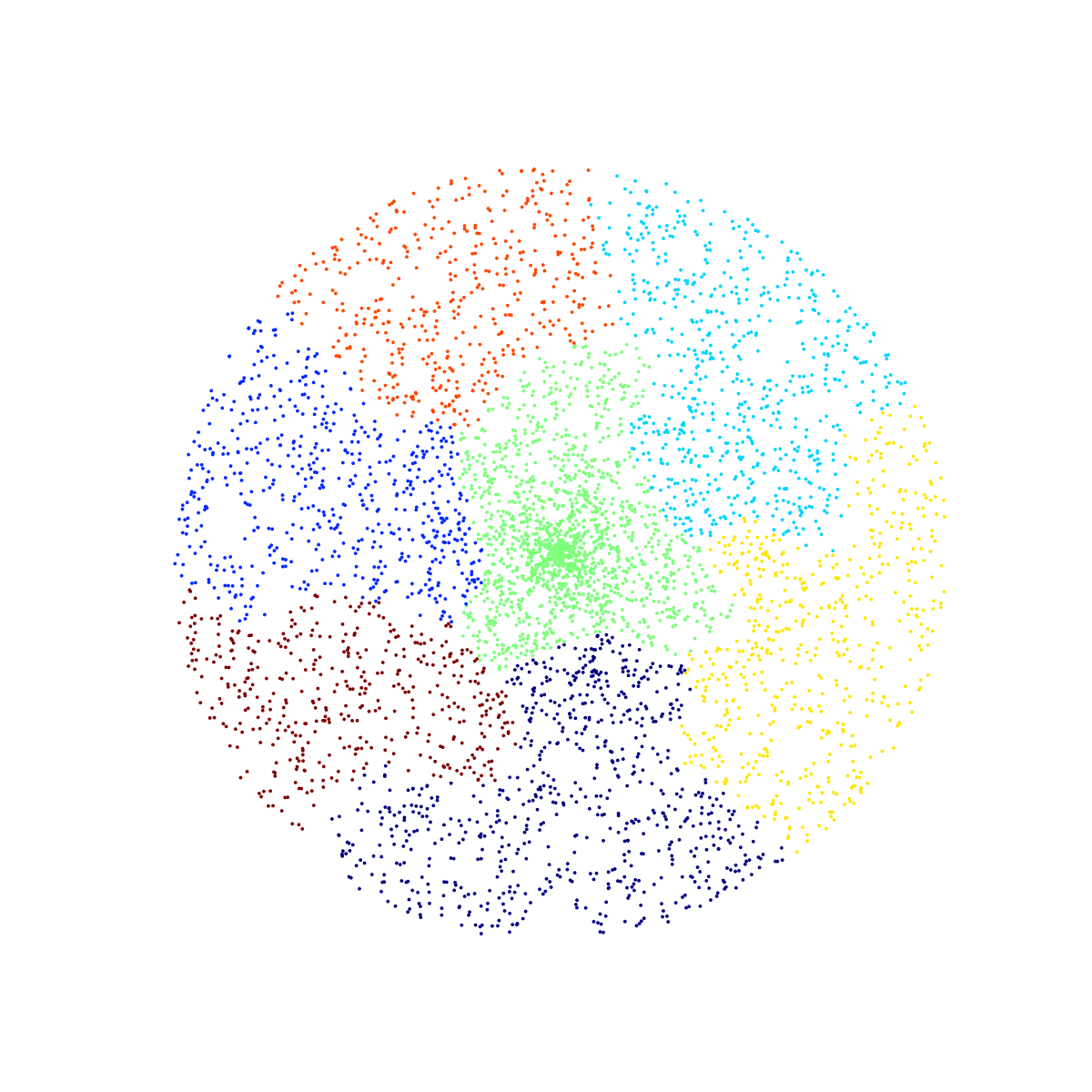}  \\
    \hline
		\textbf{(2) Torus} & \includegraphics[width=0.15\textwidth]{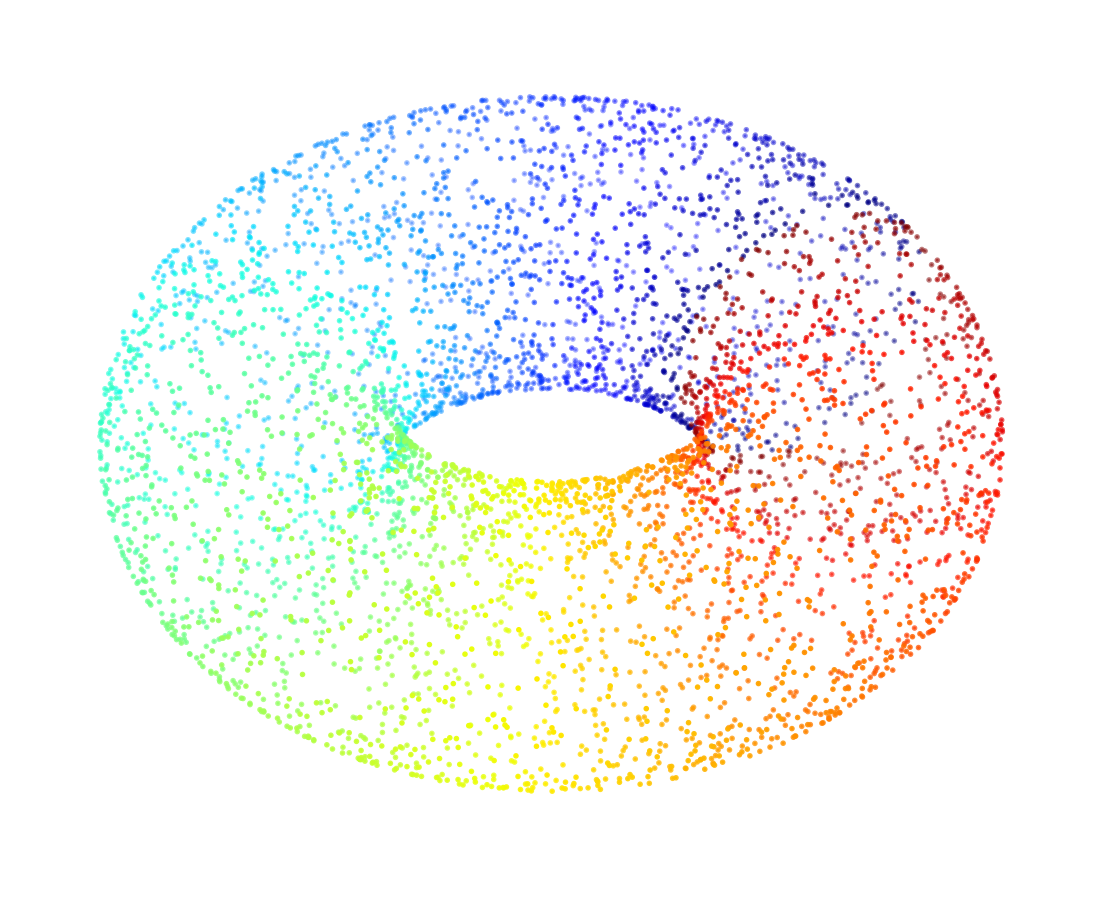} & \includegraphics[width=0.15\textwidth]{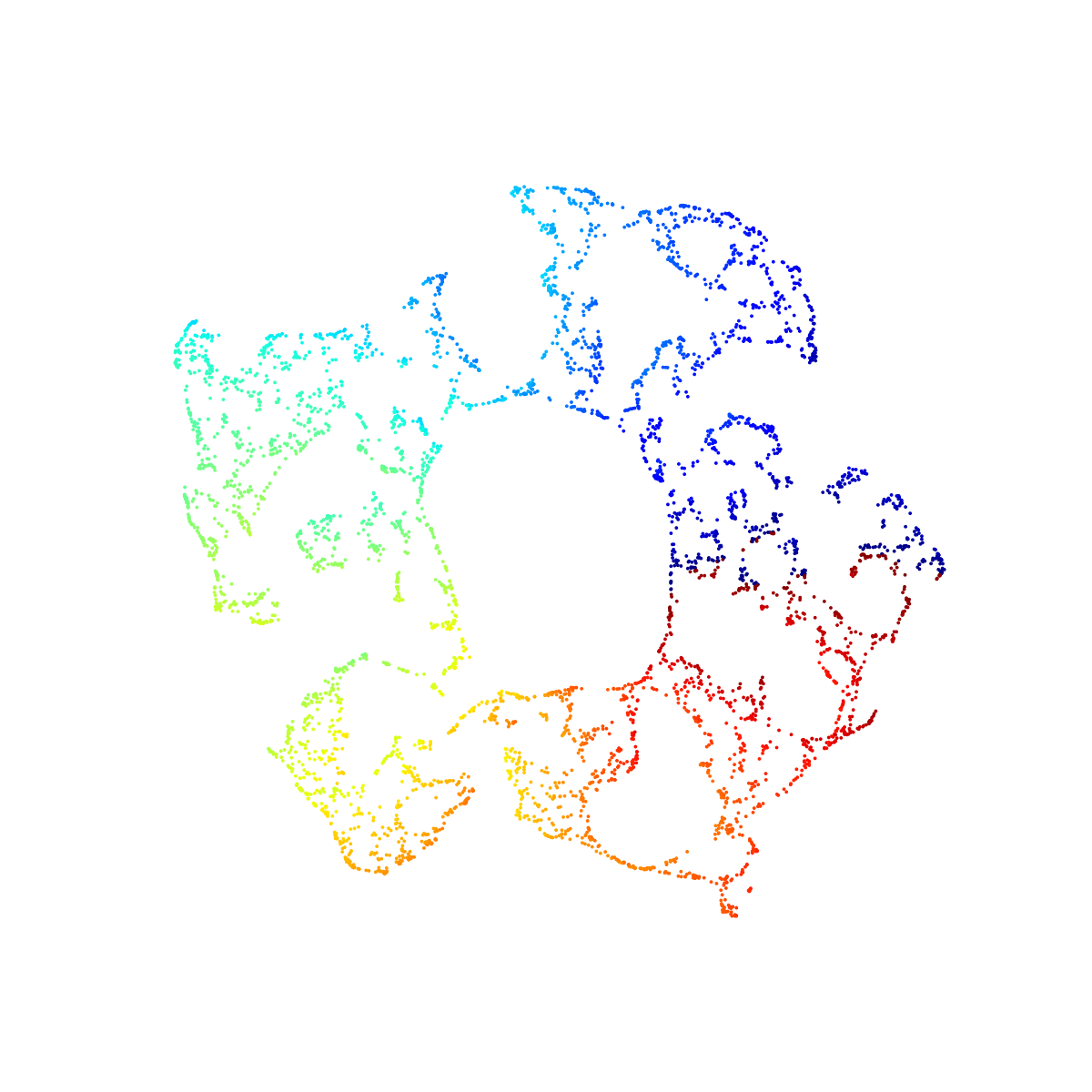} & \includegraphics[width=0.15\textwidth]{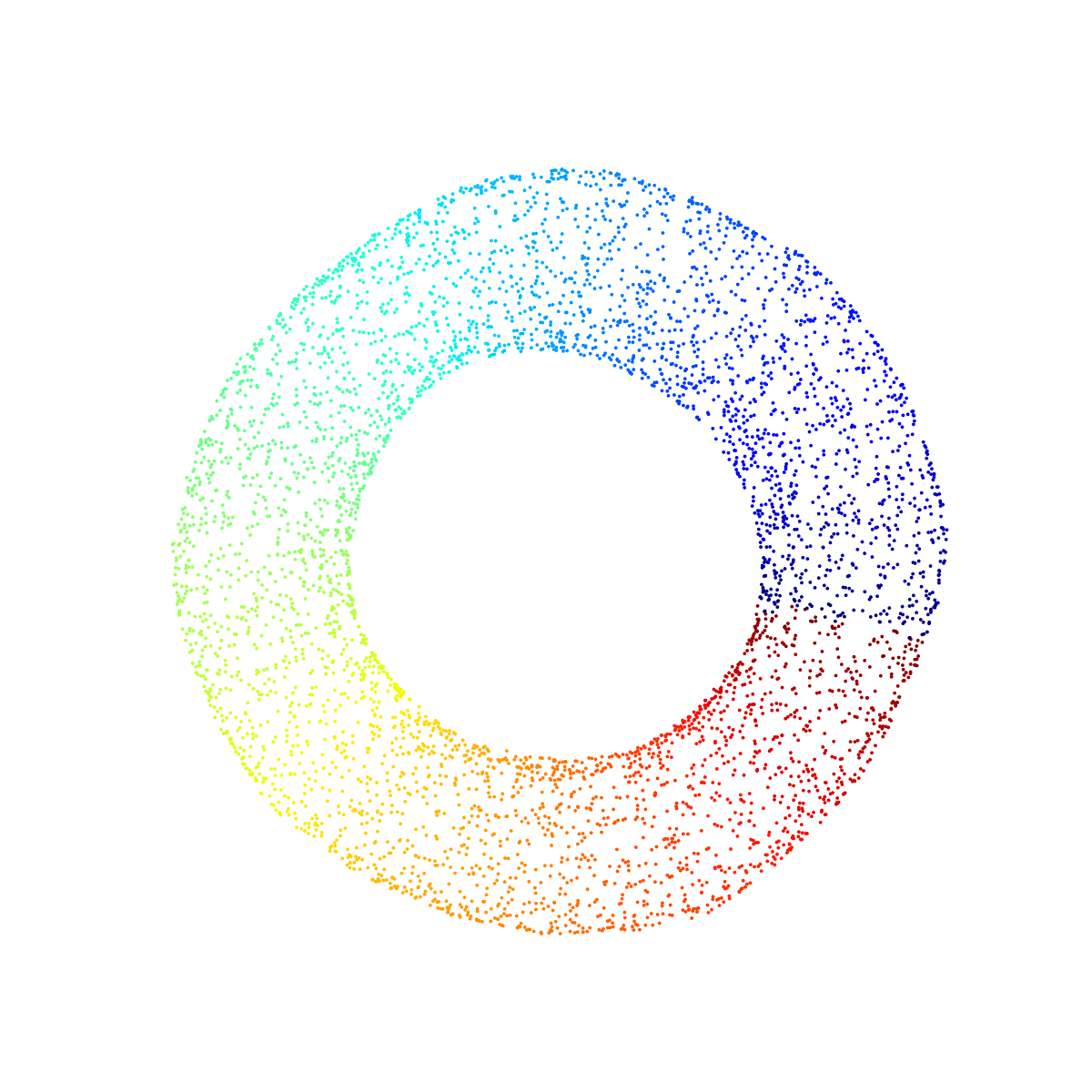} & \includegraphics[width=0.15\textwidth]{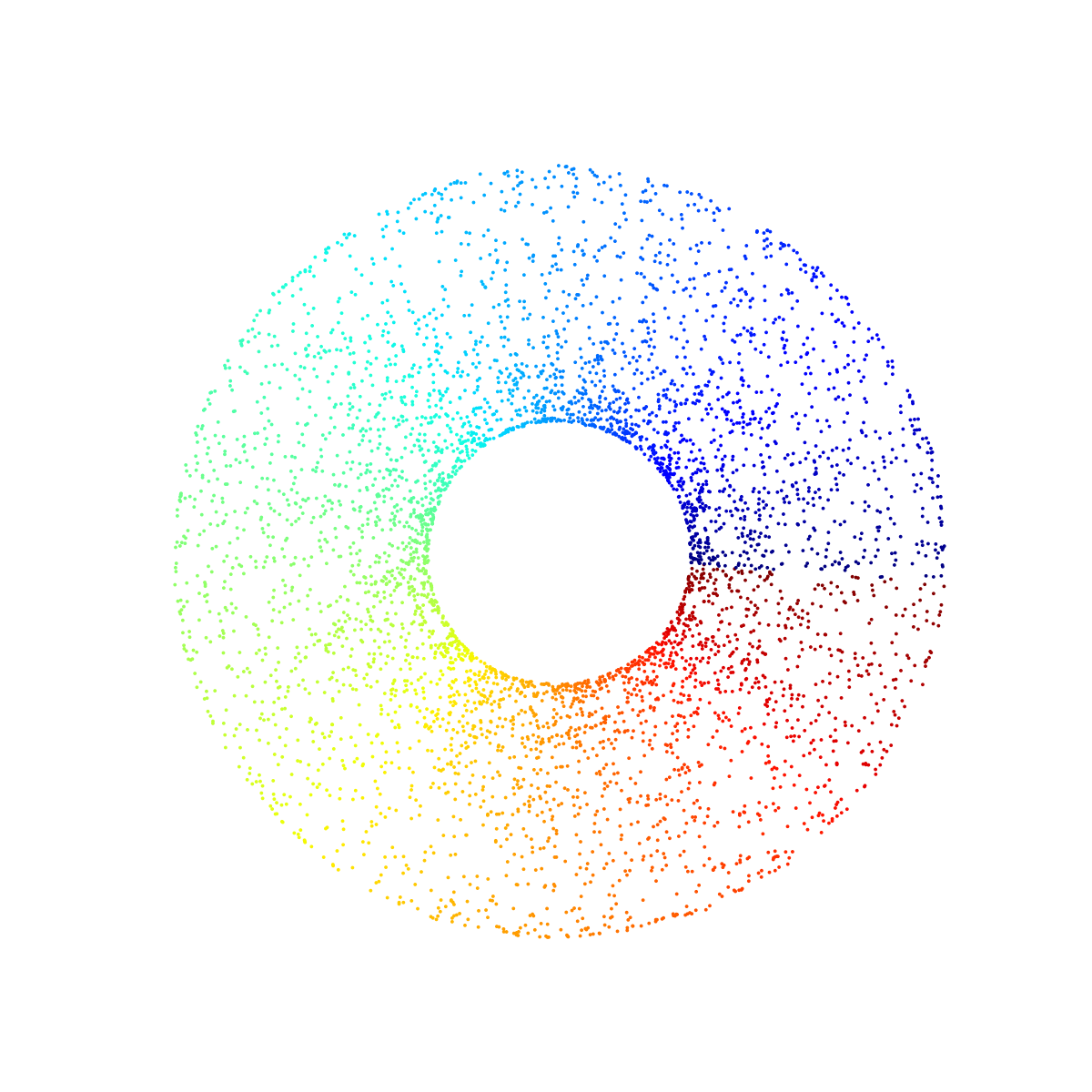} \\
    \hline
		\textbf{(3) Swiss role with hole} & \includegraphics[width=0.18\textwidth]{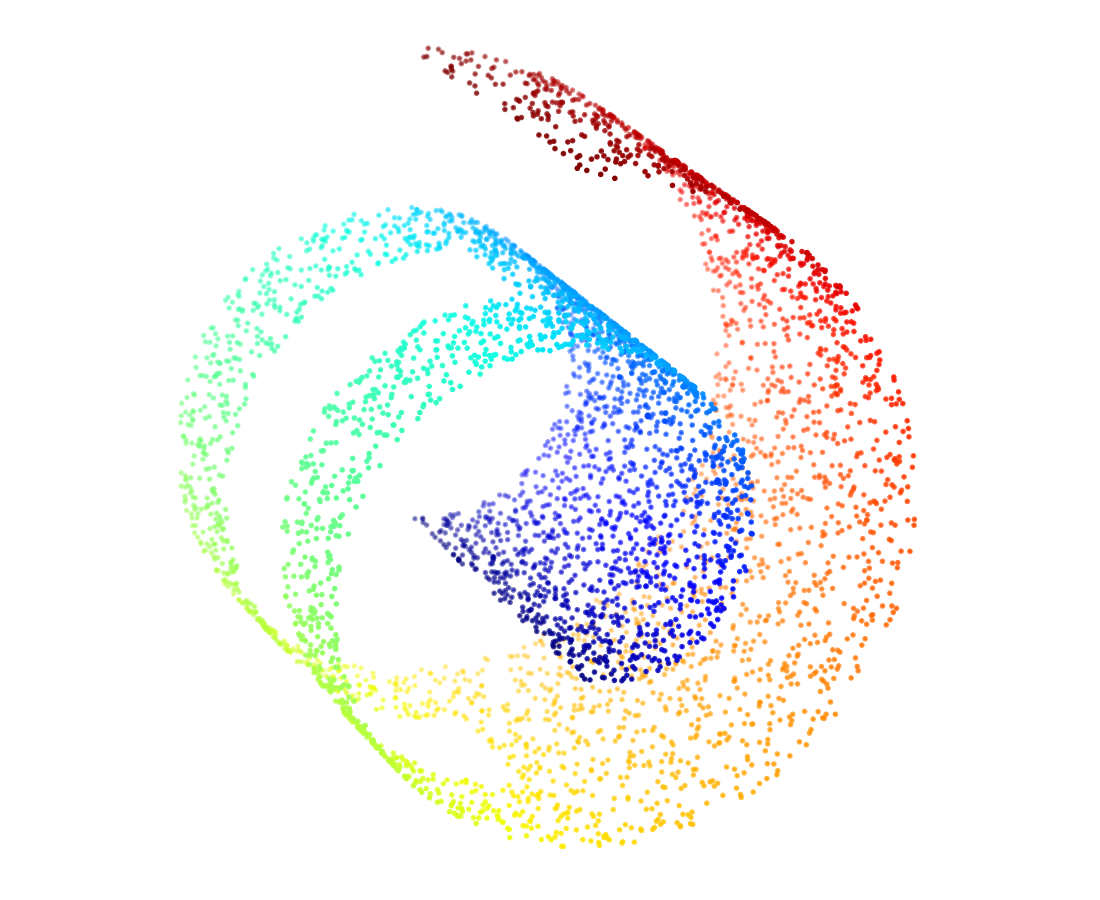} & \includegraphics[width=0.15\textwidth]{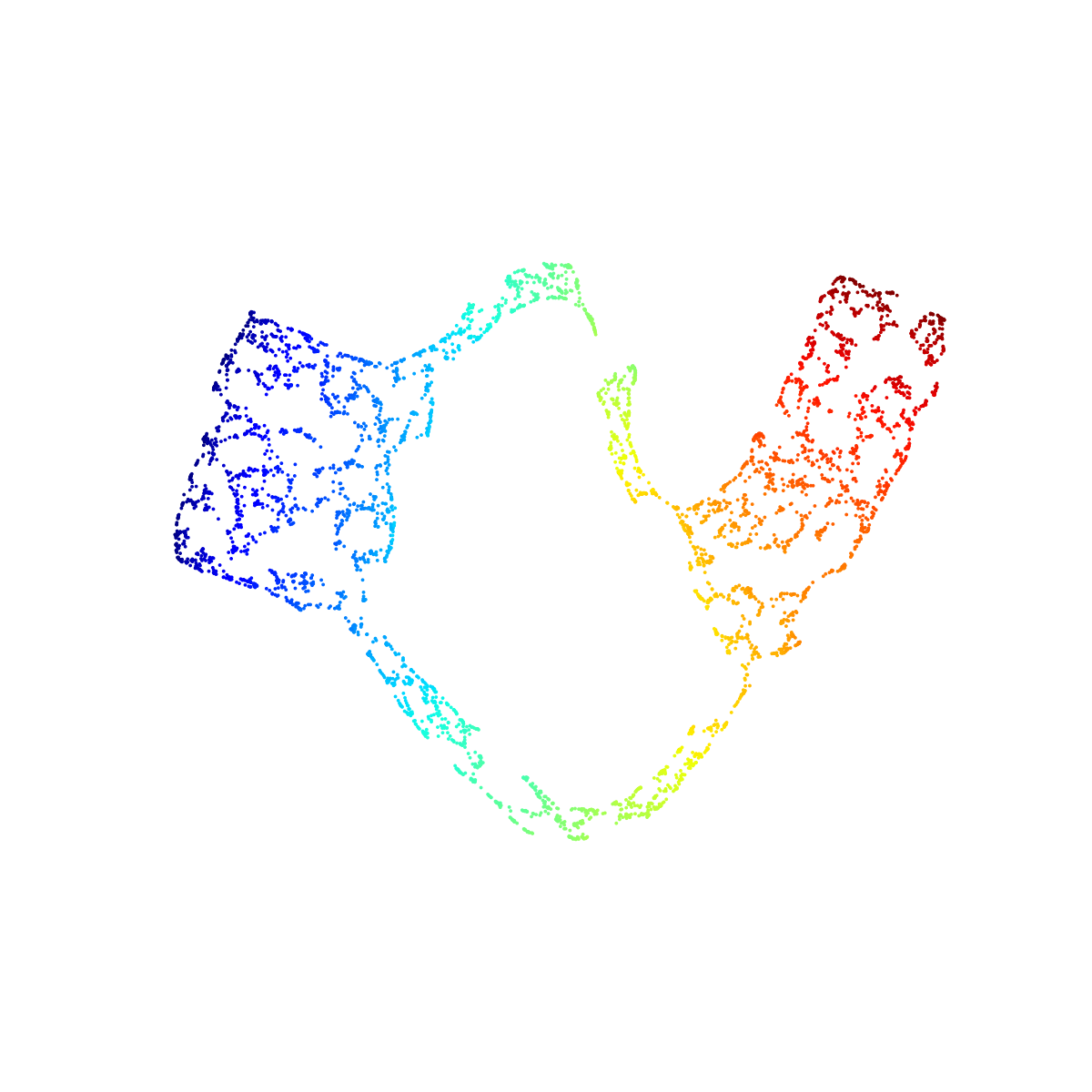} & \includegraphics[width=0.15\textwidth]{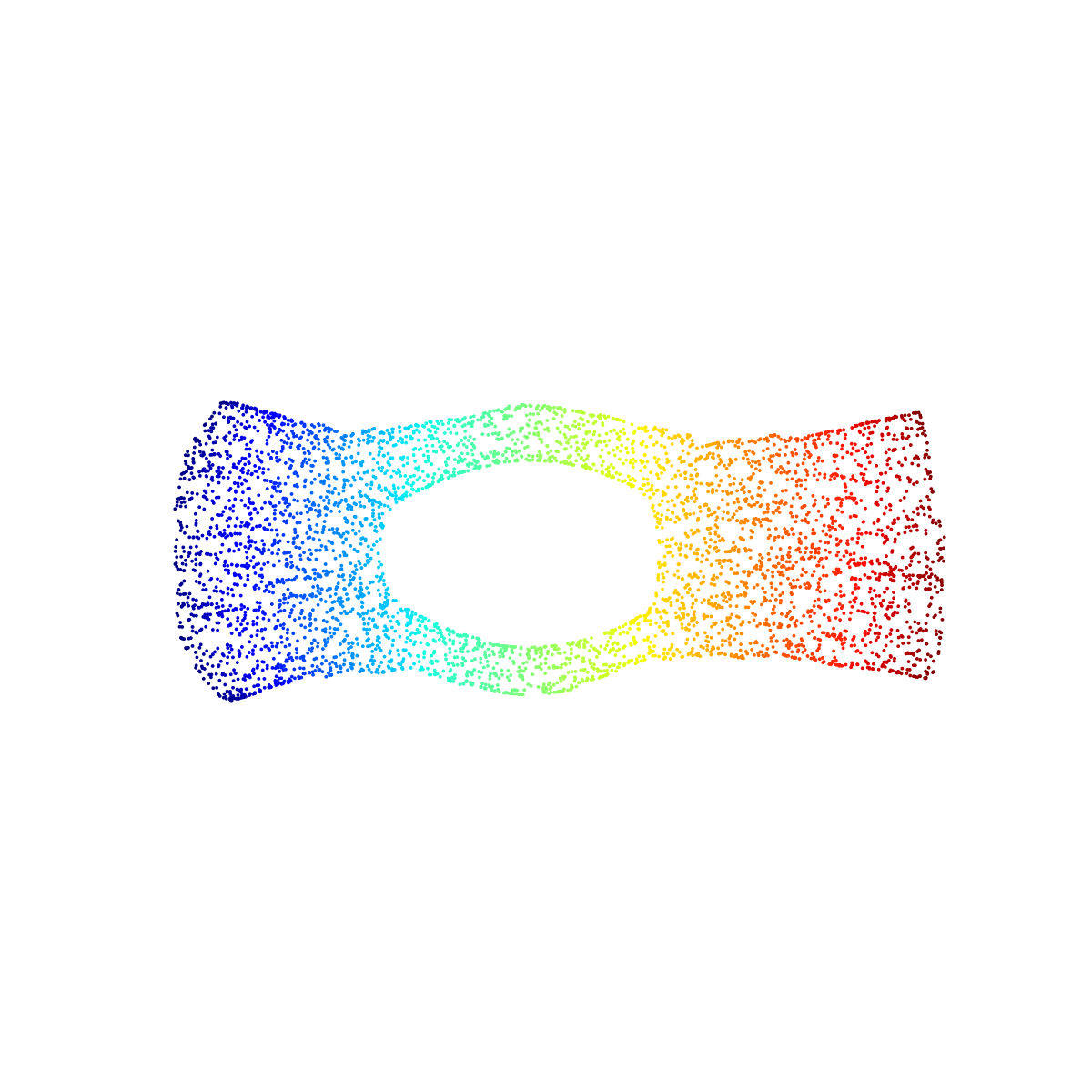} & \includegraphics[width=0.15\textwidth]{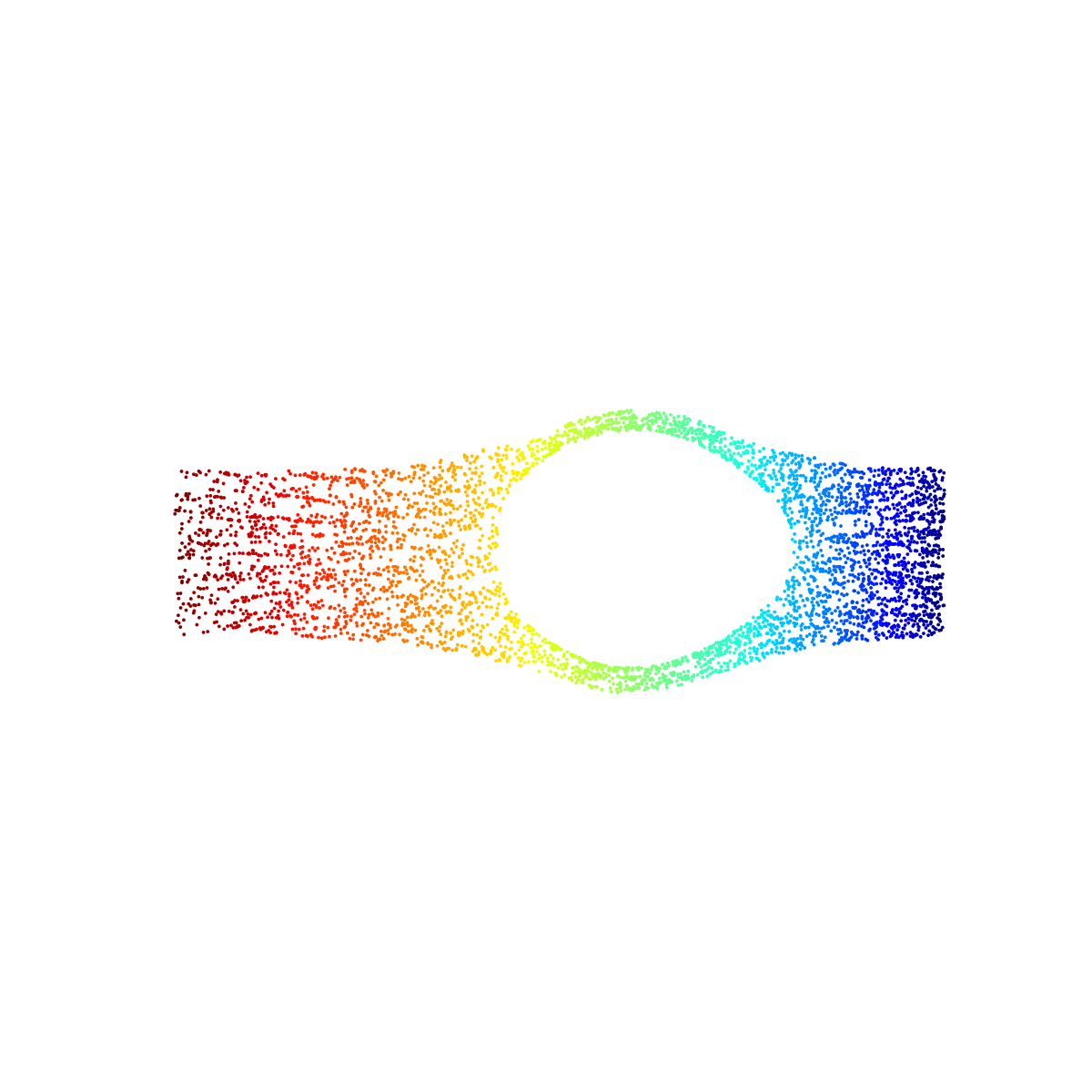} \\
		\hline 
	\end{tabular}
	\caption{Comparison of UMAP, IsUMap and Isomap on low-dimensional datasets (Hemisphere, Torus, Swiss role with hole).}
	\label{fig:comparisonUMAPandIsumap}
\end{table}

%\graphicspath{{figures/}}
\begin{table}[H]
	\centering
	\tiny
	\begin{tabular}{>{\centering\arraybackslash}m{1.2cm}|*{3}{>{\centering\arraybackslash}m{4cm}|}}
	  & \textbf{UMAP} & \textbf{IsUMap}  & \textbf{Isomap} \\
    \hline
		\textbf{(1) MNIST}
    & \includegraphics[width=0.3\textwidth]{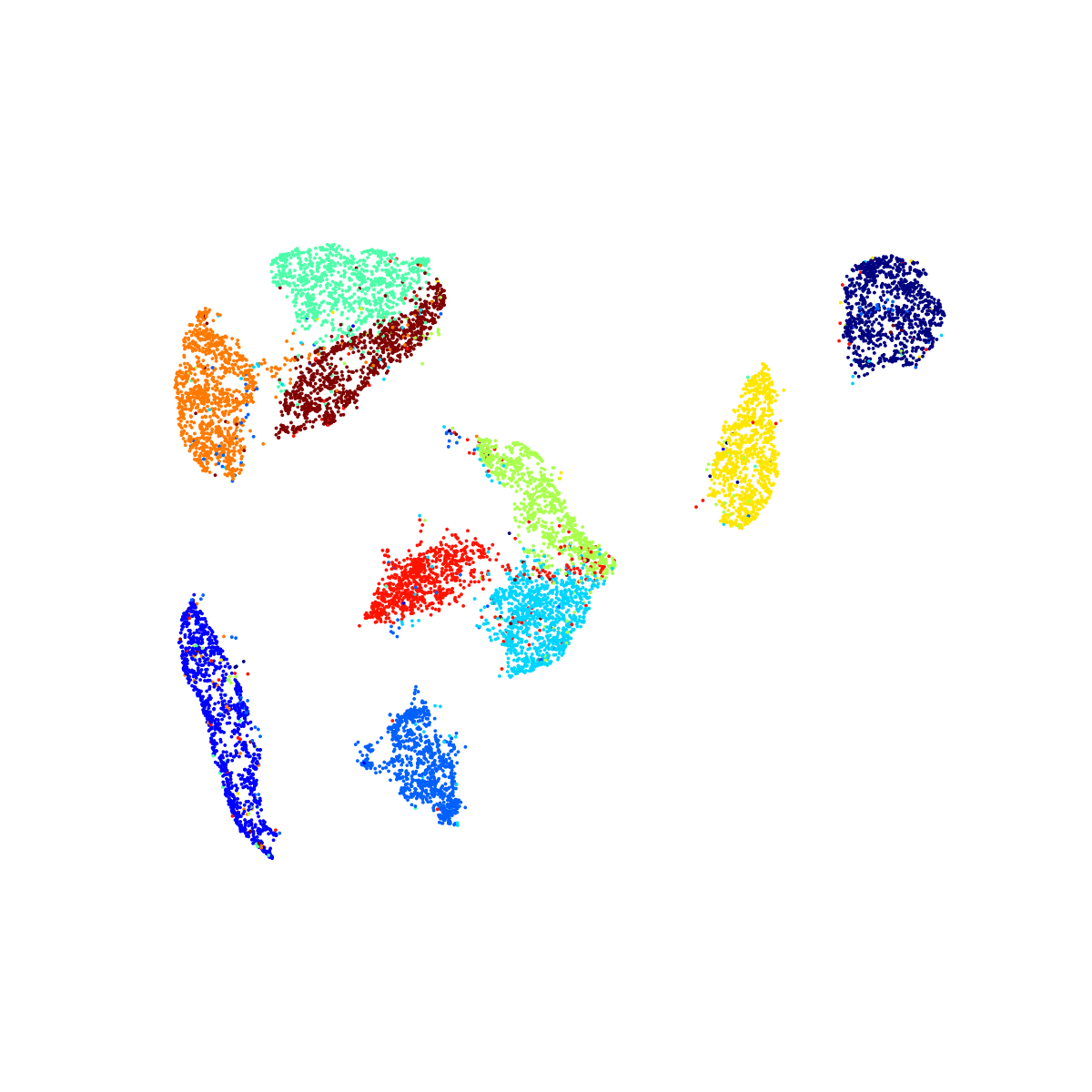} & \includegraphics[width=0.2\textwidth]{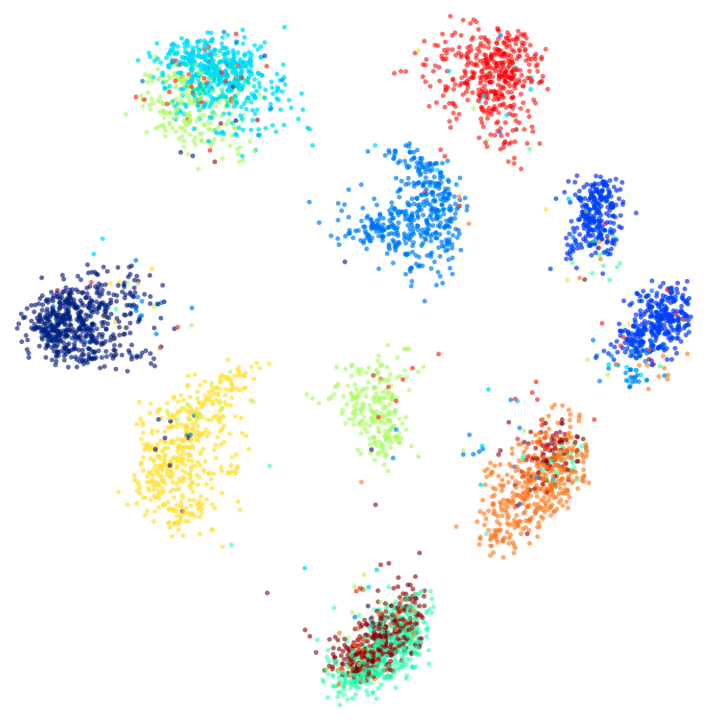} & \includegraphics[width=0.25\textwidth]{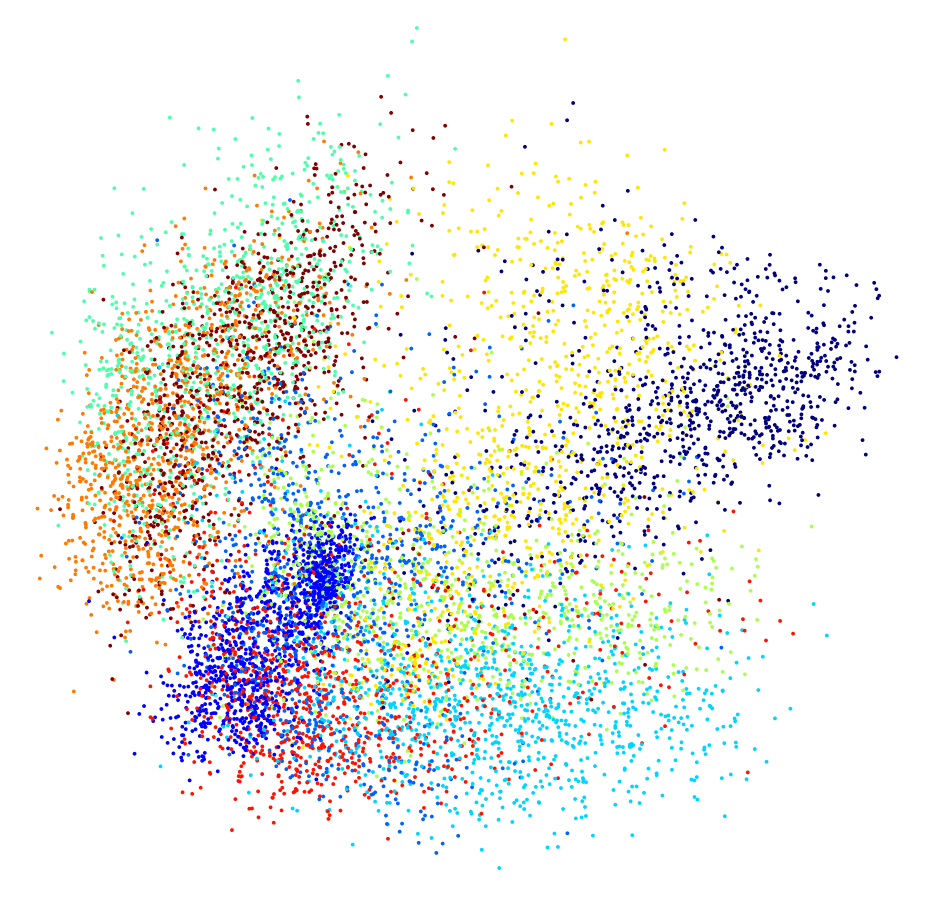}\\
    \hline
    \textbf{(2) Fashion-MNIST} & \includegraphics[width=0.28\textwidth]{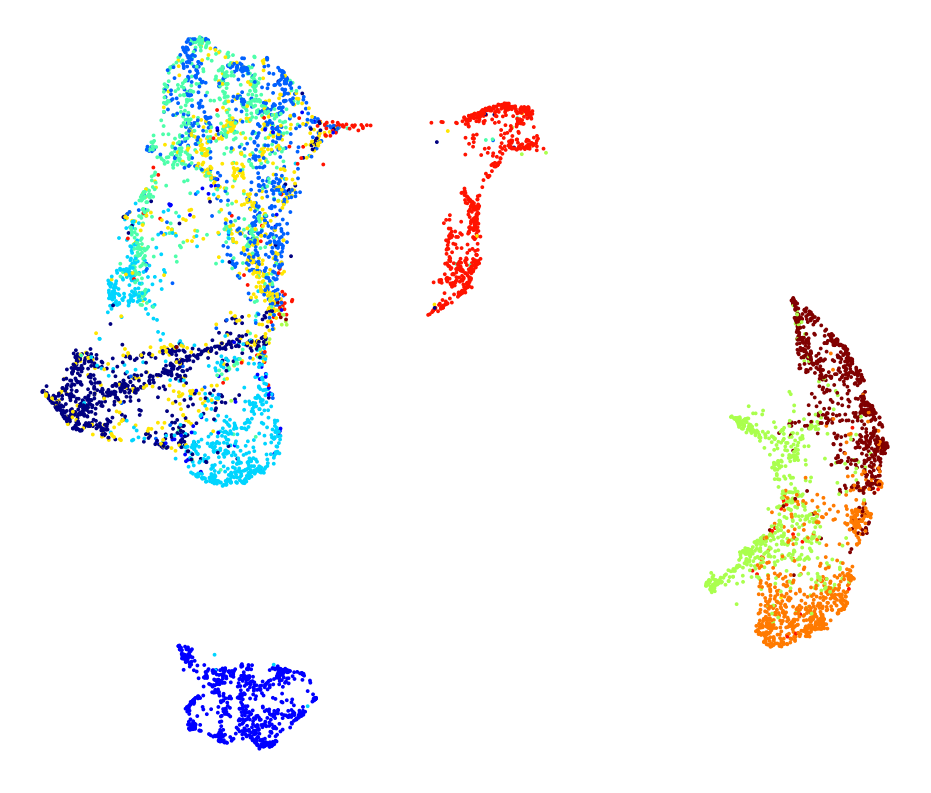} & \includegraphics[width=0.22\textwidth]{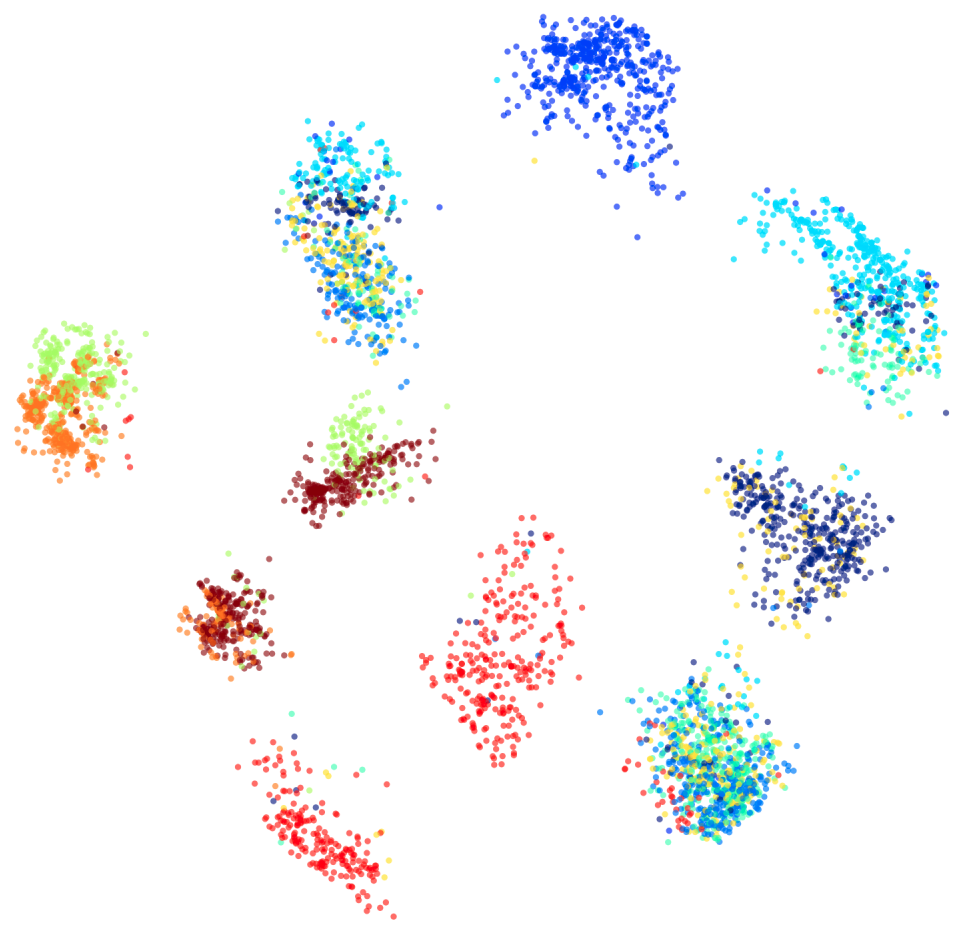} & \includegraphics[width=0.22\textwidth]{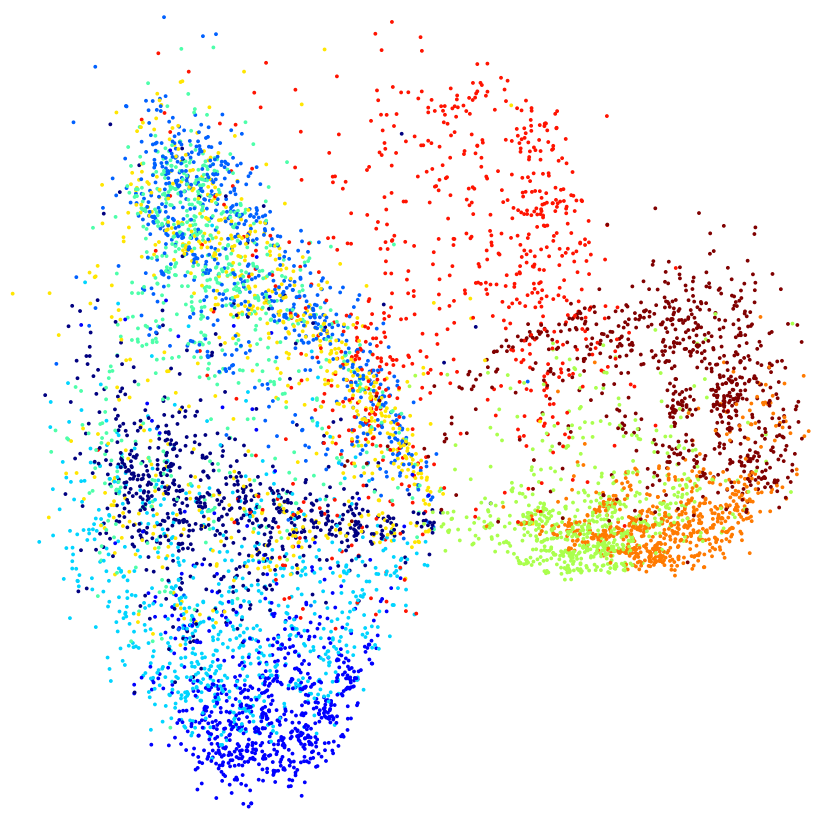}\\
		\hline  
	\end{tabular}
	\caption{Comparison of UMAP, IsUMap and Isomap on MNIST and FashionMNIST.}
	\label{fig:MNISTandFashionMNIST}
\end{table}


%\graphicspath{{figures/}}
\begin{table}[H]
	\centering
	\tiny
	\begin{tabular}{>{\centering\arraybackslash}m{1.2cm}|*{3}{>{\centering\arraybackslash}m{4cm}|}}
	  & \textbf{UMAP} & \textbf{IsUMap}  & \textbf{t-SNE} \\
    \hline
		\textbf{20-newsgroups}
    & \includegraphics[width=0.28\textwidth]{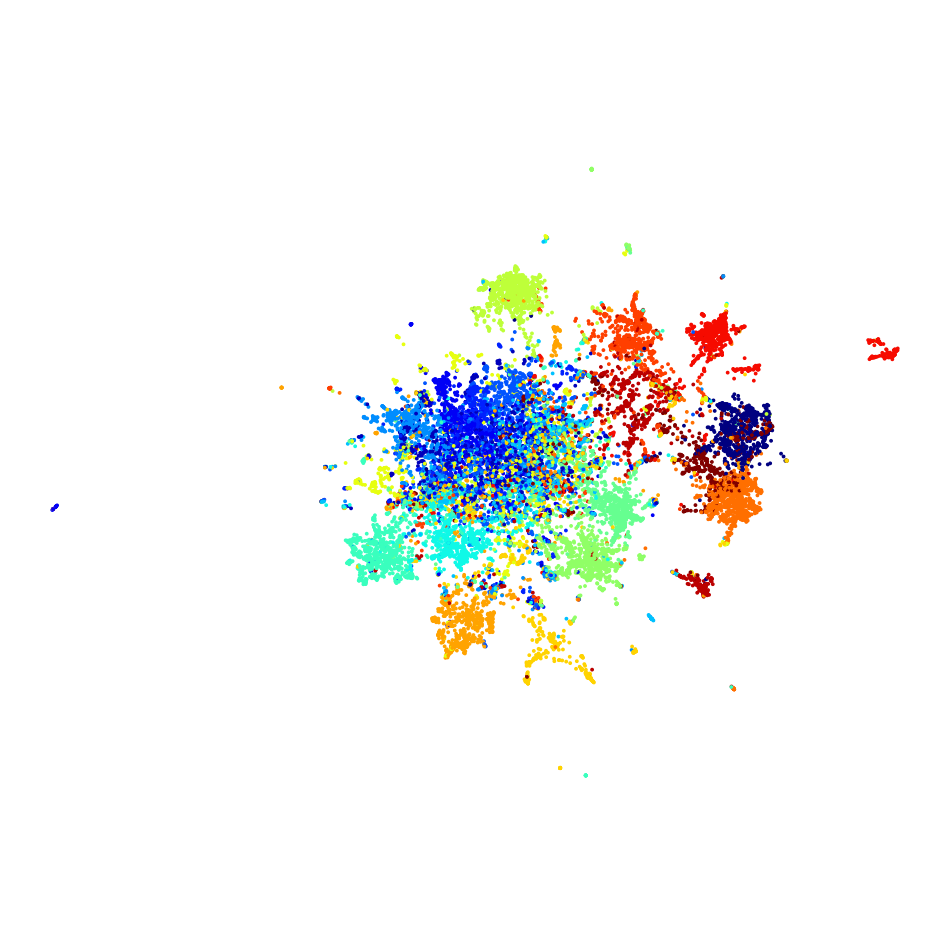} & \includegraphics[width=0.23\textwidth]{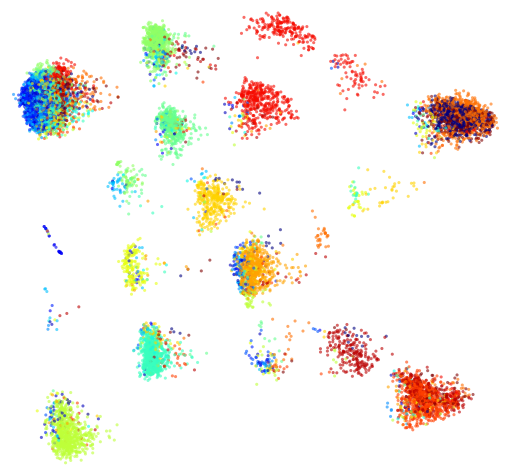} & \includegraphics[width=0.25\textwidth]{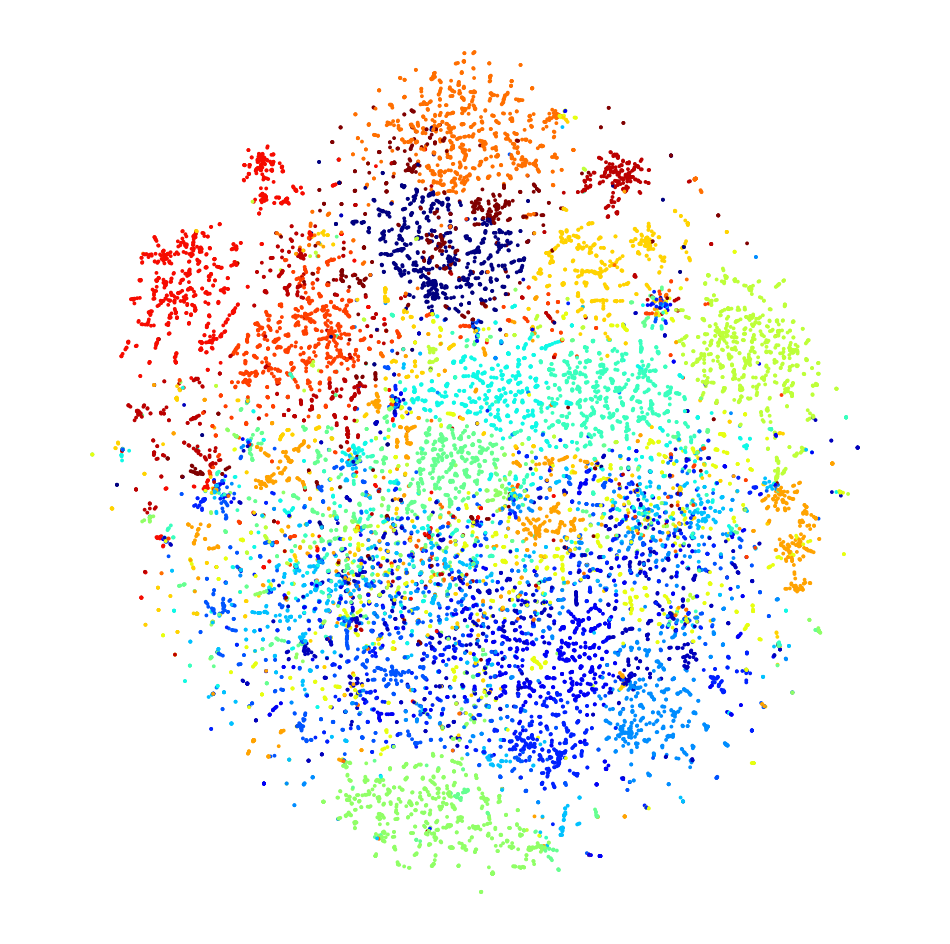}\\
		\hline 
	\end{tabular}
	\caption{Comparison of UMAP, IsUMap and t-SNE embeddings of the 20-newsgroups dataset.}
	\label{fig:newsgroups}
\end{table}

\begin{table}[H]
	\centering
	\tiny
	\begin{tabular}{>{\centering\arraybackslash}m{1.2cm}|*{2}{>{\centering\arraybackslash}m{6cm}|}}
	  & \textbf{UMAP} & \textbf{IsUMap}   \\
    \hline
		\textbf{BMMC}
    & \includegraphics[width=0.30\textwidth]{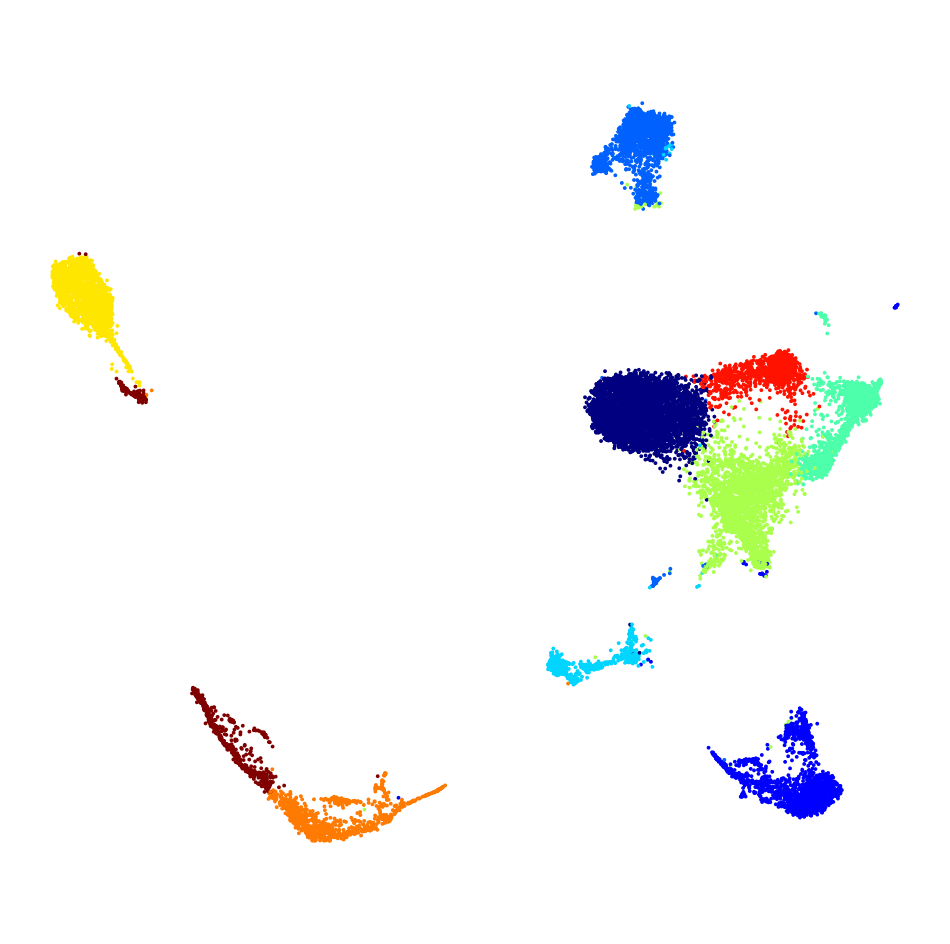} & \includegraphics[width=0.29\textwidth]{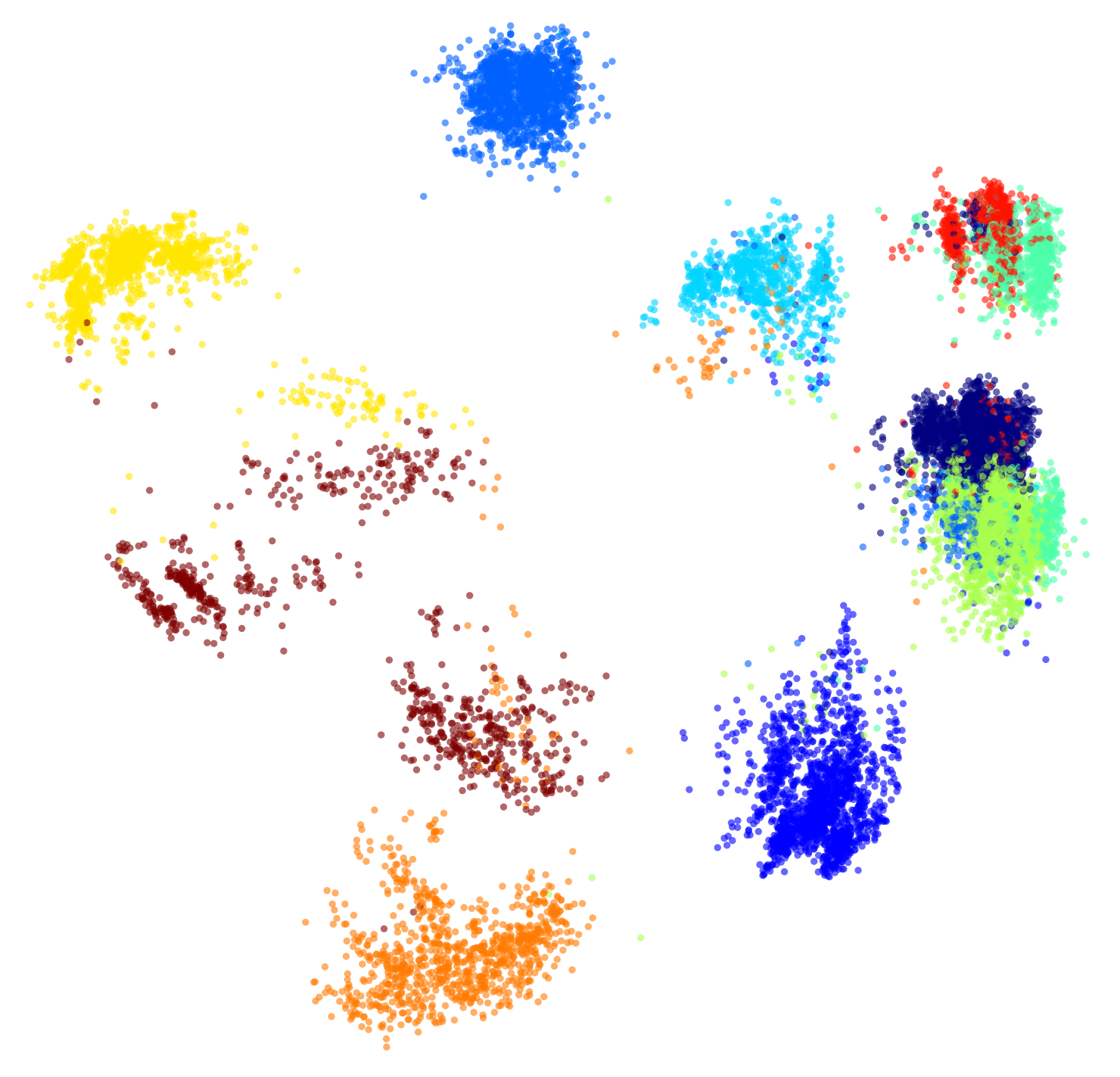}\\
		\hline 
    & \textbf{Isomap} & \textbf{t-SNE}  \\
    \hline
		\textbf{BMMC} 
    &\includegraphics[width=0.33\textwidth]{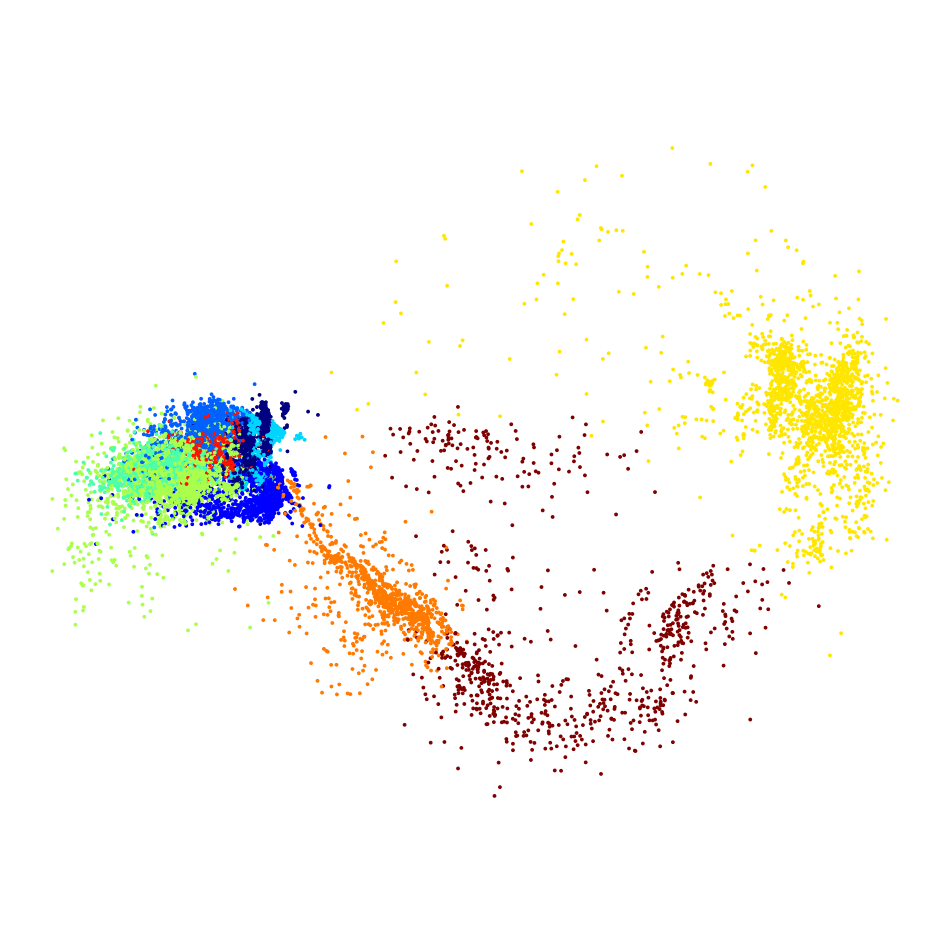} & \includegraphics[width=0.3\textwidth]{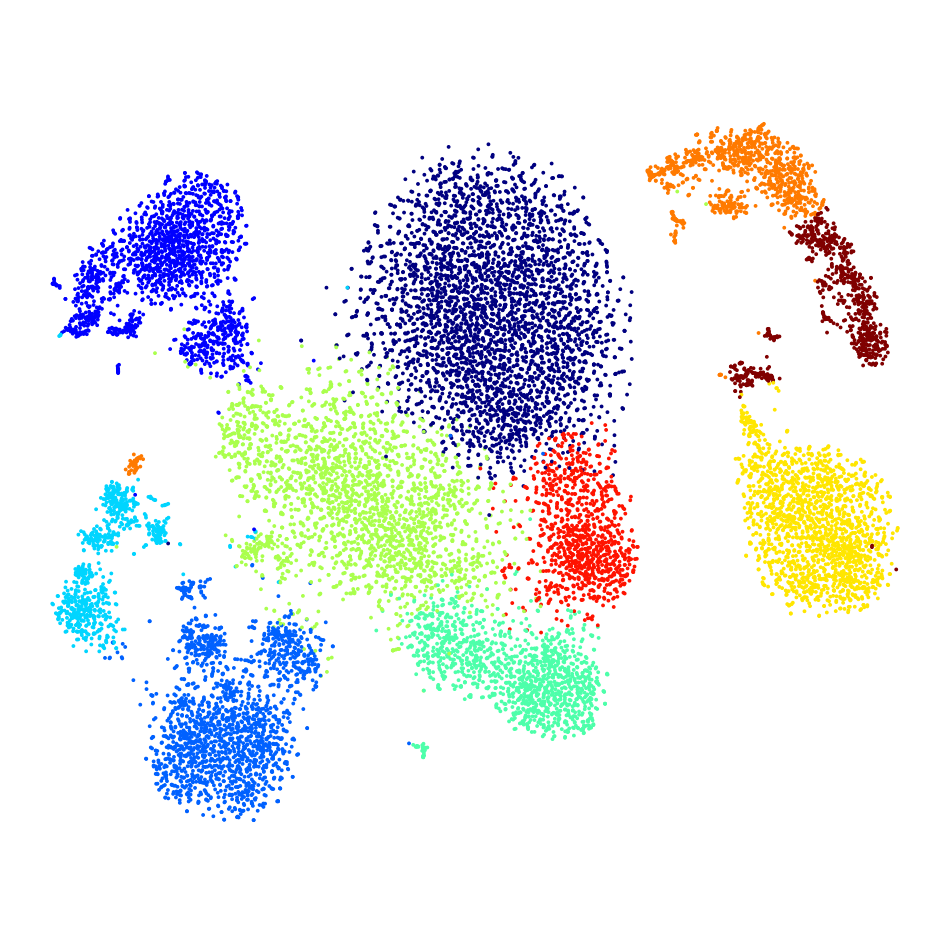}\\
		\hline 
	\end{tabular}
	\caption{Comparison of UMAP, IsUMap and t-SNE embeddings of the BMMC dataset.}
	\label{fig:BMMC}
\end{table}

%\graphicspath{{figures/}}
\begin{table}[H]
	\centering
	\tiny
	\begin{tabular}{>{\centering\arraybackslash}m{1.2cm}|*{3}{>{\centering\arraybackslash}m{4cm}|}}
	  & \textbf{UMAP} & \textbf{IsUMap} & \textbf{t-SNE} \\
    \hline
		\textbf{CRCC (Leiden)}
    & \includegraphics[width=0.27\textwidth]{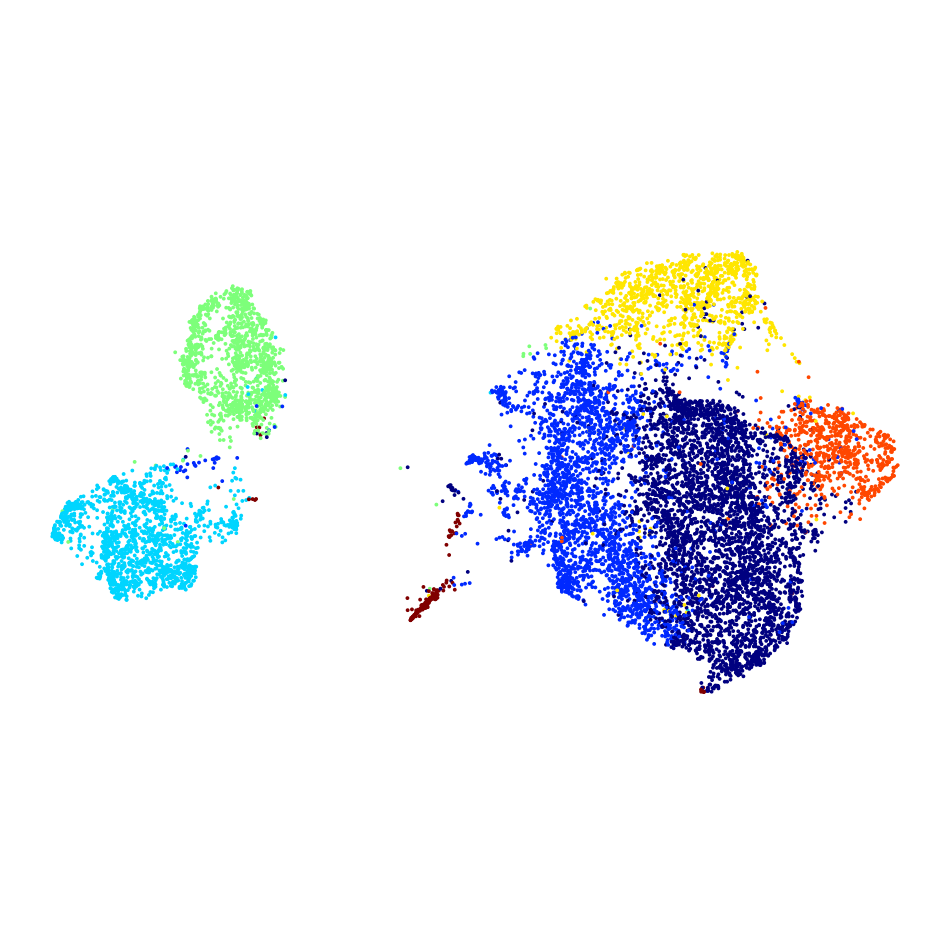} & \includegraphics[width=0.23\textwidth]{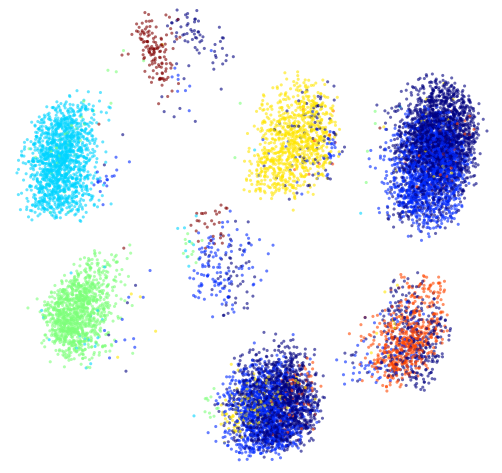} & \includegraphics[width=0.25\textwidth]{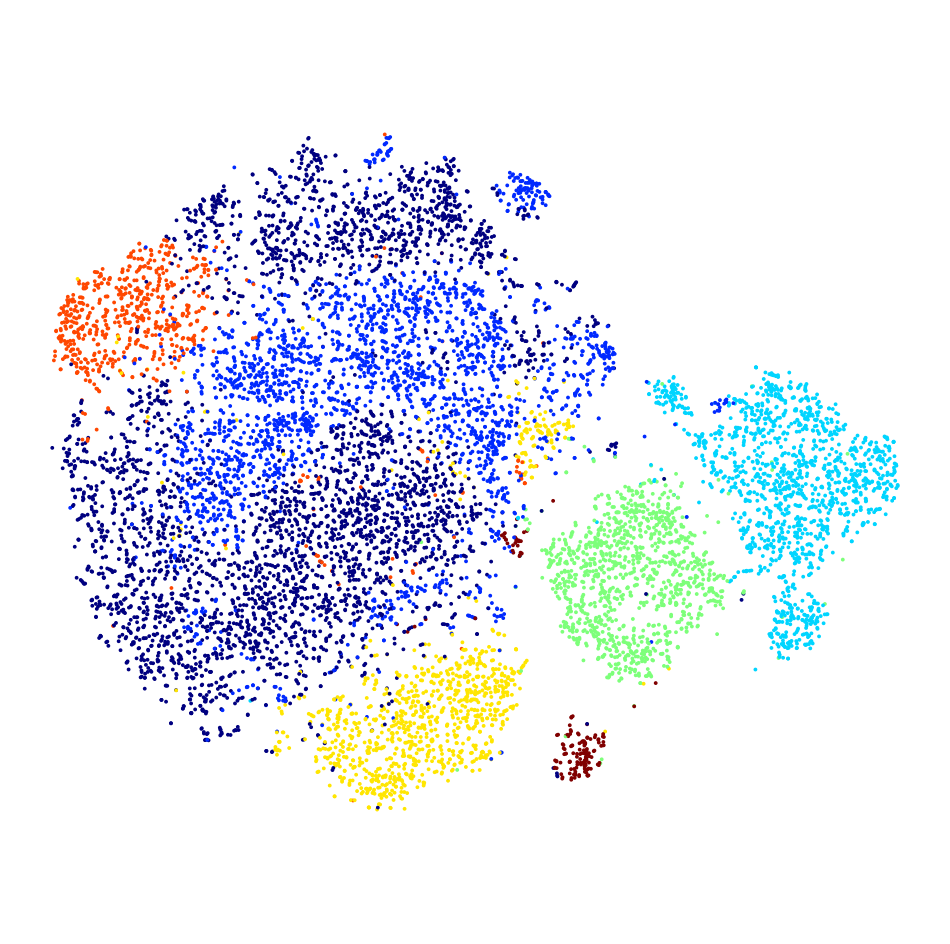}\\
    \hline
		\textbf{CRCC (SMOC2)}
    & \includegraphics[width=0.27\textwidth]{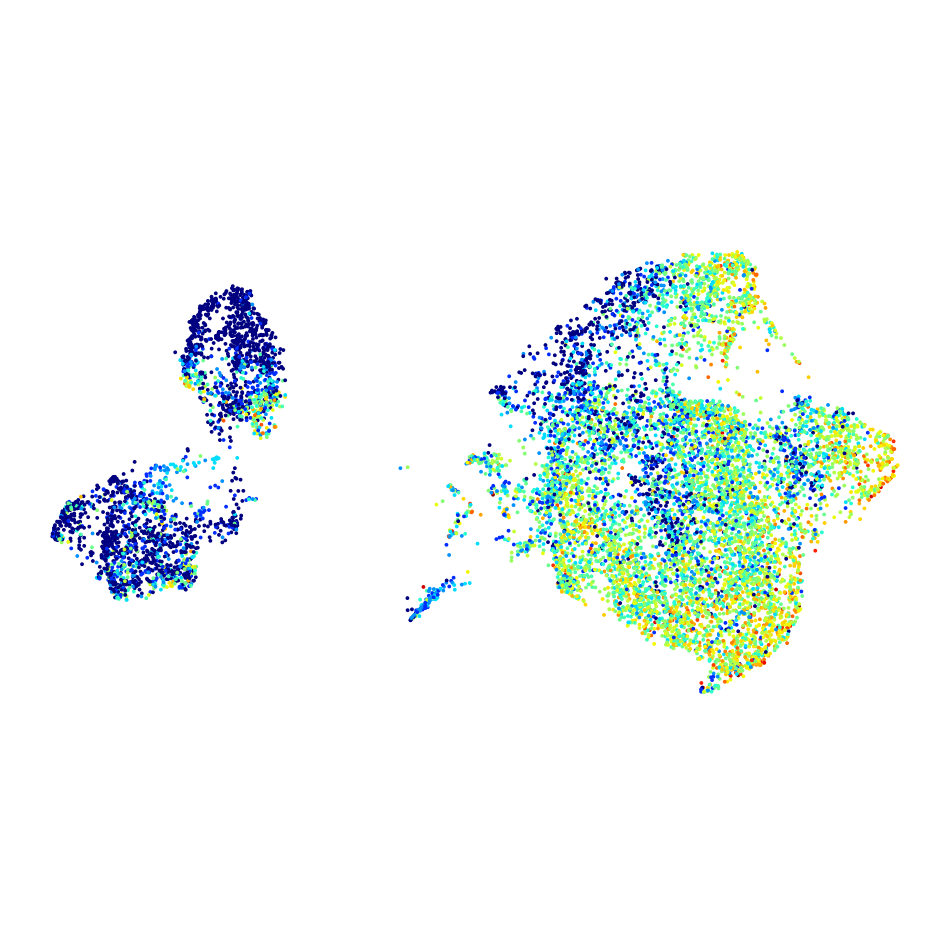} & \includegraphics[width=0.23\textwidth]{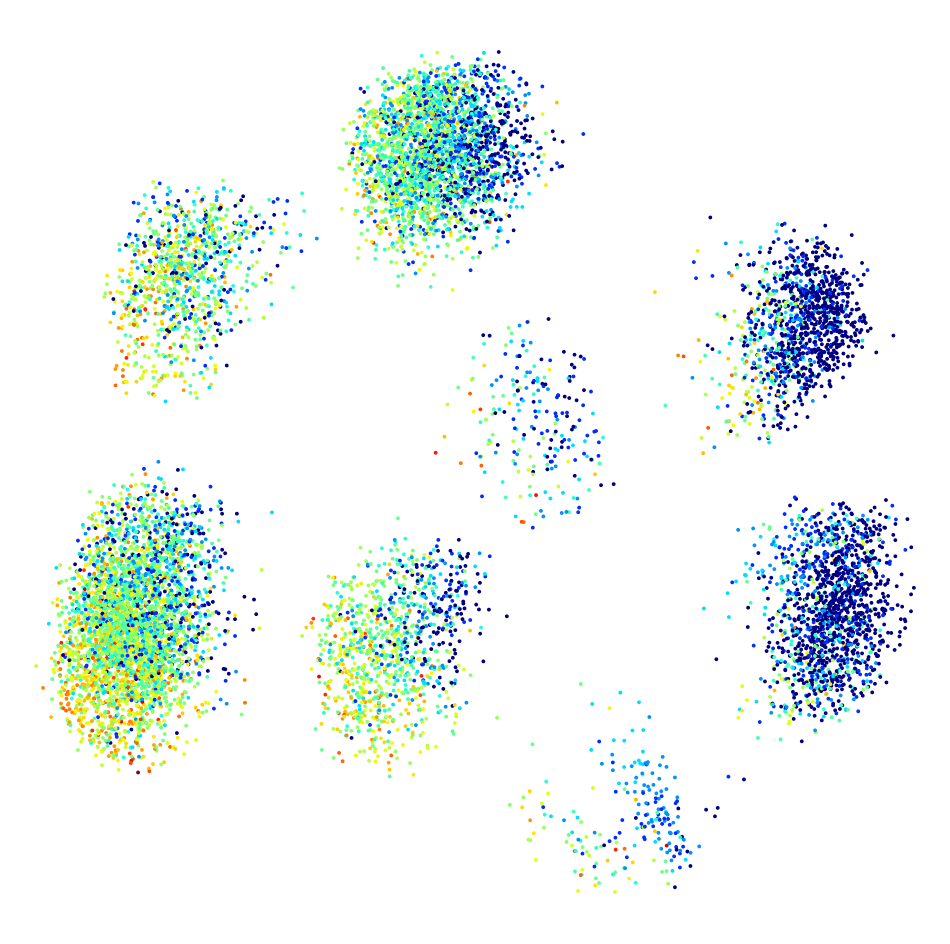} & \includegraphics[width=0.25\textwidth]{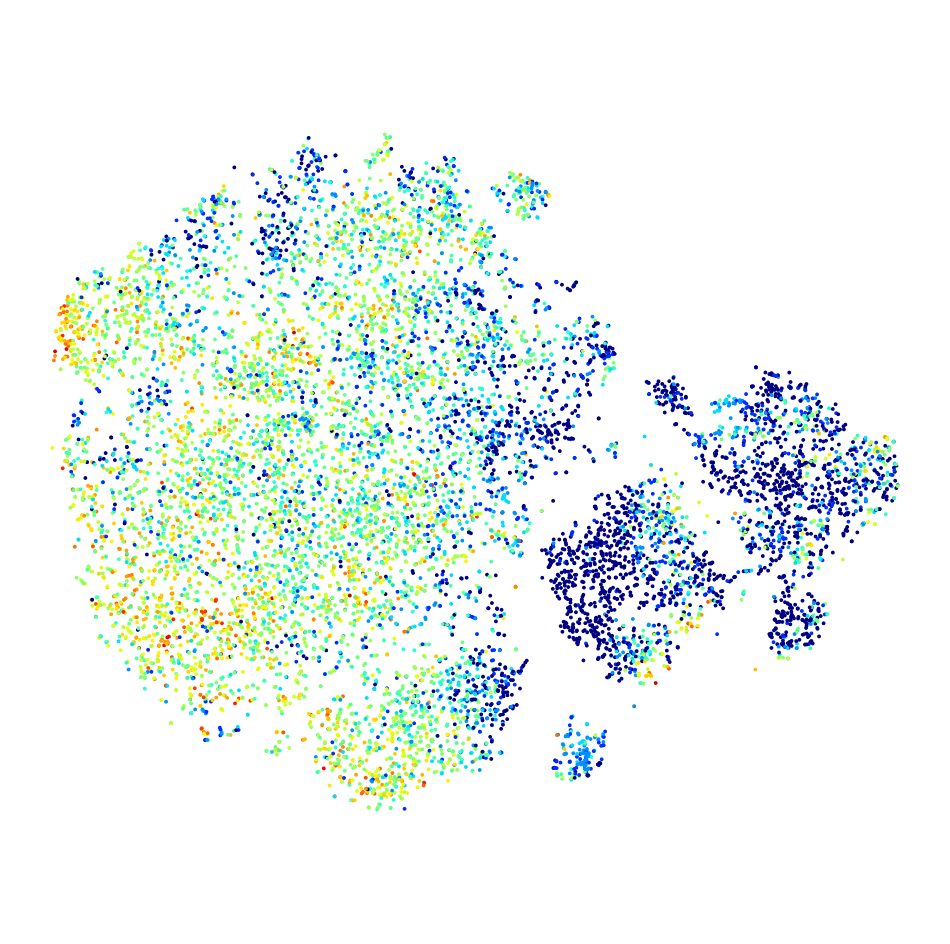}\\
    \hline
	\end{tabular}
	\caption{Comparison of UMAP, IsUMap (with Linkage clustering), Isomap and t-SNE embeddings of the CRCC scRNAseq dataset. Top row is colored by the labels of Leiden clustering and bottom row by the (log-transformed) SMOC2 gene expression count.}
	\label{fig:CRCC}
\end{table}

\subsection{Critical discussion}
\label{sec:effectsUMAP}

In light of Proposition \ref{prop:UMisInSfuz}, it is a priori surprising that the visual results of UMAP and plain Isomap differ so much. Since $\bo{UM}$ and $\bo{EPMet}$ are subcategories of $\bo{sFuz}$, only the embedding procedures differ, which represents the last step in the paths in Diagram \eqref{diag:umapAlternative}. We want to explain this difference in more detail below. There are two crucial things to point out in relation to the (force-directed graph layout) embedding used in UMAP (which have not been fully explained in the original reference \cite{McInnes18}):
\begin{enumerate}
  \item It is well-known that the \texttt{min-dist} parameter has a strong effect on the clustering abilities of UMAP. This parameter arises as follows: 
  The graph $H(Y)$ in \ref{U5} is not obtained by applying exactly the same operations to $Y$ as to $X$. Instead, the (Euclidean) distance $d_{\mathbb{R}^m}$ between the points in $Y$ is modified by subtracting a hyperparameter, called \texttt{min-dist} (substituting the nearest-neighbor distance $\rho_i$ in \eqref{eq:localDists}, while $\sigma_i$ is set to $1$), which makes the resulting metric $d_Y$ non-smooth, after which the Sing-functor $\te{Sing}^{\te{skeleton}}$ and the truncation $\te{ctr}_1$ is invoked to obtain the graph $H(Y)$ (no merging is necessary and symmetry is already ensured). 
  The edge weights on that graph correspond to the \texttt{min-dist}-modified distances, and are hence not smooth. Since smoothness is, however, required for applying gradient descent to $\mathcal{L}(Y)$ (defined in \ref{U5}), the weights $H(Y)_{ij}$ are approximated with the smooth function $\Phi(y_i,y_j):=(1+a(||y_j-y_k||^2)^b)^{-1}$ (where $a,b\in \mathbb{R}$ are chosen such that $\Phi(y_i,y_j)$ fits $H(Y)_{ij}$ as good as possible).  
  The authors of UMAP explain the necessity of $\rho_i$ and hence \texttt{min-dist} with a topological argument, namely local connectedness. However, even though subtraction of the distance of the nearest neighbor ensures that each $0$-simplex is the face
  of some $1$-simplex with membership strength $1$, this does not ensure local connectedness of the resulting graph because it might still have many disconnected components.

  We thus attribute the effect of \texttt{min-dist} mainly to the heavy tails of the approximation $\Phi$. Indeed, decreasing \texttt{min-dist} makes the tails of $\Phi$ more heavy.\footnote{The function $\Phi$ is, in fact, a slight generalization of the student-t-distribution employed in t-SNE (cf.~\cite{vandermaaten08a}), where both $a$ and $b$ are set to $1$. This distribution was employed in t-SNE precisely because it has heavier tails than the Gaussian distribution.}
  This generates a force that pushes the low-dimensional distances further apart, which is consistent with the observation that lower \texttt{min-dist} values increase the cluster distances.

  Note that, using the generalization of the adjunction that we developed in Section \ref{sec:Adjunction}, one can show that $\Phi$ gives itself rise to one such adjunction (because it is invertible, and the inverse can then be used to define $\te{Re}$, which can then be shown to be a functor) and can be understood as an alternative Sing-functor, that fits into the categorical scheme. Different Sing-functors can thus give rise to different visualizations.
  \item Even more importantly, the negative sampling strategy employed in UMAP has very strong effects on the embedding, as pointed out in \cite{Damrich2021}. In UMAP, the gradients $\Lambda_1$ and $\Lambda_2$ of the two terms $-\sum_{i,j} G_{ij}\log(\Phi_{ij})$ and $-\sum_{i,j}(1-G_{ij})\log(1-\Phi_{ij})$ of the loss $\mathcal{L}(Y)$ described in \ref{U5} (where we now replaced $H(Y)_{ij}$ by its ``approximation'' $\Phi_{ij}=\Phi(y_i,y_j)$) are computed separately with a sampling based approach to increase computational efficiency. For the first term, this means that only after sampling an edge from $i$ to $j$ with a probability proportional to $G_{ij}$, the gradient $\Lambda_1$ is computed for the distances involved. 
  The idea behind this is that most weights $G_{ij}$ in the adjacency matrix
  are $0$ because the graph was constructed as a merge of $k$
  nearest-neighborhood-star graphs, resulting in a  graph that is sparse in
  the sense that its adjacency matrix contains between $Nk$ and, if the t-conorm symmetrization process leads to an increase in the number of edges, at most $2Nk$ elements.
  For the first term, which equals $-\sum_{i,j} G_{ij}\log(\Phi_{ij})$, summands $G_{ij}\log(\Phi_{ij})$, where $G_{ij}$ is small, are themselves small anyways and their gradients, which are proportional to $G_{ij}$ as well, are also small and therefore one can employ this sampling based approach without distorting the gradients $\Lambda_1$ too much. 

  However, the same reasoning no longer works  for summands $(1-G_{ij})\log(1-\Phi_{ij})$ of the second term because $(1-G_{ij})$ is not small when $G_{ij}$ is small. One could sample edges for the second term with a probability $(1-G_{ij})$ instead but since, as explained above, most $G_{ij}$ are exactly equal to $0$, this would still result in the computation of a number of gradients $\Lambda_2$ that is proportional to $N^2$.
  In UMAP, they try to circumvent this problem by instead sampling (for each edge from $i$ to $j$ that was successfully sampled during the computation of $\Lambda_1$) an additional fixed number $m$ (where $m$ is a hyperparameter, also called `\texttt{negative\_sampling\_parameter}') of edges and computing $\Lambda_2$ for those. This actually distorts the gradient and under-represents repulsive forces as pointed out in \cite{Damrich2021} and in fact the distortion is so strong that they go as far as writing that the success that UMAP achieves in clustering tasks is to a large extent due to the effect of this undersampling, which is, however, not theoretically explained in the UMAP article. Numerically, we were able to confirm the results of \cite{Damrich2021} when re-implementing UMAP. 
\end{enumerate}

The effects described above, that were at least to some extent accidental side effects of computational approximations, seem to be very useful for clustering in practice (cf.~Table \ref{fig:comparisonUMAPandIsumap} (4,b)). However, the downside of them is that they do not have a theoretical foundation that allows to understand them in a completely predictable way. That implies that a visualization by UMAP, though often useful, might also display structures that are not inherent to the data, as can be observed in Table \ref{fig:comparisonUMAPandIsumap}, which can make their interpretation misleading.

We believe that IsUMap solves some of these issues. On the one hand, the construction of the geodesic distance matrix has a very solid theoretical grounding in Riemannian geometry and the theory of fuzzy simplicial sets that is described in the theoretical part of this work. On the other hand, the embedding procedure is very predictable and easy to understand. The first part is a simple metric embedding and the second part introduces a controlled deviation, induced by the clusters of the original space. The deviation from metric optimality can itself be visualized through the optimization path.
Furthermore, from a practical point of view, we find that the results above show that IsUMap in combination with the cluster separation optimization procedure often finds clusters that are more cleanly separated than those of UMAP, which are in turn often more cleanly separated than those of t-SNE. This can be of great help for classification when no labels are provided. 

Of course, our method also has limitations. Errors introduced during the clustering process propagate to the embedding stage and the assumption that the data lies on a sufficiently regular manifold might fail, for example. Furthermore, the cluster separation process might separate points that were formerly close in the geodesic metric space. However, in contrast to other dimension reduction methods, the visualization of the paths alleviates the last problem to a large extend because it provides knowledge about the error introduced by clustering. Finally, labels that are important for a given application might not correlate with metric distances, if the metric is not carefully chosen. The choice of the metric might itself have to be informed by the data.

\subsection{Outlook}

A possible avenue for further research from a category-theoretical point of view might be to combine the theory developed so far with ideas from functorial manifold learning and clustering as investigated in \cite{shiebler2020} and \cite{shiebler2020clustering}. Furthermore, it might be helpful to include the stochastic gradient descent procedure as well into the category-theoretical formulation, as initiated in \cite{fong_backprop_2019} and \cite{spivak_learners_2022} (which in turn is closely related to \cite{spivak_poly_book}).

For the reader interested in further empirical analysis and numerical experiments (though not yet including the cluster separation algorithm), we refer the reader to our more applied article \cite{Joharinad_Fahimi_Barth_Keck_Jost_2025}. Furthermore, our related book \cite{IsumapBook} presents additional simulations, expands on the ideas presented here and provides additional introductory material to address a broader readership.

\newpage
\section*{Declarations}
\label{sec:declarations}
\subsection*{Funding}

We acknowledge the support by "ScaDS.AI (Center for Scalable Data Analytics and Artificial Intelligence)", Dresden/Leipzig, and the "Max Planck Institute for Mathematics in the Sciences", Leipzig, and the "Max Planck School of Cognition", and the "Studienstiftung des deutschen Volkes e.V.", and the "Universitat Autònoma de Barcelona", Barcelona, Spain.

\subsection*{Competing interests}

The authors have no competing interests to declare that are relevant to the content of this article.

\subsection*{Ethics approval and consent to participate}

Not applicable.

\subsection*{Data and Code availability}

Please see \href{https://github.com/LUK4S-B/IsUMap}{https://github.com/LUK4S-B/IsUMap}.

\subsection*{Materials availability}

Not applicable.

\subsection*{Author contributions}

The article was written in the context of a bigger project about
geometric and category-theoretical methods in data analysis. In
particular, regular meetings among the authors Lukas Barth, Fatemeh
(Hannaneh) Fahimi, Parvaneh Joharinad, Jürgen Jost and Janis Keck
(listed in alphabetic order) took place, in which all ideas were discussed. 

In this project, Lukas Barth developed the category-theoretical part and
its corresponding theorems. Furthermore, he also developed and implemented the cluster separation optimization method, described in Section \ref{sec:embedding} and item 5 of Section \ref{sec:algorithmicDescription}, and conducted the corresponding simulations.
This article was mostly written by him. He also contributed
to the numerical implementation of IsUMap, testing and implementing many variations. Finally, he contributed to a detailed analysis of UMAP. Other contributions are presented in \cite{IsumapBook} and \cite{VRpaper}. 

Fatemeh (Hannaneh) Fahimi contributed to the project by providing
geometric insights and using IsUMap for many simulations of synthetic
and empirical datasets. She compared our method to other manifold learning
techniques and extended the algorithm to accommodate different merging processes. Further details of those contributions are presented in the publications
\cite{VRpaper} and \cite{IsumapBook}. 

Janis Keck made contributions to the geometric understanding and
provided statistical and probability theoretical insights, which are
also presented in \cite{IsumapBook}. Furthermore, he conceived and implemented
the first versions of the IsUMap algorithm and also contributed to later
implementations. 

Parvaneh Joharinad provided a formulation of the method based on
Vietoris-Rips filtrations, and contributed to the homological and
geometric part of our work, detailed in \cite{VRpaper} and \cite{IsumapBook}.
Furthermore, she organized the meetings. 

Jürgen Jost guided the work as supervisor, provided geometric ideas and
valuable comments. He also took part in writing, especially in \cite{IsumapBook}.

Thomas Jan Mikhail was often available during the project to discuss category-theoretical questions and contributed to this article by providing valuable feedback, proofreading everything, correcting
definitions and improving and contributing to the proofs.

\subsection*{Revision history}

\textbf{Second version.}
In the second version, we improved the theoretical exposition and clarity of several proofs and added a more explicit characterization of the essential image of the Sing functor.\\
Moreover, the cluster separation optimization algorithm was invented, implemented and evaluated on a variety of datasets, which strongly improved the experimental results as described in Section \ref{sec:application}.

\textbf{Third version.}
In the second version some proofs still missed the explicit verification of the injectivity and gluing conditions for fuzzy (simplicial) sets. In this version, we completed the proofs with these verifications. While doing so, we also derived a more useful description for the skeleton functor in Prop.~\ref{prop:leftRightAdjointOfTruncation}.\\
Furthermore, we thank David Wegmann for pointing out in \cite{wegmann2026theoryumap} that we carried over two issues from \cite{Spivak09} occurring in the definitions of fuzzy sets and of the metric realization. In this version, we corrected these problems by adapting Def.~\ref{def:fuzzySet} and Equation~\eqref{eq:spivakReObjects}.

\bibliography{bib}

@misc{mullner2011linkage,
      title={Modern hierarchical, agglomerative clustering algorithms},
      author={Daniel Müllner},
      year={2011},
      eprint={1109.2378},
      archivePrefix={arXiv},
      primaryClass={stat.ML},
}

@article{Linkage2001,
    author = {Bar-Joseph, Ziv and Gifford, David K. and Jaakkola, Tommi S.},
    title = {Fast optimal leaf ordering for hierarchical clustering},
    journal = {Bioinformatics},
    volume = {17},
    number = {suppl 1},
    pages = {S22-S29},
    year = {2001},
    month = {06},
}

@article{zadeh1965fuzzy,
  title={Fuzzy sets},
  author={Zadeh, Lotfi A},
  journal={Information and control},
  volume={8},
  number={3},
  pages={338--353},
  year={1965},
  publisher={Elsevier}
}

@book{awodey10,
  title={Category theory},
  author={Awodey, Steve},
  volume={52},
  year={2010},
  publisher={OUP Oxford}
}

@article{barr1986fuzzy,
  title={Fuzzy set theory and topos theory},
  author={Barr, Michael},
  journal={Canadian Mathematical Bulletin},
  volume={29},
  number={4},
  pages={501--508},
  year={1986},
  publisher={Cambridge University Press}
}

@misc{wegmann2026theoryumap,
      title={The Theory behind UMAP?}, 
      author={David Wegmann},
      year={2026},
      eprint={2603.03375},
      archivePrefix={arXiv},
      primaryClass={stat.ML},
      url={https://arxiv.org/abs/2603.03375}, 
}

@ARTICLE{2020SciPy-NMeth,
  author  = {Virtanen, Pauli and Gommers, Ralf and Oliphant, Travis E. and
            Haberland, Matt and Reddy, Tyler and Cournapeau, David and
            Burovski, Evgeni and Peterson, Pearu and Weckesser, Warren and
            Bright, Jonathan and {van der Walt}, St{\'e}fan J. and
            Brett, Matthew and Wilson, Joshua and Millman, K. Jarrod and
            Mayorov, Nikolay and Nelson, Andrew R. J. and Jones, Eric and
            Kern, Robert and Larson, Eric and Carey, C J and
            Polat, {\.I}lhan and Feng, Yu and Moore, Eric W. and
            {VanderPlas}, Jake and Laxalde, Denis and Perktold, Josef and
            Cimrman, Robert and Henriksen, Ian and Quintero, E. A. and
            Harris, Charles R. and Archibald, Anne M. and
            Ribeiro, Ant{\^o}nio H. and Pedregosa, Fabian and
            {van Mulbregt}, Paul and {SciPy 1.0 Contributors}},
  title   = {{{SciPy} 1.0: Fundamental Algorithms for Scientific
            Computing in Python}},
  journal = {Nature Methods},
  year    = {2020},
  volume  = {17},
  pages   = {261--272},
  adsurl  = {https://rdcu.be/b08Wh},
  doi     = {10.1038/s41592-019-0686-2},
}

@article{scikit-learn,
 title={Scikit-learn: Machine Learning in {P}ython},
 author={Pedregosa, F. and Varoquaux, G. and Gramfort, A. and Michel, V.
         and Thirion, B. and Grisel, O. and Blondel, M. and Prettenhofer, P.
         and Weiss, R. and Dubourg, V. and Vanderplas, J. and Passos, A. and
         Cournapeau, D. and Brucher, M. and Perrot, M. and Duchesnay, E.},
 journal={Journal of Machine Learning Research},
 volume={12},
 pages={2825--2830},
 year={2011}
}

@article{pytorch2019,
  author       = {Adam Paszke and
                  Sam Gross and
                  Francisco Massa and
                  Adam Lerer and
                  James Bradbury and
                  Gregory Chanan and
                  Trevor Killeen and
                  Zeming Lin and
                  Natalia Gimelshein and
                  Luca Antiga and
                  Alban Desmaison and
                  Andreas K{\"{o}}pf and
                  Edward Z. Yang and
                  Zach DeVito and
                  Martin Raison and
                  Alykhan Tejani and
                  Sasank Chilamkurthy and
                  Benoit Steiner and
                  Lu Fang and
                  Junjie Bai and
                  Soumith Chintala},
  title        = {PyTorch: An Imperative Style, High-Performance Deep Learning Library},
  journal      = {CoRR},
  volume       = {abs/1912.01703},
  year         = {2019},
  url          = {http://arxiv.org/abs/1912.01703},
  eprinttype    = {arXiv},
  eprint       = {1912.01703},
  bibsource    = {dblp computer science bibliography, https://dblp.org}
}

@misc{scanpyTutorial,
  title = {Scanpy Examples},
  author = {Scanpy},
  howpublished = {\url{https://scanpy.readthedocs.io/en/stable/tutorials/basics/clustering.html}},
  year = {2025}
}

@article{satija2015Seurat,
  title={Spatial reconstruction of single-cell gene expression data},
  author={Satija, Rahul and Farrell, Jeffrey A and Gennert, David and Schier, Alexander F and Regev, Aviv},
  journal={Nature biotechnology},
  volume={33},
  number={5},
  pages={495--502},
  year={2015},
  publisher={Nature Publishing Group US New York}
}

@article{stuart2019Seurat,
  title={Comprehensive integration of single-cell data},
  author={Stuart, Tim and Butler, Andrew and Hoffman, Paul and Hafemeister, Christoph and Papalexi, Efthymia and Mauck, William M and Hao, Yuhan and Stoeckius, Marlon and Smibert, Peter and Satija, Rahul},
  journal={cell},
  volume={177},
  number={7},
  pages={1888--1902},
  year={2019},
  publisher={Elsevier}
}

@article{Traag2019Leiden,
author={Traag, V. A.
and Waltman, L.
and van Eck, N. J.},
title={From Louvain to Leiden: guaranteeing well-connected communities},
journal={Scientific Reports},
year={2019},
month={Mar},
day={26},
volume={9},
number={1},
pages={5233},
issn={2045-2322},
doi={10.1038/s41598-019-41695-z},
url={https://doi.org/10.1038/s41598-019-41695-z}
}

@article{wolf2018scanpy,
  title={SCANPY: large-scale single-cell gene expression data analysis},
  author={Wolf, F Alexander and Angerer, Philipp and Theis, Fabian J},
  journal={Genome biology},
  volume={19},
  number={1},
  pages={15},
  year={2018},
  publisher={Springer}
}

@article{johnson1984extensions,
  title={Extensions of Lipschitz mappings into a Hilbert space},
  author={Johnson, William B and Lindenstrauss, Joram and others},
  journal={Contemporary mathematics},
  volume={26},
  number={189-206},
  pages={1},
  year={1984}
}

@article{nash1956imbedding,
  title={The imbedding problem for Riemannian manifolds},
  author={Nash, John},
  journal={Annals of mathematics},
  volume={63},
  number={1},
  pages={20--63},
  year={1956},
  publisher={JSTOR}
}

@article{whitney1936differentiable,
  title={Differentiable manifolds},
  author={Whitney, Hassler},
  journal={Annals of Mathematics},
  volume={37},
  number={3},
  pages={645--680},
  year={1936},
  publisher={JSTOR}
}

@article{olszewski2025data,
  title={Data vulnerability index for the “crowding problem” in nonlinear dimensionality reduction},
  author={Olszewski, Dominik},
  journal={Neurocomputing},
  pages={130619},
  year={2025},
  publisher={Elsevier}
}

@book{borg05,
	title={Modern multidimensional scaling: Theory and applications},
	author={Borg, Ingwer and Groenen, Patrick JF},
	year={2005},
	publisher={Springer Science \& Business Media}
}

@article{lanczos1950,
  title={An iteration method for the solution of the eigenvalue problem of linear differential and integral operators},
  author={Lanczos, Cornelius},
  year={1950},
  publisher={United States Governm. Press Office Los Angeles, CA},
  journal={N.A.}
}

@article{SMACOF,
  title={Multidimensional scaling using majorization: SMACOF in R},
  author={De Leeuw, Jan and Mair, Patrick},
  year={2011},
  journal={N.A.}
}

@article{shiebler2020,
  title={Functorial manifold learning},
  author={Shiebler, Dan},
  journal={arXiv preprint arXiv:2011.07435},
  year={2020}
}

@inproceedings{shiebler2020clustering,
  title={Functorial clustering via simplicial complexes},
  author={Shiebler, Dan},
  booktitle={TDA $\{$$\backslash$\&$\}$ Beyond},
  year={2020},
  journal={N.A.},
  pages={N.A.}
}

@Article{Orts2019,
author={Orts, F.
and Filatovas, E.
and Ortega, G.
and Kurasova, O.
and Garz{\'o}n, E. M.},
title={Improving the energy efficiency of SMACOF for multidimensional scaling on modern architectures},
journal={The Journal of Supercomputing},
year={2019},
month={Mar},
day={01},
volume={75},
number={3},
pages={1038-1050},
issn={1573-0484},
doi={10.1007/s11227-018-2285-x},
url={https://doi.org/10.1007/s11227-018-2285-x}
}

@misc{McInnes18,
  doi = {10.48550/ARXIV.1802.03426},
  url = {https://arxiv.org/abs/1802.03426},
  author = {McInnes, Leland and Healy, John and Melville, James},
  title = {UMAP: Uniform Manifold Approximation and Projection for Dimension Reduction},
  publisher = {arXiv},
  year = {2018},
  copyright = {arXiv.org perpetual, non-exclusive license},
}

@article{spivak_learners_2022,
	title = {Learners' {Languages}},
	volume = {372},
	issn = {2075-2180},
	url = {http://arxiv.org/abs/2103.01189v2},
	doi = {10.4204/EPTCS.372.2},
	language = {en},
	urldate = {2023-12-01},
	journal = {Electronic Proceedings in Theoretical Computer Science},
	author = {Spivak, David I.},
	month = nov,
	year = {2022},
	pages = {14--28},
}

@book{spivak_poly_book,
	title = {Polynomial functors - A general theory of interaction},
	author = {Spivak, David I and Niu, Nelson},
  	year = {2023},
	publisher = {Topos Institute},
}

@misc{fong_backprop_2019,
	title = {Backprop as {Functor}: {A} compositional perspective on supervised learning},
	shorttitle = {Backprop as {Functor}},
	url = {http://arxiv.org/abs/1711.10455},
	language = {en},
	urldate = {2023-11-19},
	publisher = {arXiv},
	author = {Fong, Brendan and Spivak, David I. and Tuyéras, Rémy},
	month = may,
	year = {2019},
	note = {Number: arXiv:1711.10455
arXiv:1711.10455 [cs, math]},
}

@misc{Spivak09,
  title={METRIC REALIZATION OF FUZZY SIMPLICIAL SETS},
  author={David I. Spivak},
  year={2009},
  url={https://api.semanticscholar.org/CorpusID:35883069},
  note={\href{http://www.dspivak.net/metric\_realization090922.pdf}{http://www.dspivak.net/metric\_realization090922.pdf}}
}

@book{MacLane78,
  title={Categories for the Working Mathematician},
  author={Saunders {Mac Lane}},
  isbn={9780387984032},
  series={Graduate Texts in Mathematics},
  year={1978},
  publisher={Springer}
}

@article{Leinster16,
  doi = {10.48550/ARXIV.1612.09375},
  
  url = {https://arxiv.org/abs/1612.09375},
  
  author = {Leinster, Tom},
  
  title = {Basic Category Theory},
  journal = {arXiv},
  publisher = {arXiv},
  
  year = {2016},
  
  copyright = {Creative Commons Attribution Non Commercial Share Alike 4.0 International}
}

@article{Carlsson09,
author = {Carlsson, Gunnar},
year = {2009},
month = {04},
pages = {255-308},
title = {Topology and Data},
volume = {46},
journal = {Bulletin of the  American Mathematical Society},
doi = {10.1090/S0273-0979-09-01249-X}
}

@book{Gabriel67,
  title={Calculus of fractions and homotopy theory},
  author={Gabriel, P. and Zisman, M.},
  series={Ergebnisse der Mathematik und ihrer Grenzgebiete},
  url={https://books.google.de/books?id=y94gwgEACAAJ},
  year={1967},
  publisher={Springer}
}

@book{Goerss09,
  title={Simplicial Homotopy Theory},
  author={Goerss, P.G. and Jardine, J.F.},
  isbn={9783034601894},
  series={Modern Birkh{\"a}user Classics},
  url={https://books.google.de/books?id=ED1bVh5K-5YC},
  year={2009},
  publisher={Birkh{\"a}user Basel}
}

@misc{nlab:simplicialIdentities,
  author = {{nLab}},
  title = {simplicial identities},
  url = {https://ncatlab.org/nlab/revision/simplicial+identities/14},
  month = jul,
  year = 2023
}

@article{Joharinad_Fahimi_Barth_Keck_Jost_2025, title={IsUMap: Manifold Learning and Data Visualization Leveraging Vietoris-Rips Filtrations}, volume={39}, url={https://ojs.aaai.org/index.php/AAAI/article/view/33946}, DOI={10.1609/aaai.v39i17.33946}, abstractNote={This work introduces IsUMap, a novel manifold learning technique that enhances data representation by integrating aspects of UMAP and Isomap with Vietoris-Rips filtrations and metric realization of one-parameter filtrations of simplicial complexes. Inferring topological information from combinatorial models which have been built according to metric relations (Vietoris-Rips complexes) has proven useful in topological data analysis and general machine learning applications. This encourages the use of such objects for geometric inference.
We extend this research direction by proposing a clear theoretical pipeline that not only provides a comprehensive guide for assigning a (triangulated) metric space to every admissible one- parameter filtration of simplicial complexes but also offers a method for merging these objects. With this, our method presents a systematic and detailed construction of a metric representation for locally distorted metric spaces that captures complex data structures more accurately than the
previous schemes. Our approach addresses limitations in existing methods by accommodating non-uniform data distributions and intricate local geometries. We validate its performance through extensive experiments on examples with known geometries and in applications to data, in particular from computational biology.}, number={17}, journal={Proceedings of the AAAI Conference on Artificial Intelligence}, author={Joharinad, Parvaneh and Fahimi, Hannaneh and Barth, Lukas Silvester and Keck, Janis and Jost, Jürgen}, year={2025}, month={Apr.}, pages={17699-17706} }

@misc{VRpaper,
      title={IsUMap: Manifold Learning and Data Visualization leveraging Vietoris-Rips filtrations}, 
      author={Lukas Silvester Barth and Fatemeh and Fahimi and Parvaneh Joharinad and Jürgen Jost and Janis Keck},
      year={2024},
      eprint={2407.17835},
      archivePrefix={arXiv},
      primaryClass={cs.LG},
      url={https://arxiv.org/abs/2407.17835}, 
}

@book{IsumapBook,
  title = {Data visualization with category theory and geometry},
  author={Barth, Lukas Silvester and Fahimi, Fatemeh (Hannaneh) and Joharinad, Parvaneh and Jost, Jürgen and Keck, Janis},
  publisher={Mathematics of Data Series, Springer},
  year={2025}
}

@article{Tenenbaum00,
  title={A global geometric framework for nonlinear dimensionality reduction},
  author={Tenenbaum, Joshua B and Silva, Vin de and Langford, John C},
  journal={science},
  volume={290},
  number={5500},
  pages={2319--2323},
  year={2000},
  publisher={American Association for the Advancement of Science}
}

@techreport{Bernstein00,
  title={Graph approximations to geodesics on embedded manifolds},
  author={Bernstein, Mira and De Silva, Vin and Langford, John C and Tenenbaum, Joshua B},
  year={2000},
  institution={Citeseer}
}

@inproceedings{Pearson1901,
  title={On lines and planes of closest fit to systems of points in space},
  author={Karl Pearson},
  booktitle={Proceedings of the 17th ACM SIGACT-SIGMOD-SIGART symposium on Principles of database systems (SIGMOD)},
  pages={19},
  year={1901}
}

@article{Damrich2021,
  title={On UMAP's true loss function},
  author={Damrich, Sebastian and Hamprecht, Fred A},
  journal={Advances in Neural Information Processing Systems},
  volume={34},
  pages={5798--5809},
  year={2021}
}

@article{Joharinad19,
author = {P. Joharinad and J. Jost},
title = {{Topology and curvature of metric spaces}},
doi = {10.1016/j.aim.2019.106813},
pages = {106813},
year = {2019},
volume = {356},
issn = {0001-8708},
journal = {Advances in mathematics}
}

@article{douze2024faiss,
      title={The Faiss library},
      author={Matthijs Douze and Alexandr Guzhva and Chengqi Deng and Jeff Johnson and Gergely Szilvasy and Pierre-Emmanuel Mazaré and Maria Lomeli and Lucas Hosseini and Hervé Jégou},
      year={2024},
      eprint={2401.08281},
      archivePrefix={arXiv},
      primaryClass={cs.LG},
      journal={arXiv}
}

@article{moon2017phate,
  title={PHATE: a dimensionality reduction method for visualizing trajectory structures in high-dimensional biological data},
  author={Moon, Kevin R and van Dijk, David and Wang, Zheng and Chen, William and Hirn, Matthew J and Coifman, Ronald R and Ivanova, Natalia B and Wolf, Guy and Krishnaswamy, Smita},
  journal={BioRxiv},
  volume={120378},
  year={2017},
  publisher={Cold Spring Harbor Laboratory}
}

@inproceedings{tang2016largeVis,
  title={Visualizing large-scale and high-dimensional data},
  author={Tang, Jian and Liu, Jingzhou and Zhang, Ming and Mei, Qiaozhu},
  booktitle={Proceedings of the 25th international conference on world wide web},
  pages={287--297},
  year={2016}
}

@article{vandermaaten08a,
  author  = {Laurens van der Maaten and Geoffrey Hinton},
  title   = {Visualizing Data using t-SNE},
  journal = {Journal of Machine Learning Research},
  year    = {2008},
  volume  = {9},
  number  = {86},
  pages   = {2579--2605},
  url     = {http://jmlr.org/papers/v9/vandermaaten08a.html}
}

@inproceedings{luecken2021sandbox,
  title={A sandbox for prediction and integration of DNA, RNA, and proteins in single cells},
  author={Luecken, Malte D and Burkhardt, Daniel Bernard and Cannoodt, Robrecht and Lance, Christopher and Agrawal, Aditi and Aliee, Hananeh and Chen, Ann T and Deconinck, Louise and Detweiler, Angela M and Granados, Alejandro A and others},
  booktitle={35th conference on neural information processing systems (NeurIPS 2021) track on datasets and benchmarks},
  year={2021}
}

@article{Krauss2025,
author = {Krau\ss{}, Dana and Moreno-Viedma, Veronica and Adachi-Fernandez, Emi and de S\'{a} Fernandes, Cristiano and Genger, Jakob-Wendelin and Fari, Ourania and Blauensteiner, Bernadette and Kirchhofer, Dominik and Bradaric, Nikolina and Gushchina, Valeriya and Fotakis, Georgios and Mohr, Thomas and Abramovich, Ifat and Mor, Inbal and Holcmann, Martin and Bergthaler, Andreas and Haschemi, Arvand and Trajanoski, Zlatko and Winkler, Juliane and Gottlieb, Eyal and Sibilia, Maria},
title = {EGFR controls transcriptional and metabolic rewiring in KRAS G12D colorectal cancer},
journal = {EMBO Molecular Medicine},
volume = {17},
number = {6},
pages = {1355-1392},
doi = {https://doi.org/10.1038/s44321-025-00240-4},
url = {https://www.embopress.org/doi/abs/10.1038/s44321-025-00240-4},
eprint = {https://www.embopress.org/doi/pdf/10.1038/s44321-025-00240-4},
year = {2025}
}

@article{Bentley1975,
author = {Bentley, Jon Louis},
title = {Multidimensional binary search trees used for associative searching},
year = {1975},
issue_date = {Sept. 1975},
publisher = {Association for Computing Machinery},
address = {New York, NY, USA},
volume = {18},
number = {9},
issn = {0001-0782},
url = {https://doi.org/10.1145/361002.361007},
doi = {10.1145/361002.361007},
journal = {Commun. ACM},
month = {sep},
pages = {509–517},
numpages = {9}
}

@book{Burago01,
    author = {D. Burago and  Yu. Burago and S. Ivanov},
     title = {A course in metric geometry},
     publisher = {AMS},
     address = {},
      year = {2001},
}

@Article{Torgerson1952,
author={Torgerson, Warren S.},
title={Multidimensional scaling: I. Theory and method},
journal={Psychometrika},
year={1952},
month={Dec},
day={01},
volume={17},
number={4},
pages={401-419},
issn={1860-0980},
doi={10.1007/BF02288916},
url={https://doi.org/10.1007/BF02288916}
}

@article{Belkin03,
     title = {Laplacian Eigenmaps for Dimensionality Reduction and Data Representation},
     author = {Belkin, M. and Niyogi, P.},
     journal = {Neural Computation},
     volume = {15},
     year = {2003},
     pages = {1373--1396},
     }
\bibliographystyle{apalike}

\newpage
\section{Appendix}
\label{sec:appendixx}
\appendix

\section{Preliminaries}
\label{sec:Preliminaries}

\bp{prop:density}{(\textbf{Density theorem})
Every presheaf $X\in\fc{C}$ of any category $\bo{C}$ is the colimit of representables:
\begin{equation}
    \begin{split}
    X\simeq \te{colim}(\bo{y}\circ P_X),
    \end{split}
\end{equation}
which means it is the colimit of the diagram
\bt \bo{El}(X)\ar{r}{P_X}\& \bo{C} \ar{r}{\bo{y}} \& \fc{C}\et, where $\bo{El}(X)$ is the category of elements of $X$.}
\pr{c.f. \cite[p. 155, Theorem 6.2.17]{Leinster16}.}

\bp{prop:RepsPreserveLimits}{(\textbf{Representables preserve limits}) Let $\bo{A}$
be locally small and $D:\bo{I}\to\bo{A}$ a diagram. Then
\begin{equation*}
    \begin{split}
    &\text{Hom}_{\bo{A}}(A,\text{lim}_{\bo{I}} (D))\simeq \text{lim}_{\bo{I}} [\text{Hom}_{\bo{A}}(A,D)]\\
    \text{ and }\quad&
    \text{Hom}_{\bo{A}}(\text{colim}_{\bo{I}} (D),A)\simeq \text{lim}_{\bo{I}^{\text{op}}} [\text{Hom}_{\bo{A}}(D,A)]
    \end{split}
\end{equation*}
}
\pr{See, for instance, \cite[Proposition 6.2.2] {Leinster16}.}

\bp{prop:kanExtension}{
  (\textbf{Kan extensions}) 
  Let $F : \mathbf{A} \to \bo{B}$ be a functor between small categories,
  and let $B\in\bo{B}$. Assume that $(F \Ra B)$ is a comma category and
  $P_B : (F\Ra B) \to \bo{A}$ the canonical projection functor.
  suppose $X : \bo{A} \ra \bo{S}$ is a functor from $\bo{A}$ to a category $\bo{S}$
  with small colimits. Then
  \begin{enumerate}
    \item
      For each $B\in\bo{B}$, if $\te{Lan}_FX(B)$ is the colimit of the diagram
      \begin{equation}
        \begin{split}
          \bt
          (F\Ra B) \ar{r}{P_B} \& \bo{A} \ar{r}{X} \& \bo{S},
          \et
        \end{split}\label{diagram:KanExtension}
      \end{equation}
      then $\te{Lan}_F X : \bo{B}\to\bo{S}$ defines a functor.
    \item
      For functors
      $Y : \bo{B}\ra\bo{S}$, there is a canonical bijection between natural transformations
      $\te{Lan}_FX \to Y$ and natural transformations $X \to Y \circ F$.
    \item
      For any category $\bo{S}$ with small colimits, the functor
      \begin{equation}
        \begin{split}
          -\circ F: [\bo{B},\bo{S}]\to[\bo{A},\bo{S}].
        \end{split}
      \end{equation}
      has a left adjoint, $\te{Lan}_F$, which is called the \textit{left Kan-extension} along $F$.
    \item
      When $\bo{S}$ has both small limits and colimits, the
      functor $- \circ F$ has both left and right adjoints.
      The right adjoint is $\te{Ran}_F$, the \textit{right Kan-extension} along $F$.  $\te{Ran}_FX(B)$ is defined as
      the limit of the diagram 
      \begin{equation}
        \begin{split}
          \bt
          (B \Ra F) \ar{r}{P^B} \& \bo{A} \ar{r}{X} \& \bo{S}
          \et
        \end{split}\label{diagram:RightKanExtension}
      \end{equation}
      where the objects of $(B\Ra F)$ are pairs $(A,h:B\ra F(A))$ and the morphisms
      are maps $f:A\ra A'$ s.t. $h'=F(f)\circ h$. The projection
      is defined by $P^B(A,h):=A$ and $P^B(f):=f$.
  \end{enumerate}
}
\pr{
  See \cite[p. 158, exercise 6.6.25]{Leinster16}. 
}
\bp{prop:adjunctionsPreserveLimits}{
  Let $F\dashv G$ be an adjunction. We also say that $F$ is left
  adjoint to $G$ and $G$ right adjoint to $F$. Then $F$ and $G$ preserve all colimits and all limits respectively.
}
\pr{
  See e.g.~\cite[p. 159, theorem 6.3.1]{Leinster16}.
}

\bp{prop:LeftAdjointFullyFaithfulUnitIso}{
    The right (left) adjoint of an adjunction is fully faithful if and only if the counit (unit) is an isomorphism.
}
\pr{
    Given an adjunction $ F \dashv U $, with $ F : \bo{A} \to \bo{B} $ and counit $ \epsilon : FU \to \te{id} $, this follows directly from the commutativity of the diagram
    \begin{center}
        \begin{tikzcd}[ampersand replacement=\&]
            \& {\te{Hom}(FU(X),Y)} \arrow[d, "\cong"] \\
            {\te{Hom}(X,Y)} \arrow[r, "U"'] \arrow[ru, "\epsilon^*_X"] \& {\te{Hom}(U(X),U(Y))}                 
        \end{tikzcd}
    \end{center}
    for all $ X, Y \in \bo{A} $ where the vertical isomorphism comes from the adjunction.
}

\bp{prop:colimitsAreCoEqualizersOfCoproducts}{
    If a category has all small coproducts and small coequalizers, then it has all small colimits.
    The reason is that, given a small diagram $D:\bo{I}\to \bo{A}$, its colimit can be expressed as the coequalizer \bt \coprod_{I\in\bo{I}} D(I) \ar{r}{e} \&  L \et of the following diagram
    \begin{equation}
    \bt
    {\displaystyle \coprod_{u:J\ra K\text{ in }\bo{I}} D(J)}  \ar[yshift=0.3em]{r}{s}\ar[yshift=-0.3em,']{r}{t} \&
    {\displaystyle \coprod_{I\in\bo{I}} D(I) }
    \et,
    \label{diagram:equalizerProduct}
\end{equation}
where the $u$-th summand of $s$ is the composite $\bt 
D(J) \ar{r}{Du} \& D(K) \ar{r}{i_K} \&
{\coprod_{I\in\bo{I}} D(I)}
\et$ and of $t$ is $\bt
D(J) \ar{r}{i_J} \&
{\prod_{I\in\bo{I}} D(I) }
\et$. ($i_K$ denotes a coproduct injection. The universal property of the coproduct then ensures that defining $s$ and $t$ for each summand defines them completely.)
}
\pr{This is the dual of \cite[Proposition 5.1.26 (a)]{Leinster16}.}

\section{Proofs of section \ref{sec:DiscreteApprox}}
\label{sec:ProofsOfSection:sec:DiscreteApprox}

\subsection{Proof of Proposition \ref{prop:colimitUM}}
\label{proof:prop:colimitUM}

The idea is to use the fact that every colimit can be constructed as the coequalizer of coproducts by Proposition \ref{prop:colimitsAreCoEqualizersOfCoproducts}. We thus use two lemmata about the exact form of the coproducts and coequalizers and then combine those results at the end to obtain an explicit description of the colimit in $\bo{UM}$ and $\bo{EPMet}$.

To begin, recall that (see, for example, \cite[chapter III.3]{MacLane78}), in $\bo{Sets}$, the coproduct $\coprod$ is the disjoint union $\bigcupdot$ and the coequalizer $Y$ of the diagram
\bt
A  \ar[yshift=0.3em]{r}{f}\ar[yshift=-0.3em,']{r}{g} \&
    X 
\et
is 
$Y = X~ / \sim$,
where the equivalence relation $\sim$ is generated by 
\begin{equation}
\begin{split}
x\sim x'\quad\text{ iff }\quad \exists a:x=f(a)\text{ and }x'=g(a).
\end{split}
\label{eq:coeqRelation}
\end{equation}

In \cite{Spivak09}, a brief characterization of the coproduct and coequalizer in $\bo{UM}$ was given in his existence proof of co-completeness of $\bo{UM}$. One can slightly modify the constructions (in particular for the coequalizer) to obtain the following two lemmata.

\lm{prop:coprodUber}{
    Let $A$ be a set and for all $a \in A$, let $(X_a , d_a)$ either be objects of $\bo{UM}$ or of $\bo{EPMet}$. We abuse notation and denote by $F:\bo{UM}\to \bo{Sets}$ and $F:\bo{EPMet}\to \bo{Sets}$ the forgetful functors that map either a uber-metric space or an extended pseudo-metric space $(Y,d)$ to its underlying set $Y$. Let $X:=\coprod_{a\in A}F(X_a,d_a)$
    be the coproduct, which is the disjoint union, in the category of sets. Let $d_{\coprod}$ denote the metric such that for all $x, x' \in X$,
    \begin{equation}
        \begin{split}
            d_{X}(x,x') := \begin{cases}
                d_a (x, x' ),&\text{if }\exists a~:~x,x'\in F(X_a,d_a),\\
                \infty,&\text{else.}
            \end{cases}
        \end{split}
        \label{eq:coprodMetric}
    \end{equation}
    Then the coproduct in $\bo{UM}$ or $\bo{EPMet}$ is given given by $(X,d_X)$, i.e.~$\coprod_{a\in A}(X_a,d_a) \simeq (X,d_X)$.
}
\lm{prop:coequalizerUber}{
    Suppose that \bt
            A  \ar[yshift=0.3em]{r}{f}\ar[yshift=-0.3em,']{r}{g} \&
                X\ar{r}{e} \&  Y 
            \et
    is a coequalizer diagram of sets. Then $Y \simeq X/\sim$ where the equivalence relation $\sim$ is generated by
    \begin{equation}
        \begin{split}
            x \sim x'\quad\text{iff}\quad \exists a \in A \text{ s.t. } x =
            f (a)\te{ and }x' = g(a).
        \end{split}
        \label{eq:coeqRelation2}
    \end{equation}
    If $(X,d_X)$ is an object of $\bo{UM}$, we define a metric $d_{\sim}$ on $Y$ by
    \begin{equation}
        \begin{split}
            d_{\sim} ([x], [x']) = \inf(d_X (p_1 , q_1 ) + \cdots + d_X (p_n , q_n )),
        \end{split}
        \label{eq:simMetric}
    \end{equation}
    where the infimum is taken over all pairs of sequences $(p_1 , \cdots , p_n ),~ (q_1 , \cdots , q_n )$ of
    elements of $X$, such that
    \begin{equation}
        \begin{split}
            p_1 \sim x,\quad q_n \sim x',\quad \text{ and }\quad p_{i+1} \sim q_i ~\te{ for all } ~1 \le i \le n - 1.
        \end{split}
    \end{equation}
    Then $(Y, d_{\sim} )$ is a coequalizer in $\bo{UM}$.\\
    If $(X,d_X)$ is instead an object of $\bo{EPMet}$, then we define a metric $d_{\sim}^{\te{EPMet}}$ on $Y$ by
    \begin{equation}
      \begin{split}
        d_{\sim}^{\te{EPMet}}([x],[x']) := 
        \begin{cases}
          \infty,&\te{if }\inf_{y\in[x],y'\in [x']}=\infty,\\
          d_{\sim}\qquad\te{(as in \eqref{eq:simMetric})},&\te{else.}
        \end{cases}
      \end{split}
      \label{eq:simMetricEPMet}
    \end{equation}
    and $(Y, d_{\sim}^{\te{EPMet}} )$ is a coequalizer in $\bo{EPMet}$. 
}

With those at hand, given a small diagram $D:\bo{I}\to\bo{UM}$, or $D:\bo{I}\to\bo{EPMet}$, one proceeds to construct the objects in diagram \eqref{diagram:equalizerProduct}. Due to Prop. \ref{prop:coequalizerUber}, we know that the coequalizer of this diagram is obtained as the coequalizer of the corresponding diagram in the category of sets,
\begin{equation}
    \bt
    {\displaystyle \coprod_{u:J\ra K\text{ in }\bo{I}} FD(J)}  \ar[yshift=0.3em]{r}{s}\ar[yshift=-0.3em,']{r}{t} \&
    {\displaystyle \coprod_{I\in\bo{I}} FD(I) }
    \et
    \label{eq:equalizerSetsColim}
\end{equation}
which is defined in eq.~\eqref{eq:coeqRelation}, 
equipped with the metric $d_{\sim}$ defined in eq.~\eqref{eq:simMetric}, or $d_{\sim}^{\te{EPMet}}$ defined in eq.~\eqref{eq:simMetricEPMet}, where $X$ is now given by $X:=\coprod_{I\in\bo{I}} FD(I)$. By Proposition \ref{prop:coprodUber}, we know that $d_X$ must hence be defined as in equation \eqref{eq:coprodMetric}. 

It only remains to simplify the definition of the equivalence relation. Since we know that the coequalizer of the diagram \eqref{eq:equalizerSetsColim} is corresponding to the colimit of $FD$ in $\bo{Sets}$ and we know the colimit to be isomorphic to $X/\sim$ with $\sim$ defined in \eqref{eq:colimSetsFD}, we can check that this simple relation is consistent with the more convoluted relation obtained by directly applying eq.~\eqref{eq:coeqRelation} to the diagram \eqref{eq:equalizerSetsColim}. Indeed, we obtain for all $x,x'\in X$, 
\begin{equation}
    \begin{split}
        x\sim x'
        \quad\text{iff}\quad \exists y\in \coprod_{u:J\ra K\text{ in }\bo{I}} FD(J)~:~x=s(y),~x'=t(y),
    \end{split}
\end{equation}
and one can rewrite the condition as follows:
\begin{equation}
    \begin{split}
        &\exists y\in \coprod_{u:J\ra K\text{ in }\bo{I}} FD(J)~\te{ s.t. }~x=t(y),~x'=s(y)\\
        \Lr \quad& \exists u:J\to K~\te{ s.t. }~\exists z \in FD(J)~\te{ s.t. }~ ~x=i_J(z),~x'=i_K(FDu(z)) \\
        \Lr \quad& \exists u:J\to K~\te{ s.t. }~(u,FD(u))(x=(J,z))=(K,FDu(z))=x'.
    \end{split}
\end{equation}
Since $i_J$ and $i_K$ are injective, it is common to identify $x$ with $z$ and $x'$ with $FDu(z)$ and hence just write, with a slight abuse of notation, $FDu(x)=x'$.
This completes the proof.
\hfill$\blacksquare$

\section{Proofs of section \ref{sec:FuzzySimplicialSets}}

\subsection{Proof of Proposition \ref{prop:isoClassicalAndSimplicialFuzzySets}}
\label{proof:prop:isoClassicalAndSimplicialFuzzySets}

We define the functor $C : \bo{Fuz} \to \bo{cFuz}$ as follows: to a fuzzy set $S \in \bo{Fuz}$ we assign the classical fuzzy set $(S(0), \xi_S)$, where $\xi(s) := \sup\left\{ a \in [0,1] \, | \, s \in \mathrm{im}(S(0\to a)) \right\}$. Since $S(0\to a)$ is injective, $s \in \mathrm{im}(S(0\to a))$ corresponds to a unique element in $S(a)$ and we therefore permit ourselves henceforth to slightly abuse and simplify notation by writing $s\in S(a)$ instead of $s\in \mathrm{im}(S(0\to a))$. To a morphism $f : S \to T$ of fuzzy sets we assign the map $(S(0), \xi_S) \to (T(0), \xi_T) $, given by $f(0)$, or $f$ for short. This is indeed a map of classical fuzzy sets: for all $s \in S(0)$ we have
\begin{equation}
     \begin{split}
    \xi_T(f(s)) = \sup\left\{ a \in [0,1] \, | \, f(s) \in T(a) \right\} \ge \sup\left\{ a \in [0,1] \, | \, s \in S(a) \right\} = \xi_S(s),
    \end{split}
\end{equation}
since by definition $s \in S(a)$ implies that $f(s) \in T(a)$.

The inverse functor $C^{-1} : \bo{cFuz} \to \bo{Fuz}$ assigns to each classical fuzzy set $(X,\eta)$ the fuzzy set $\eta^{-1}[-,1]$, which evaluated at $a \in [0,1]$ is the set $ \eta^{-1}[a,1]$. To a map of classical fuzzy sets $g : (X,\eta) \to (Y,\zeta)$ we assign the map of fuzzy sets, whose $a$-th component is $ g : \eta^{-1}[a,1] \to \zeta^{-1}[a,1]$.

We now show that $C$ and $C^{-1}$ are indeed inverse to each other. Given a fuzzy set $S$ we have
\begin{equation}
     \begin{split}
    (C^{-1}C(S))(a) = C^{-1}(S(0), \xi_S) = \xi_S^{-1}[a,1] = \Big\{x \in S(0) \, \Big| \, \sup\{b \in [0,1] \, | \, x \in S(b) \} \in [a,1] \Big\}.
    \end{split}
\end{equation}
Now, if $ x \in S(a) $, then $ \sup\{b \in [0,1] \, | \, x \in S(b)\} \ge a $ so that $x \in  (C^{-1}C(S))(a) $. Conversely, if $x \in (C^{-1}C(S))(a) $, then by definition $ \sup\{b \in [0,1] \, | \, x \in S(b) \} \ge a $ and
\begin{equation}
     \begin{split}
    x \in \bigcap_{\substack{b \in [0,1] \\ \text{s.t. } x \in S(b)}} S(b) \simeq S\Big(\sup\{b \in [0,1] \, | \, x \in S(b) \}\Big) \subset S(a).
    \end{split}
\end{equation}
where the isomorphism is given by the gluing condition (note that the set $\{b \in [0,1] \, | \, x \in S(b)\}$ is non-empty, since by definition $x \in S(0)$). Thus, $ (C^{-1}C(S))(a) \simeq S(a) $. Naturality in $a$, as in the diagram
\begin{equation}
\begin{tikzcd}
	{S(b)} & {(C^{-1}C(S))(b)} \\
	{S(a)} & {(C^{-1}C(S))(a)}
	\arrow[from=1-1, to=1-2]
	\arrow[from=1-1, to=2-1]
	\arrow[from=1-2, to=2-2]
	\arrow[from=2-1, to=2-2]
\end{tikzcd}
\end{equation}
essentially boils down to the functoriality of $S$, namely that the composite restriction $S(b) \hookrightarrow S(a) \hookrightarrow S(0)$ equals the restriction $S(b) \hookrightarrow S(0)$. The naturality in $S$, as in the diagram
\begin{equation}
\begin{tikzcd}
	{S(a)} & {(C^{-1}C(S))(a)} \\
	{T(a)} & {(C^{-1}C(T))(a)}
	\arrow[from=1-1, to=2-1]
	\arrow[from=1-1, to=1-2]
	\arrow[from=1-2, to=2-2]
	\arrow[from=2-1, to=2-2]
\end{tikzcd}
\end{equation}
essentially boils down to the naturality of $S \to T$, namely that $ S(a) \hookrightarrow S(0) \to T(0)$ equals $ S(a) \to T(a) \hookrightarrow T(0)$.

On the other hand, given a classical fuzzy set $(X,\eta)$ we have
\begin{equation}
    \begin{split}
    CC^{-1}(X,\eta) &= C(\eta^{-1}[-,1]) \\
    &= (\eta^{-1}[0,1], \xi), \quad \text{where } \xi(x) = \sup\{ a \in [0,1] \, | \, x \in \eta^{-1}[a,1] \} \\
    &= (X, \xi), \quad \text{where } \xi(x) = \sup\{ a \in [0,1] \, | \, \eta(x) \in [a,1] \} \\
    &= (X,\eta).
\end{split}
\end{equation}
This is automatically natural in $(X,\eta)$, as we have an equality on the nose.

\hfill$\blacksquare$

\subsection{Proof of Proposition \ref{prop:simplexFaceStrength}}
\label{proof:prop:simplexFaceStrength}

Let $ S : \Delta^{\te{op}} \to \bo{Fuz} $ be a fuzzy simplicial set. Given any morphism $ f : n \to m $ in $ \Delta $ we get by functoriality a natural transformation $ S(f) : S(n) \to S(m) $. This means, for all $ a \in [0,1] $ there exists a morphisms $ S(f)_a : S(n,a) \to S(m,a) $. Now, recall that an element $ x \in S(n,a) $ is an $n$-simplex of strength at least $ a $. Here the strength of a simplex is defined by passing over to classical fuzzy sets as described in Proposition \eqref{prop:isoClassicalAndSimplicialFuzzySets}. It follows that the image $ S(f)_a(x) \in S(m,a) $ is again of strength at least $ a $ and as a result, by taking suprema, we get $ \te{str}(x) \le \te{str}(S(f)_a(x)) $.

The first part of the proposition follows by choosing $ f $ to be a face map. If, on the other hand $ f $ is a degeneracy $ \sigma : n \to n-1 $, then there exists a face map $ \delta : n-1 \to n $ such that $ \sigma \circ \delta = \te{id}_{n-1} $. The above argument then implies
\begin{align}
    \te{str}(x) \le \te{str}(S(\sigma)_a(x)) \le \te{str}(S(\delta)_a(S(\sigma)_a(x))) = \te{str}(x).
\end{align}
from which we conclude $ \te{str}(x) = \te{str}(S(\sigma)_a(x)) $. \hfill $\blacksquare$

\subsection{Proof of Proposition \ref{prop:Fsiso}}
\label{proof:prop:Fsiso}

Since $\bo{sFuz}=[\Delta^{\te{op}},\bo{Fuz}]$ is a functor category, we can invoke the equivalence $C$ on the codomain catagory to obtain $[\Delta^{\te{op}},\bo{Fuz}]\simeq [\Delta^{\te{op}},\bo{cFuz}]$.
\hfill$\blacksquare$

\subsection{Proof of Proposition \ref{prop:cFuzCompat}}
\label{proof:prop:cFuzCompat}

This follows directly from Proposition \ref{prop:simplexFaceStrength} and the definition of $\xi$ in \eqref{eq:maxStrength}.
\hfill$\blacksquare$

\section{Proofs of section \ref{sec:nFuz}}

\subsection{Proof of Proposition \ref{prop:leftRightAdjointOfTruncation}}
\label{proof:prop:leftRightAdjointOfTruncation}

As proven in Proposition \ref{prop:kanExtension}, items 3 and 4, about Kan-extensions, any precomposition functor, i.e.~any functor of the form $-\circ F$, has both a left and right adjoint if $\bo{S}$ has small limits and colimits. Below we use this theorem to construct two adjunctions between the presheaf categories $[(\bo{Delta}\times \bo{I})^{\te{op}},\bo{Set}]$ and $[(\bo{n}\times \bo{I})^{\te{op}},\bo{Set}]$ and thereafter proceed to derive natural isomorphisms to more explicit descriptions and to show that these adjunctions additionally restrict appropriately to adjunctions between $\bo{sFuz}$ and $\bo{nFuz}$ (meaning that if one starts at either side of the adjunctions with (truncated) fuzzy simplicial sets, then the application of the adjunction again yields fuzzy simplicial sets). 
For the presheaf categories, we can choose
\begin{equation}
    \begin{split}
        \bo{A}=(\bo{n}\times\bo{I})^{\te{op}},\quad \bo{B}=(\Delta\times\bo{I})^{\te{op}},\quad\bo{S}=\bo{Sets},\quad F=I_n,\quad X\in\bo{nFuz}
    \end{split}
\end{equation}
and directly apply Proposition \ref{prop:kanExtension} because $\bo{Sets}$ does have all small limits and colimits. As a consequence, we immediately obtain the result that $-\circ I_n=\te{tr}_n$ has both left and right adjoints, given by the left and right Kan-extensions respectively. By Proposition \ref{prop:adjunctionsPreserveLimits}, this implies that $\te{tr}_n$ preserves all limits and colimits.

\textbf{\large{Skeleton functor}}

Let us now take a closer look at the skeleton functor. Recall first that the comma category $(I_n\Ra (m,a))$ has objects which are pairs $(f:I_n(k,b)\to (m,a),x\in N(k,b))$, where $(k,b)\in(\bo{n}\times\bo{I})$ and $f:I_n(k,b)\to (m,a)$ is a morphism in $(\Delta\times \bo{I})^{\te{op}}$, while morphisms between $(f,x)$ and $(f',x')$ are maps $\gamma:(k',b')\to (k,b)$ in $\bo{n}\times \bo{I}$ such that $f'\circ I_n(\gamma)=f$. To simplify notation, we use that, for all $k\le n$, $(k,b)=I_n(k,b)$ and work in $\Delta\times \bo{I}$ instead of $(\Delta\times \bo{I})^{\te{op}}$.
We then have
\begin{equation} \label{eq:LanSkeletonFormula}
    \begin{split}
        \te{Lan}_{I_n}N(m,a) &=\te{colim}(\bt (I_n\Ra (m,a)) \ar{r}{\te{proj}} \& (\bo{n}\times \bo{I})^{\te{op}} \ar{r}{N} \& \bo{Sets} \et) \\
        &\simeq \smashoperator[r]{\coprod_{f:(m,a) \to (k\le n,b)}} \ N(k,b) \bigg\slash \sim
    \end{split}
\end{equation}
where the equivalence relation is generated by
\begin{equation}
    \label{eq:equivRelationSk}
  (f:(m,a)\to(k,b), x) \sim (f':(m,a)\to(k',b'), x') \quad \iff \quad 
  \begin{matrix}
    \exists \, \gamma: (k',b') \to (k,b)\text{  such that  }  \\[2mm]
    \gamma \circ f'=f\te{ and } N(\gamma)(x)=x'
  \end{matrix}
\end{equation}
and where we used the explicit form of the colimit in the category of sets in the last step.

We can further simplify this expression by exploiting the following identifications:
\begin{enumerate}
    \item Eq.~\eqref{eq:equivRelationSk} allows us to identify any pair $(f=(\sigma,i_{ab}):(m,a)\to(k,b),~x\in N(k,b))$ via the morphism $\gamma=(\te{id}_k,i_{ab})$ with the pair $((\sigma,i_{aa}),~N(\te{id}_k,i_{ab})(x))$. 
    Consider the induced map 
    \begin{align} \label{eq:Restrict_to_a}
        \begin{split}
            B(N,m,a):\smashoperator[r]{\coprod_{f:(m,a) \to (k,b)}} \ N(k,b) \bigg\slash \sim \quad &\longrightarrow \smashoperator[r]{\coprod_{\substack{f:m \to k \\ k \le n}}} \ N(k,a) \bigg\slash \sim \\
            \Big[\Big((\sigma,i_{ab}),x\Big)\Big] &\longmapsto \Big[\Big(\sigma,~N(\te{id}_k,i_{ab})(x)\Big)\Big]
        \end{split}
    \end{align}
    where in the codomain $\sim$ is generated by $(\sigma: m \to k,x)\sim (\sigma':m \to k',x')$ iff $\exists \gamma:k'\to k$ such that $\gamma\circ \sigma'=\sigma$ and $N(\gamma,i_{aa})(x)=x'$. We show that this map
    is a well-defined bijection: 
    \begin{enumerate}
        \item It is well-defined because if $((\sigma,i_{ab}),x)\sim ((\sigma',i_{ab'}),x')$ is a generator of the equivalence relation, i.e.~it is witnessed by $\gamma=(\tau,i_{b'b})$ (where both $b\le a$ and $b'\le a$), then 
        \begin{equation}
            \begin{split}
                (\sigma,~N(\te{id}_k,i_{ab})(x))\sim (\sigma,~N(\te{id}_k,i_{ab'})(x'))
            \end{split}
        \end{equation}
        is simply witnessed by $\tau$.
        \item It is surjective because we can just choose $b=a$ (and then $x=x'$).
        \item It is injective because if $[(\sigma,N(\te{id}_k,i_{ab})(x))]=[(\sigma',N(\te{id}_k,i_{ab'})(y))]$, then by construction we also have
        \begin{equation}
            \begin{split}
                [((\sigma,i_{ab}),x)]&=[(\sigma,~N(\te{id}_k,i_{ab})(x))]=[(\sigma',~N(\te{id}_k,i_{ab'})(y))]
                =[((\sigma',i_{ab'}),y)].
            \end{split}
        \end{equation}
    \end{enumerate}
    The codomain in Equation~\eqref{eq:Restrict_to_a} can be extended to a presheaf on $\Delta\times I$ in the obvious way (defining the action on morphisms $(g,i_{a'a}):(m',a')\to (m,a)$ by $(\sigma:m\to k,~x\in N(k,a))\mapsto (\sigma\circ g,~ N(\te{id}_k,i_{a'a})(x))$).
    One checks that this makes the isomorphism in Equation~\eqref{eq:Restrict_to_a} natural in $(m,a)$. For the naturality in $N$, we consider, for any morphism $N\to N'$, the Diagram
    \begin{equation}
        \begin{split}
            \bt
                \smashoperator[r]{\coprod_{f:(m,a) \to (k,b)}} \ N(k,b) \bigg\slash \sim \ar{rr}{B(N,m,a)} \ar{d}{\te{sk}_n(N\to N')(m,a)}\&\& \smashoperator[r]{\coprod_{\substack{f:m \to k \\ k \le n}}} \ N(k,a) \bigg\slash \sim \ar[dashed]{d}\\
                \smashoperator[r]{\coprod_{f:(m,a) \to (k,b)}} \ N'(k,b) \bigg\slash \sim \ar{rr}{B(N',m,a)}\&\& \smashoperator[r]{\coprod_{\substack{f:m \to k \\ k \le n}}} \ N'(k,a) \bigg\slash \sim
            \et
        \end{split}
        \label{eq:naturalityOfSkInN}
    \end{equation}
    Since we proved that $B(N,m,a)$ (given by \eqref{eq:Restrict_to_a}) is a bijection for all $N,m,a$, it is invertible. We thus simply define the dashed arrow in the above Diagram as the composition $B(N',m,a)\circ \te{sk}_n(N\to N')(m,a) \circ B(N,m,a)$. By construction, this is makes the diagram commute.
    As a result,
    \begin{equation}
        \begin{split}
            \te{sk}_n(N)(m,a) &\simeq \smashoperator[r]{\coprod_{\substack{f: m \to k \\ k \le n}}} \ N(k,a) \bigg\slash \sim~,
        \end{split}
        \label{eq:skOnlyDependsOnA}
    \end{equation}

    \item We can make use of epi-mono factorizations to simplify \eqref{eq:skOnlyDependsOnA} even further. Recall that all morphisms in $\Delta$ factor uniquely as a surjection followed by an injection. In other words $\Delta$ admits an epi-mono factorization.
    
    For any $f:m \to k\le n$ in $\Delta$, we write the epi-mono factorization as \bt m \ar{r}{\te{epi}(f)} \& \te{im}(f) \ar{r}{\te{mon}(f)} \& k \et. We can use this to identify any pair $(f:m\to k\le n,~x\in N(k,a))$ via $\gamma=\te{mon}(f)$ with the pair $(\te{epi}(f),~N(\te{mon}(f),i_{aa})(x))$. This sets up a map
    \begin{align}
        \label{eq:epiMonoEquivMap}
        \begin{split}
            {\coprod_{\substack{f: m \to k \\ k \le n}}} \ N(k,a) \bigg\slash \sim~ &\longrightarrow {\coprod_{\substack{f:m \twoheadrightarrow k \\ k\le n}}} \ N(k,a) \bigg\slash \sim~, \\
            [(f,x)] &\longmapsto \Big[\Big(\te{epi}(f),N\big(\te{mon}(f),i_{aa}\big)(x)\Big)\Big]
        \end{split}
    \end{align}
    where the equivalence relation in the codomain $\sim$ is generated by: $(f:m\to k,x\in N(k,a))\sim (f':m\to k',x'\in N(k',a))$ iff $\exists \gamma:k'\to k$ such that $\gamma\circ f'=f$ in $\Delta$ and $N(\gamma,i_{aa})(x)=x'$. We show that this map
    is well-defined and a bijection.
    \begin{enumerate}
        \item To see that it is well-defined, first note that in $\Delta$, if there is a commutative diagram square
        \begin{equation}
                \begin{tikzcd}
                    A & X \\
                    B & Y
                    \arrow["g", from=1-1, to=1-2]
                    \arrow["e"', two heads, from=1-1, to=2-1]
                    \arrow["m", tail, from=1-2, to=2-2]
                    \arrow["d", dashed, from=2-1, to=1-2]
                    \arrow["h", from=2-1, to=2-2]
                \end{tikzcd}
            \label{eq:squareDiagonalMap}
        \end{equation}
        then there exists a unique lift (the map $d$ stylized with a dashed body), making the diagram commute.\footnote{
            We briefly explain why $\Delta$ has this property: 
            Note that epis in $\Delta$ always admit at least one section (because any section of a monotone surjection, which exists set-theoretically, is automatically monotone).
            Pick any section $s:B\to A$ to define $d=gs$.
            Then the lower triangle commutes because $h=hes=mgs=md$.
            We can use this to show that the upper triangle also commutes: We have $mde=he=mg$ and because $m$ is mono this implies $de=g$.
            To show uniqueness, we assume there was another such lift $d'$ that makes everything commute. Then, by commutativity of the lower square, $md=h=md'$ and since $m$ is mono, $d=d'$.  
        }
        Now if $(f,x)\sim (f',x')$ is a generator of the equivalence relation, i.e.~it is witnessed by $\gamma:k'\to k$ in $\bo{n}$, the following diagram commutes in $\Delta$:
        \begin{equation}
            \begin{split}
                \bt
                k 
                \&\&\&\& k' \ar[']{llll}{\gamma} \\
                \& \te{im}(f) \ar[']{ul}{\te{mon}(f)} 
                \&\& \te{im}(f') \ar{ur}{\te{mon}(f')} \ar[']{ll}{\te{im}(\gamma)}  \& \\
                \&\&m \ar["\te{epi}(f)"{description}]{ul}\ar[',"\te{epi}(f')"{description}]{ur} \ar[',bend right=40]{uurr}{f'} \ar[bend left=40]{uull}{f} \&\&
                \et
            \end{split}
            \label{eq:epiMonoFactInDelta}
        \end{equation}
        where $\te{im}(\gamma)$ is defined as the unique lift in diagram~\eqref{eq:squareDiagonalMap} upon setting $e = \te{epi}(f')$ and $m = \te{mon}(f)$.
        As a consequence, whenever $(f,x)\sim (f',x'=N(\gamma,i_{aa})(x))$ is witnessed by $\gamma$, we also obtain that
        \begin{equation}
            \begin{split}
                (\te{epi}(f),N(\te{mon}(f),i_{aa})(x)) \sim (\te{epi}(f'),N(\te{mon}(f'),i_{aa})(x'))
            \end{split}
        \end{equation}
        is witnessed by $\te{im}(\gamma)$ because the small lower triangle in \eqref{eq:epiMonoFactInDelta} commutes and
        \begin{equation}
            \begin{split}
                N(\te{im}(\gamma),i_{aa})(N(\te{mon}(f),i_{aa})(x)) &= N(\te{mon}(f'),i_{aa})(N(\gamma,i_{aa})(x))\\
                &= N(\te{mon}(f'),i_{aa})(x').
            \end{split}
        \end{equation}
        As a result, \eqref{eq:epiMonoEquivMap} is well-defined.
        \item \eqref{eq:epiMonoEquivMap} is surjective because we can just choose $f$ to be equal to $\te{epi}(f)$. 
        \item To prove injectivity, we must show that if 
        \begin{equation}
            \begin{split}
                [(\te{epi}(f),N(\te{mon}(f),i_{aa})(x))]=[(\te{epi}(f'),N(\te{mon}(f'),i_{aa})(x'))],
            \end{split}
            \label{eq:injConditionEpiMonFact}
        \end{equation}
        then also $[(f,x)]=[(f',x')]$. But by construction the left side of \eqref{eq:injConditionEpiMonFact} equals $[(f,x)]$ and the right side $[(f',x')]$, implying their equality.
    \end{enumerate}

        We extend Equation~\eqref{eq:epiMonoEquivMap} to an isomorphism natural in $(m,a)$. This requires extending the right-hand side to a presheaf on $\Delta \times I$. Note that the functorial action on a morphism in $\Delta \times I$ is forced by naturality. To obtain formulas for explicit computations, we briefly spell out the consequences. On a morphism $i_{ab}$ the action is the obvious one.
        For the simplicial component, recall that
        the action of the functor in the domain on a map $h : m' \to m$ is given by precomposition: $(f,x) \mapsto (fh,x)$. The functorial action of the codomain on $h$ is then given as follows:
        start by taking the precomposite and then apply epi-mono factorization to it. In more detail, if $(f,x)$ is a pair with $f : m \twoheadrightarrow k $ and $x \in N(k,a)$, then the image of $(f,x)$ is $(m' \twoheadrightarrow k', x')$ as in the diagram
    \begin{equation}
        \begin{tikzcd}[row sep=small, column sep=small]
            {x\in N(k,a)} && {x' \in N(k',a)} \\
            k && {k'} \\
            & m \\
            && {m'}
            \arrow[maps to, from=1-1, to=1-3]
            \arrow[tail, from=2-3, to=2-1]
            \arrow["f", two heads, from=3-2, to=2-1]
            \arrow[two heads, from=4-3, to=2-3]
            \arrow["h", from=4-3, to=3-2]
        \end{tikzcd}
    \end{equation}
    By construction this ensures the naturality of Equation~\ref{eq:epiMonoEquivMap} in $(m,a)$. For the naturality in $N$, we proceed in analogy to what we did in Diagram \eqref{eq:naturalityOfSkInN}: We simply define the action on a morphism $N\to N'$ by what is forced upon us by naturality and the fact that \eqref{eq:epiMonoEquivMap} is a bijection for all $N,m,a$.

    As a conclusion, we obtain
    \begin{equation}
        \begin{split}
            \te{sk}_n(N)(m,a) &\simeq \smashoperator[r]{\coprod_{\substack{f:m \twoheadrightarrow k \\ k\le n}}} \ N(k,a) \bigg\slash \sim~,
        \end{split}
        \label{eq:skEpiFormulation}
    \end{equation}

    \item We can remove all remaining gluing conditions 
    in \eqref{eq:skEpiFormulation} 
    by restricting to the class of non-degenerate simplices.
    
    Recall that a simplex $x'\in N(k',a)$ is non-degenerate iff there is no triple $(k',~\gamma:k'\to k,~x\in N(k,a))$, such that $\gamma$ is a surjective non-identity morphism and $N(f,i_{aa})(x)=x'$. 
    Consider a pair $(f, x)$ where $f : m \twoheadrightarrow k $ and $x \in N(k,a)$. There exists a unique non-degenerate simplex $x_* \in N(k_*, a)$ and a unique surjection $s : k \to k_*$\footnote{We explain why the surjection is unique. Suppose we are given a nondegenerate simplex $x$ in some simplicial set $X$ and two surjections $s,s'$ in $\Delta$ such that $X(s)(x) = X(s')(x)$. If $s \not = s'$, then $s$ has a section $d$, which is not a section of $s'$. It follows that $x = X(sd)(x) = X(d)X(s)(x) = X(d)X(s')(x) = X(s'd)(x)$. But, $s'd$ is by assumption not the identity, and thus admits a nontrivial factorization as a surjection $s''$ followed by an injection $d'$. Continuing the calculation we have $x = X(d's'')(x) = X(s'')X(d')(x)$, by which we have written $x$ as a degeneracy, contradicting our assumption. Thus, we must have had $s = s'$.} such that $x = N(s,x_*)$, as in the diagram on the left below.
    \begin{equation}
        \label{eq:nonDegenerateEquivClassesDiagram}
    \begin{tikzcd}
        {x\in N(k,a)} & {x_* \in N(k_*,a)} \\
        k & {k_*} \\
        & m
        \arrow[maps to, from=1-2, to=1-1]
        \arrow["s", two heads, from=2-1, to=2-2]
        \arrow["f", two heads, from=3-2, to=2-1]
        \arrow["sf"', two heads, from=3-2, to=2-2]
    \end{tikzcd}
    \qquad\quad
    \begin{tikzcd}
        k &&&& {k'} \\
        & {k_*} && {k'_*} \\
        && m
        \arrow["\gamma", from=1-1, to=1-5]
        \arrow["s", two heads, from=1-1, to=2-2]
        \arrow["t"', two heads, from=1-5, to=2-4]
        \arrow["f", bend left=25, two heads, from=3-3, to=1-1]
        \arrow["g"', two heads, bend right=25, from=3-3, to=1-5]
        \arrow["sf"', two heads, from=3-3, to=2-2]
        \arrow["tg", two heads, from=3-3, to=2-4]
    \end{tikzcd}
    \end{equation}
    This sets up a map
    \begin{align}
        \begin{split}
            \smashoperator[r]{\coprod_{\substack{f:m \twoheadrightarrow k\\ k\le n}}} \ N(k,a) \bigg\slash \sim~\quad \longrightarrow \smashoperator[r]{\coprod_{\substack{f:m \twoheadrightarrow k\\ k\le n}}} \ N^{\te{nd}}(k,a) \qquad\quad
            (f,x) \longmapsto (sf, x_*).
        \end{split}
        \label{eq:sk_nd_last_step}
    \end{align}
    where $N^{\te{nd}}(k,a)$ is the subset of elements of $N(k,a)$ that are non-degenerate. We now show that this map is well-defined and bijective.

    \begin{enumerate}
        \item Suppose we have two pairs $(f,x)$ and $(g,y)$ with $f$ and $g$ surjective and related by a morphism $\gamma : k \to k'$, as in Diagram \ref{eq:nonDegenerateEquivClassesDiagram} above on the right. The diagram also features the surjections $s$ and $t$ which restrict $x$ and $y$ onto nondegenerates $x_*$ and $y_*$ respectively (as explained in the previous paragraph). By the commutativity of the outer diagram, $\gamma$ is surjective. Then $ t\gamma $ is another surjection which restricts $x$ to a nondegenerate, namely $y_*$. By uniqueness $x_* = y_*$ and $ s = t \gamma$, which implies $ sf = tg $. More generally, this argument shows that if any two pairs $(f,x)$ and $(g,y)$, with $f,g$ surjections and $x$, $y$ nondegenerate, represent the same class in the domain of Equation~\eqref{eq:sk_nd_last_step}, then $(sf,x_*) = (tg,y_*)$. This proves the well-definedess.
        
        \item Surjectivity follows from the fact that each element $(f_*,x_*)$ in the codomain has $[(f_*,x_*)]$ in its preimage (this follows from the uniqueness in the presentation of any simplex as the restriction of a nondegenerate simplex along a surjective simplicial operator).
        
        \item Injectivity can be read off Diagram \ref{eq:nonDegenerateEquivClassesDiagram} on the right (with $\gamma$ removed): given two $(f,x)$ and $(g,y)$, if $(sf,x_*)=(tg,y_*)$, then $(f,x)$ and $(g,y)$ are related by a zigzag and thus belong to the same equivalence class.
    \end{enumerate}    

    We extend Equation~\eqref{eq:skNonDegFormulation} to a natural isomorphism, by extending the codomain to a presheaf on $\Delta\times I$. Again the definition is forced upon us by naturality and we briefly explain the consequences of fixing the action on morphisms in this way to obtain a computationally more explicit description. 
    We begin with a lemma.

    \lm{lem:DegeneracyIsPreservedAcrossDifferentStrengths}{
        For all $a \le b$ and all $k\in \Delta$, given an $x \in N(k,b)$ we have that, $x$ is degenerate if and only if $N(\text{id}, i_{ab})(x) \in N(k,a)$ is degenerate.
    }
    \textbf{Proof of Lemma \ref{lem:DegeneracyIsPreservedAcrossDifferentStrengths}:}
    Assume first that $y:= N(\te{id},i_{ab})(x)$ is degenerate, i.e.~there exists a tuple $k_* \in \bo{n}$, $ y_* \in N(k_*,a)$ and a surjective map $s : k \to k_*$, such that $ y = N(s,a)(y_\ast) $, as in the bottom left of the diagram
    \[\begin{tikzcd}[column sep=6em, row sep=3em]
        {N(k_*,b)} & {N(k,b)} & {N(k_*,b)} \\
        {N(k_*,a)} & {N(k,a)} & {N(k_*,a)}
        \arrow[from=1-1, to=1-2]
        \arrow[bend left=20, equals, from=1-1, to=1-3]
        \arrow[from=1-1, to=2-1]
        \arrow["{x \longmapsto x_*}"', from=1-2, to=1-3]
        \arrow["{\rotatebox{90}{$x$} \, \longmapsto \rotatebox{90}{$y$}}", sloped, from=1-2, to=2-2]
        \arrow[from=1-3, to=2-3]
        \arrow["{y_* \longmapsto y}", from=2-1, to=2-2]
        \arrow[bend right=20, equals, from=2-1, to=2-3]
        \arrow[from=2-2, to=2-3]
    \end{tikzcd}\]
    The surjection $k \to k_*$ has a section $k_* \to k$ in $\bo{n}$, so that $s$ has a retraction, as in the above diagram. A diagram chase shows that $x$ is a degenerate simplex of $x_*$. In particular $x$ is degenerate.

    Assume now that $x$ is degenerate, meaning there exists some $k_* \in \bo{n}$, an $x_* \in N(k_*,b)$ and some surjection $s : k \to k_*$ such that $x = N(s,b)(x_*)$. Then the commutativity of the left square in the above diagram shows that $y$ is degenerate. \hfill$\blacksquare$

    As a result of Lemma \ref{lem:DegeneracyIsPreservedAcrossDifferentStrengths}, the map $N(k,b) \to N(k,a)$ restricts to nondegenerate simplices. In particular, there is a commutative diagram
    \begin{equation} \label{diag:nondegenerate_naturality}
    \begin{tikzcd}
        {N^\te{nd}(k,b)} & {N^\te{nd}(k,a)} \\
        {N(k,b)} & {N(k,a)}
        \arrow[from=1-1, to=1-2]
        \arrow[from=1-1, to=2-1, tail]
        \arrow[from=1-2, to=2-2, tail]
        \arrow[from=2-1, to=2-2]
    \end{tikzcd}
    \end{equation}

    This explains the functorial action of the codomain of Equation~\eqref{eq:sk_nd_last_step} in $a$. As for the functorial action in $m$, consider a map $h : m' \to m$ in $\Delta$. The action of the domain of Equation~\eqref{eq:sk_nd_last_step} on $h$ has already been explained: it is given by the precomposition with $h$ and then applying epi-mono factorization to extract a surjection. The action on the codomain on $h$ is defined as follows: for a pair $(f,x)$ with $f : m \to k$ epi and $x$ nondegenerate we again precompose with $h$ and then apply epi-mono factorization, as in the left triangle in the diagram
    \begin{equation}
        \begin{tikzcd}[row sep=small, column sep=small]
            {x\in N(k,a)} && {x' \in N(k',a)} & {x'_* \in N(k'_*,a)} \\
            k && {k'} & {k'_*} \\
            & m \\
            && {m'}
            \arrow[maps to, from=1-1, to=1-3]
            \arrow[maps to, from=1-4, to=1-3]
            \arrow[tail, from=2-3, to=2-1]
            \arrow[two heads, from=2-3, to=2-4]
            \arrow["f", two heads, from=3-2, to=2-1]
            \arrow[two heads, from=4-3, to=2-3]
            \arrow[two heads, from=4-3, to=2-4]
            \arrow["h", from=4-3, to=3-2]
        \end{tikzcd}
        \label{eq:ActionOnMorphismsOfCoproduct}
    \end{equation}
    The simplex $x'$, which $x$ restricts to along $k' \rightarrowtail k$, may be degenerate. So we do as before, and consider the unique nondegenerate simplex $x'_*$ which restricts to $x'$, as in the triangle on the right-hand side in the above diagram. The pair $(m' \twoheadrightarrow k'_*, x'_*)$ is the image of $(f,x)$. By definition it lives in the same equivalence class as $(fh,x)$ in the domain of Equation~\eqref{eq:sk_nd_last_step}. 
    By construction, this definition ensures the naturality of Equation~\eqref{eq:skNonDegFormulation} in $m$. 
    Similarly the action on morphisms $N\to N'$ is also forced on us by naturality (using a Diagram similar to \eqref{eq:naturalityOfSkInN}) and makes \eqref{eq:skNonDegFormulation} natural in $N$ by construction.

    All in all we have shown that any equivalence class in Equation~\eqref{eq:skEpiFormulation} has a unique representative with a nondegenerate simplex. As a result, we obtain
    \begin{equation}
        \begin{split}
            \te{sk}_n(N)(m,a) & \simeq \smashoperator[r]{\coprod_{\substack{f:m \twoheadrightarrow k\\ k\le n}}} \ N^{\te{nd}}(k,a),
        \end{split}
        \label{eq:skNonDegFormulation}
    \end{equation}
\end{enumerate}

The right-hand side in Eq.~\eqref{eq:skNonDegFormulation} is still a functor from $[(\bo{n}\times \bo{I})^{\te{op}},\bo{Set}]$ to $[(\bo{Delta}\times \bo{I})^{\te{op}},\bo{Set}]$. The formulation as a coproduct over surjections to non-degenerate simplices now additionally allows us to prove that, whenever $N\in\bo{nFuz}$, then $\te{sk}_n(N)$ is a fuzzy simplicial set, i.e.~$\te{sk}_n(N)(m)$ satisfies the injectivity and gluing condition for each $m$ in $\Delta$:

\textbf{Injectivity:} We must show that $\te{sk}_n(N)(m,b)\to \te{sk}_n(N)(m,a)$ is injective. Since a coproduct of a family of injections yields an injection, it suffices to show that $N^\te{nd}(m,b) \to N^\te{nd}(m,a)$ is injective. In diagram~\eqref{diag:nondegenerate_naturality} the vertical maps and the bottom horizontal map are injective, and thus $N^\te{nd}(m,b) \to N^\te{nd}(m,a)$ must be injective as well. 

\textbf{Gluing condition:} Let $B \subset I$ be a nonempty subset. We must show that $\lim_{b\in B}\te{sk}_n(N)(m,b)\simeq \te{sk}_n(N)(m,\sup B)$ for all $m\in\Delta$. 
Note that $B$ is connected as a category because we only consider non-empty sets in Def.~\ref{def:fuzzySet} and because it fulfills the property that any two objects in it can be connected by a zigzag (in this case an inclusion) of morphisms.
Since coproducts commute with connected limits in sets,\footnote{Let $S$ be a set, $B$ be a connected category, and $F : S\times B \to \te{Sets}$  a functor. The canonical map
\begin{equation} \label{eq:connected_lim_coproduct}
    \bigsqcup_{s \in S} \lim_{b \in B} F(s,b) \longrightarrow \lim_{b \in B} \bigsqcup_{s \in S} F(s,b)
\end{equation}
is given by $(s, \{ x_b \}_{b \in B}) \mapsto \{(s,x_b)\}_{b \in B} $, where $\{x_b\}_{b \in B}$ is a compatible family of elements $x_b \in F(s,b)$. The injectivity of this map is immediate. To see the surjectivity, consider a compatible family $\{(s_b, x_b)\}_{b\in B}$ in the codomain. 
The functorial action of $\bigsqcup_{s \in S}F(s,b)$ on some morphism $f : b \to b'$ maps an element $(s, x)$ with $x \in F(s,b)$ to $(s,F(\te{id}_s,f)(x)) $, where $F(\te{id}_s,f)(x) \in F(s,b')$. Compatibility means, for all maps $b \to b'$, the pair $(s_{b'}, x_{b'})$ is mapped to $(s_b, x_b)$, which, by the previous sentence forces $s_b = s_{b'}$. Since $B$ is connected, $s_b = s_{b'}$ for all $b,b' \in B$, which shows that $\{(s_b, x_b)_{b\in B} \}$ is in the image of Equation~\eqref{eq:connected_lim_coproduct}.
} 
we can conclude that
\begin{equation}
    \begin{split}
        \lim_{b\in B} \te{sk}_n(N)(m,b) \simeq  \lim_{b\in B}\smashoperator[r]{\coprod_{\substack{f:m \twoheadrightarrow k \\k\le n}}} \ N^{\te{nd}}(k,b) \simeq \coprod_{\substack{f:m \twoheadrightarrow k \\k\le n}}  \lim_{b\in B} N^{\te{nd}}(k,b).
    \end{split}
    \label{eq:commutationOfLimAndCoprod}
\end{equation}
It thus only remains to show that  $\lim_{b\in B} N^{\te{nd}}(k,b)\simeq N^{\te{nd}}(k,\sup B)$. To do this, consider the diagram
\begin{equation}
    \begin{split}
        \bt 
        N^{\te{nd}}(k, \sup B)\ar[tail]{rrr} \ar[tail]{d} \&\&\& \lim_{b\in B} N^{\te{nd}}(k,b)\ar[tail]{d} \\
        N(k,\sup B) \ar[tail]{rrr}{\simeq} \&\&\& \lim_{b\in B} N(k,b)
        \et
    \end{split}
\end{equation}
The bottom horizontal map is an isomorphism since $N\in\bo{nFuz}$, and the vertical maps are injective. By commutativity the top horizontal map is injective. To prove surjectivity, consider an element in $\lim_{b\in B} N^\text{nd}(k,b)$. Such an element is a sequence $\{x_b\}_{b\in B}$ such that $x_{b'}|_b = x_b$ for all $b\le b'$ in $B$. By injectivity we may view this family as an element in $\lim_{b \in B}N(k,b)$. The isomorphism then gives us an $x_{\sup B} \in N(k,\sup B)$ which restricts to $x_{b}$ for each $b \in B$. By Lemma \ref{lem:DegeneracyIsPreservedAcrossDifferentStrengths}, $x_{\sup B}$ lives in $N^\text{nd}(k, \sup B)$, and by the commutativity of the square it is in the preimage of $(x_b)_{b\in B} \in \lim_{b \in B} N^{\text{nd}}(k,b)$.

\bigskip

\textbf{\large{Coskeleton functor}}

For the right Kan-extension, we again start by applying Proposition \ref{prop:kanExtension} to obtain another adjunction from $[(\bo{n}\times \bo{I})^{\te{op}},\bo{Set}]$ to $[(\bo{Delta}\times \bo{I})^{\te{op}},\bo{Set}]$. In this case, we do not even have to compute the limit involved in this construction explicitly because we can use the following argument. Suppose that $\te{cosk}_n$ is right adjoint to $\te{tr}_n$. Then we must have 
\begin{equation}
    \begin{split}
        \te{Hom}_{\bo{PSh}(\Delta\times\bo{I})}(\te{tr}_n(S),N)\simeq \te{Hom}_{\bo{PSh}(\Delta\times\bo{I})}(S,\te{cosk}_n(N)).
    \end{split}
\end{equation}
In particular, this must be true for $S=\Delta^m_{a}$. But then, by
Yoneda,
\begin{equation}
    \begin{split}
        \te{Hom}_{\bo{PSh}(\Delta\times\bo{I})}(\Delta^m_{a},\te{cosk}_n(N))\simeq \te{cosk}_n(N)(m,a).
    \end{split}
\end{equation}
As a result, we define 
\begin{equation}
    \begin{split}
        \te{cosk}_n(N)(m,a) := \te{Hom}_{\bo{PSh}(\bo{n}\times\bo{I})}(\te{tr}_n(\Delta^m_{a}),N)
    \end{split}
    \label{eq:defCoskn}
\end{equation}
and one checks that this indeed defines a right adjoint. By uniqueness (up to isomorphism) of adjoints, we must have
\begin{equation}
    \begin{split}
        \te{Ran}_{I_n}N(m,a)\simeq \te{Hom}_{\bo{PSh}(\bo{n}\times\bo{I})}(\te{tr}_n(\Delta^m_{a}),N).
    \end{split}
\end{equation}

It remains to check that this adjunction can be restricted to an adjunction between $\bo{nFuz}$ and $\bo{sFuz}$, i.e.~that for each $N\in\bo{nFuz}$, we have that $\te{cosk}_n(N)(m,a)$ is a fuzzy simplicial set, which in turn means that the injectivity and gluing conditions are satisfied for each $m$.

\textbf{Injectivity:} For injectivity, we must show that
    \begin{equation}
        \begin{split}
            \te{cosk}_n(N)(\te{id}_m,i_{ab}): \te{Hom}_{\bo{PSh}(\bo{n}\times\bo{I})}(\te{tr}_n(\Delta^m_{b}),N) \to \te{Hom}_{\bo{PSh}(\bo{n}\times\bo{I})}(\te{tr}_n(\Delta^m_{a}),N),
        \end{split}
        \label{eq:coskInjCond}
    \end{equation}
    is injective for all $N\in\bo{nFuz}$. Since $\te{cosk}_n(N)(\te{id}_m,i_{ab})$ acts by precomposition with
    \begin{equation}
        \begin{split}
            \te{tr}_n(\bo{y}(\te{id}_m,i_{ab})):\te{tr}_n(\Delta^m_{a})\to \te{tr}_n(\Delta^m_{b}),
        \end{split}
        \label{eq:epiCond}
    \end{equation}
    this means we must show that for all $N\in\bo{nFuz}$ and for all natural transformations $f,g:\te{tr}_n(\Delta^m_{b})\to N$, if $f\circ \te{tr}_n(\bo{y}(\te{id}_m,i_{ab}))=g\circ \te{tr}_n(\bo{y}(\te{id}_m,i_{ab}))$, then $f=g$.

    Consider the natural transformation in Equation~\eqref{eq:epiCond}. Evaluated at $(k,c) \in \bo{n} \times I$, where $c \le a$, it reduces to the identity $\Hom_\Delta(k,m) \to \Hom_\Delta(k,m)$, which is epi (surjective). (When $ a < c \le b$, evaluation at $(k,c)$ reduces to $\emptyset \to \text{Hom}_\Delta(k,m)$, which is not surjective). Epiness then implies that $f(k,c) = g(k,c)$ for all $(k,c) \in \bo{n}\times I$ with $c \le a$.

    For all $ c > a $, consider the naturality square
    \begin{equation}
        \begin{split}
            \bt
                \te{tr}_n(\Delta^m_{b})(k,c) 
                \ar{rrr}{\te{tr}_n(\Delta^m_{b})(\te{id}_k,i_{ac})}
                \ar[']{d}{f(k,c)} 
                \&\&\& 
                \te{tr}_n(\Delta^m_{b})(k,a) 
                \ar{d}{f(k,a)} \\
                N(k,c) 
                \ar{rrr}{N(\te{id}_k,i_{ac})} 
                \&\&\& 
                N(k,a)
            \et
        \end{split}
        \label{eq:diagInjCosk}
    \end{equation}
    Since $N \in \bo{nFuz}$, the bottom horizontal map is injective. Since $f$ and $g$ agree at $(k,a)$, the commutativity of the diagram guarantees that they also agree at all $(k,c)$ for $ c > a $.

\textbf{Gluing condition:} Consider a subset $B \subset I$ with supremum $a \in I$. We use the shorthand $[-,-]$ for $\Hom(-,-)$. Recall that the set of natural transformations has a presentation in terms of ends: given any two functors $F,G : \mathcal C \to \mathcal D$,
    \begin{equation}
        [F,G] \simeq \int_{C \in \mathcal C} [FC,GC].
    \end{equation}
    Recall also that ends commute with limits,  and that $[-,-]$ commute with limits and ends in the second variable. For $N \in \bo{nFuz} $ we have
    \begin{align}
        \lim_{b \in B} \te{cosk}(N)(m,b) &\simeq \lim_{b \in B} [\te{tr}_n(\Delta^m_b), N] \\ \nonumber
        &\simeq \lim_{b \in B} \int_{\substack{k \in \bo{n} \\ c \in I}} [\te{tr}_n(\Delta^m_b)(k,c), N(k,c)] \\ \nonumber
        &\simeq \lim_{b \in B} \int_{\substack{k \in \bo{n} \\ c \in I}} \Big[ [k,m]\times [c,b], N(k,c)\Big] \\ \nonumber
        &\simeq \lim_{b \in B} \int_{\substack{k \in \bo{n} \\ c \in I}} \Big[ [k,m], \ \big[ [c,b], N(k,c) \big] \Big], \qquad \text{Cartesian closedness} \\ \nonumber
        &\simeq \lim_{b \in B} \int_{k \in \bo{n} } \int_{c \in I} \Big[ [k,m], \ \big[ [c,b], N(k,c) \big] \Big], \qquad \text{Fubini} \\ \nonumber
        &\simeq \lim_{b \in B} \int_{k \in \bo{n} } \Big[ [k,m], \int_{c \in I} \big[ [c,b], N(k,c) \big] \Big] \\ \nonumber
        &\simeq \lim_{b \in B} \int_{k \in \bo{n} } \Big[ [k,m], [\bo{y}(b), N(k,-)] \Big] \\ \nonumber
        &\simeq \lim_{b \in B} \int_{k \in \bo{n} } \Big[ [k,m], N(k,b) \Big], \qquad \text{Yoneda} \\ \nonumber
        &\simeq \int_{k \in \bo{n} } \Big[ [k,m], \lim_{b \in B} N(k,b) \Big] \\ \nonumber
        &\simeq \int_{k \in \bo{n} } \Big[ [k,m], N(k,a) \Big] \\ \nonumber
        &\simeq \te{cosk}(N)(m,a)
    \end{align} 
    where in the last line we traced all the steps from the beginning until Yoneda in reverse.

This completes the proof. \hfill$\blacksquare$

\subsection{Proof of Proposition \ref{cor:trcoskId}}
\label{proof:cor:trcoskId}

Since $ I_n : \bo{n} \to \Delta $ is fully faithful (by definition), by abstract nonsense the unit of the adjunction $ \te{sk}_n \dashv \te{tr}_n $ and the counit of $ \te{tr}_n \dashv \te{sk}_n $ are isomorphisms (cf.~\cite[Corollary X.3.3]{MacLane78}). By Proposition \ref{prop:LeftAdjointFullyFaithfulUnitIso}, this is equivalent to saying that $ \te{sk}_n $ and $ \te{cosk}_n $ are fully faithful. Next we give a more explicit proof, justifying the last equation in the proposition.

We begin with the left Kan-extension for which we construct an explicit isomorphism, starting from the formulation given in Eq.~\eqref{eq:LanSkeletonFormula}. Assume $ m \le n $ and let $ (\sigma, x) \in \coprod_{\sigma : I_n(k,b) \to (m,a) } N(k,b) $ where $ \sigma : I_n(k,b) \to (m,a) $ is a morphism in $ (\Delta \times I)^{\te{op}} $ and $ x \in N(k,b) $. Since $ m \le n $, by fully faithfulness of $ I_n $, the morphism $ \sigma $ lives in the subcategory $ (\bo{n}\times I)^{\te{op}} $ and by application of $ N $ we get a morphism $ N(\sigma) : N(k,b) \to N(m,a) $. This defines a map
\begin{align}
    \tilde \phi : \left( \smashoperator[r]{\coprod_{\sigma : I_n(k,b) \to (m,a)}} \ N(k,b) \right) \to N(m,a); \quad (\sigma, x) \mapsto N(\sigma)(x).
\end{align}
Now, suppose we have two elements in the coproduct such that $ (\sigma, x) \sim(\sigma', x')$ as in equation \eqref{eq:LanSkeletonFormula}. This means there exists a $ \gamma : (k,b) \to (k',b') $ such that $ \sigma' \circ I_n(\gamma) = \sigma $, where all morphisms again live in the subcategory $ (\bo{n}\times I)^{\te{op}} $. Applying $ N $ to this diagram yields $ \tilde \phi(\sigma,x) = \tilde \phi(\sigma',x') $, thus $ \tilde \phi $ factors through $ \te{Lan}_{I_n} N(m,a) $ and this gives a candidate for the map $ \phi : \te{Lan}_{I_n} N(m,a) \to N(m,a) $.

The map $\phi $ is surjective, since for each $ x \in N(m,a) $, we have $ (\te{id}_{(m,a)}, x) \in \te{Lan}_{I_n}(m,a) $ in its preimage. To show that it is injective, we show that, whenever we assume $\phi(\sigma,x)=\phi(\sigma',x')$, then we must have $(\sigma,x)\sim(\sigma',x')$. To this end, note that $(\sigma,x)$ is in the same equivalence class as $(\te{id}_{(m,a)},N(\sigma)(x))$ (to see this, choose $\gamma=\sigma$ in eq.~\eqref{eq:equivRelationSk}).
As a consequence,
\begin{equation}
    \begin{split}
        (\sigma,x) &\sim (\te{id}_{(m,a)},N(\sigma)(x)) \ov{by definition of $\phi$} (\te{id}_{(m,a)},\phi(\sigma,x)) \\
        & \ov{by assumption} (\te{id}_{(m,a)},\phi(\sigma',x')) = (\te{id}_{(m,a)},N(\sigma')(x')) \sim (\sigma',x').
    \end{split}
\end{equation}

As for the right Kan extension, this follows directly from applying the Yoneda lemma to the formula given in Proposition \eqref{prop:leftRightAdjointOfTruncation}.
\hfill$\blacksquare$

\section{Proofs of section \ref{sec:Adjunction}}

\subsection{Proof of Proposition \ref{prop:KanExtension}}
\label{proof:prop:KanExtension}

The assumption that $\te{Re}$ is a functor, and that $\bo{M}$ has all small colimits, allows us to use the theorem of Kan extensions (Prop. \ref{prop:kanExtension}) that asserts that \eqref{eq:ReExtended} defines a functor because we can specialize the functors in the Kan extension theorem to the following ones
\begin{equation}
    \begin{split}
      (F &=\bo{y}):(\bo{A}=\bo{C})\to (\bo{B}=\bo{PSh}(\bo{C})), \\
      (X&=\text{Re}\circ \bo{y}) : (\bo{A}=\bo{C}) \to (\bo{S}=\bo{M})
    \end{split}
    \label{eq:choicesForExtensionGeneral}
\end{equation}
The exact definition of $\te{Re}$ does not matter as long as it is any functor between $\bo{y}(\bo{C})$ and $\bo{M}$.

Next, we use the following series of natural isomorphisms to show that 
\begin{equation}
    \begin{split}
        \text{Hom}_{\bo{M}}(\text{Re}(S),Y)\simeq \text{Hom}_{\bo{PSh}(\bo{C})}(S,\text{Sing}(Y))
    \end{split}
\end{equation}
for any $S\in \bo{PSh}(\bo{C})$ and $Y\in\bo{M}$.
\begin{equation*}
    \begin{split}
    &\text{Hom}_{\bo{M}}(\text{Re}(S),Y)
    \ovs{\ref{prop:density}}
    \text{Hom}_{\bo{M}}(\text{Re}(\text{colim}(\bo{y}(P_S)),Y)\\
    &\ov{def}
        \text{Hom}_{\bo{M}}(\text{colim}(\text{Re}(\bo{y}(P_S)),Y)
        \ovs{\ref{prop:RepsPreserveLimits}} \text{lim}[\text{Hom}_{\bo{M}}(\text{Re}(\bo{y}(P_S),Y)]\\
        &\ov{def}  \text{lim} [\text{Sing}(Y)(P_S) ]
        \ovs{Yoneda} \text{lim}[\text{Hom}_{\bo{PSh}(\bo{C})}(\bo{y}(P_S),\text{Sing}(Y))] \\
        &\ovs{\ref{prop:RepsPreserveLimits}}
        \text{Hom}_{\bo{PSh}(\bo{C})}(\text{colim}(\bo{y}(P_S)),\text{Sing}(Y)) \\
        &\ovs{\ref{prop:density}}
        \text{Hom}_{\bo{PSh}(\bo{C})}(S,\text{Sing}(Y)).
    \end{split}
\end{equation*}

\hfill$\blacksquare$

\subsection{Proof of Proposition \ref{prop:spivakReNonExpansive}}
\label{proof:prop:spivakReNonExpansive}

$\te{Re}^p$ preserves identities because
\begin{equation}
    \begin{split}
        \te{Re}^p(\bo{y}(\te{id}_n,i_{aa}))(x) &= 
        \left( \parallel x_{\te{id}_n^{-1}(0)} \parallel_p ,~\cdots~,\parallel x_{\te{id}_n^{-1}(n)} \parallel_p\right)\\
        &= \left( \parallel x_{0} \parallel_p ,~\cdots~,\parallel x_{n} \parallel_p\right)= \left( x_{0},~\cdots~,x_{n}\right).
    \end{split}
\end{equation}
Now let $f:l\to m$, $g:m\to n$ be morphisms in $\Delta$ and $i_{ab}:a\to b$ and $i_{bc}:b\to c$ morphisms in $\bo{I}$. Then $\te{Re}^p$ preserves their composition because
\begin{equation}
    \begin{split}
        \te{Re}^p(\bo{y}(g,i_{bc}))\circ \te{Re}^p(\bo{y}(f,i_{ab}))(x) &= 
        \te{Re}^p(\bo{y}(g,i_{bc}))\underbrace{( \parallel x_{f^{-1}(0)} \parallel_p ,~\cdots~,\parallel x_{f^{-1}(m)} \parallel_p)}_{=:(y_1,\cdots,y_m)}\\
        &= \left( \parallel y_{g^{-1}(0)} \parallel_p ,~\cdots~,\parallel y_{g^{-1}(n)} \parallel_p\right)
    \end{split}
\end{equation}
and each term in turn is of the form 
\begin{equation}
    \begin{split}
        \parallel y_{g^{-1}(k)} \parallel_p &= \left(\sum_{a\in g^{-1}(k)}y_a^p\right)^{1/p}
        = \left(\sum_{a\in g^{-1}(k)}\left[\left(\sum_{b\in f^{-1}(a)}x_b^p\right)^{1/p}\right]^p\right)^{1/p}  \\
        &= \left(\sum_{a\in g^{-1}(k)}\sum_{b\in f^{-1}(a)}x_b^p\right)^{1/p} =  \left(\sum_{b\in f^{-1}(g^{-1}(k))}x_b^p\right)^{1/p}
        = \parallel x_{(g\circ f)^{-1}(k)} \parallel_p,
    \end{split}
\end{equation}
yielding
\begin{equation}
    \begin{split}
        \te{Re}^p(\bo{y}(g,i_{bc}))\circ \te{Re}^p(\bo{y}(f,i_{ab}))(x) &=\te{Re}^p(\bo{y}(g\circ f,i_{bc}\circ i_{ab})).
    \end{split}
\end{equation}

To see that it is non-expansivene, note that this is immediate for face maps. For degeneracies, it suffices to consider the case $ a = b $, since $ \log b / \log a \le 1 $. Furthermore, due to cancellations on both sides, it suffices to consider the case $ s : [1] \to [0] $. In this case, the task is reduced to showing that
\begin{equation}
  \parallel x - y \parallel_p \, \ge \Big| \parallel x \parallel_p - \parallel y \parallel_p \Big|,
\end{equation}
which holds since $ \Vert - \Vert_p$ is a norm.

\hfill$\blacksquare$

\subsection{Proof of Proposition \ref{prop:singFactorsThroughSfuz}}
\label{proof:prop:singFactorsThroughSfuz}

We now check that for a fixed $n \in \Delta$, the functor $\text{Sing}_n$ is a fuzzy set. First of all, all the functions $\text{Re}^p(\bo{y}(\text{id}, i_{ab}))$ are surjective and therefore epi. As a result, precomposition by these maps is a monomorphism (i.e.~injective). Lastly we must show that canonical map
\begin{equation}
  \Hom(\text{Re}^p(\Delta^n_{a}), Y) \to \lim_{b \in B} \Hom(\text{Re}^p(\Delta^n_{b}), Y)
\end{equation}
is a bijection for every non-empty set $B \subset [0,1]$ with supremum $a$.
When $B$ is the singleton, then this is automatic. So assume $B$ is not a singleton set.
An element in the image is a sequence of function $(f_b : \text{Re}^p(\Delta^n_b) \to Y)$ such that for $ 0 < b' \le b $ in $\bo{I}$ we have $f_{b'}(x) = f_b\left(x\right) $. Given any fixed nonzero $b_\ast \in B$, define $ f_a : \text{Re}^p(\Delta^n_a) \to Y$ by $ f_a(x) : = f_{b_\ast}\left( x\right) $. In fact, this definition is forced, and by construction $f_a$ restricts to $f_b$ for all $b \in B$. Moreover, $f_a$ is non expansive. To see this we compute
\begin{align}
  d_Y (f_a (x), f_a (y)) &= d_Y \left( f_b\left(x \right), f_b\left( y \right) \right) \\ \nonumber
  &\le d_b \left( (x), (y) \right) \\ \nonumber
  &= -\log b \parallel  x - y \parallel_p \\ \nonumber
  &= \frac{\log b}{\log a}  (-\log a) \parallel x - y \parallel_p \\ \nonumber
  &= \frac{\log b}{\log a} d_a(x,y).
\end{align}
Thus $ \frac{\log b}{\log a} d_a(x, y) \ge d_Y(f_a (x), f_a ( y))$. But this is true for all $b \in B$ and thus also when we take the supremum over $B$. Since $\log$ is monotone, this proves the claim.

\hfill$\blacksquare$

\subsection{Proof of Proposition \ref{prop:umapReIsFunWithCodUM}}
\label{proof:prop:umapReIsFunWithCodUM}

Note that the map defined in \eqref{eq:reDeltaMorphisms} is non-expansive because $i_{ab}$ exists iff $a\le b$, which in turn implies that $-\log(a)\ge -\log(b)$ which implies, for $i\ne j$ and $\sigma(i)\ne \sigma(j)$, that
\begin{equation}
    \begin{split}
        d_b(\te{Re}^{\te{skeleton}}(\sigma,i_{ab})(x_i,x_j))&=d_b(x_{\sigma(i)},x_{\sigma(j)})=-\log(b)\\&\le -\log(a)=d_a(x_i,x_j),
    \end{split}
\end{equation}
while for $i=j$ or $\sigma(i)=\sigma(j)$, the inequality holds trivially as well. Since $\te{Re}^{\te{skeleton}}:\bo{y}(\Delta\times\bo{I})\to\bo{UM}$ preserves composition and identities, it is a functor.
The proof for $\bo{EPMet}$ is analogous.
\hfill$\blacksquare$

\subsection{Proof of Proposition \ref{prop:skeletonSingFactorsThroughSfuz}}
\label{proof:prop:skeletonSingFactorsThroughSfuz}

We must check, that for a fixed $n \in \Delta$, the $\text{Sing}_n$ functor is a fuzzy set. As set maps, the functions $ \te{Re}^{\te{skeleton}}(\text{id}_n, i_{ab}) $ are bijective, and therefore epi as non-expansive maps of metric spaces. It follows that precomposition by $ \te{Re}^{\te{skeleton}}(\text{id}_n, i_{ab}) $ is mono (i.e. injective). Next we must show that the canonical map
\begin{equation}
  \Hom(\text{Re}^p(\Delta^n_{a}), Y) \to \lim_{b \in B} \Hom(\text{Re}^p(\Delta^n_{b}), Y)
\end{equation}
is bijective for every non-empty set $B \subset [0,1]$ with $\sup(B) = a$.
An element on the right-hand side is a sequence of function $(f_b)_{b \in B}$ such that $f_b(x) = f_{b'}(x)$ for all $ b' \le b \in B $. Given some fixed $b_\ast \in B$, define $f_a : \text{Re}^p(\Delta^n_{a}) \to Y $ by setting $f_a(x) = f_{b_\ast}(x) $. Then, by construction $f_a$ restricts to $f_b$ for all $b \in B$. Moreover, $f_a$ is non-expansive, for, given any two $x_i \not = x_j$ we have
\begin{align}
  d_Y(f_a(x_i), f_a(x_j)) &= d_Y(f_b (x_i), f_b(x_j)) \\ \nonumber
  &\le d_b( x_i, x_j) \\ \nonumber
  &= -\log b.
\end{align}
Since this is true for all $b \in B$, it must be true for the supremum,~i.e.~$d_y(f_a(x_i), f_b(x_j)) \ge -\log a = d_a(x_i, x_j) $.
\hfill$\blacksquare$

\subsection{Proof of Proposition \ref{prop:metricRealizationOfIdentity}}
\label{proof:prop:metricRealizationOfIdentity}

We only show the proof for the skeleton adjunction because the proofs for Spivak's adjunction is analogous.

As a first step, we check what the metric realization functor does to diagrams of the form 
\begin{equation}
    \begin{split}
        \bt
        \begin{pmatrix}
            m,a\\
            s
        \end{pmatrix}
          \ar[yshift=0.7em]{r}{(\te{id},i_{ab})}\& \ar[yshift=-0.85em]{l}{S(\te{id},i_{ab})} 
          \begin{pmatrix}
              m,b\\
              s'
          \end{pmatrix}
        \et
    \end{split}
    \label{eq:34}
\end{equation}
(Note that since $S(\te{id},i_{ab})$ is injective, $s'$ can be identified with $s$ in $S(m,a)$, i.e.~up to inclusion $s'$ equals $s$.)
Applying $D_S$, defined in eq.~\eqref{eq:ReExtended}, results in $D_S((m,a),s)=\re \circ \bo{y}(m,a)=\re(\Delta^m_{a})$,
and similarly $D_S((m,b),s')=\re(\Delta^m_{b})$. We also have 
\begin{equation}
    \begin{split}
        D_S 
        \begin{pmatrix}
            \bt \ar[yshift=0.7em]{r}{(\te{id},i_{ab})}\& \ar[yshift=-0.85em]{l}{S(\te{id},i_{ab})} \et
        \end{pmatrix}(x_i\in \re(\Delta^m_{a}))
        &= \re \circ \bo{y}(\te{id},i_{ab})(x_i) \ov{\ref{eq:reDeltaMorphisms}}
        x_{\te{id}(i)}=x_i.
    \end{split}
\end{equation}
As a consequence,
\begin{equation}
    \begin{split}
        D_S((\te{id},i_{ab}),S)(\re(\Delta^m_{a}))=\re(\Delta^m_{b}).
    \end{split}
    \label{eq:identi}
\end{equation}
If we now take the colimit of the diagram 
\begin{equation}
    \begin{split}
        \bt[huge]
        \re(\Delta^m_{a})\ar{r}{D_S((\te{id},i_{ab}),S)}\& \re(\Delta^m_{b}),
        \et
    \end{split}
    \label{eq:colimDiagSimpleSimplices}
\end{equation}
which is the image under $D_S$ of \eqref{eq:34}, we obtain, by Prop. \ref{prop:colimitUM}, a set 
\begin{equation}
    \begin{split}
        X/\sim &= F(\re(\Delta^m_{a}))\cupdot F(\re(\Delta^m_{b})) / \sim,\\
        \te{where }&x\sim x'\te{ iff }x'=FD_S((\te{id},i_{ab}),S).
    \end{split}
\end{equation}
And due to \eqref{eq:identi}, we have $X/\sim ~\simeq F(\re(\Delta^m_{b}))$,
while the metric is given by $d(x,y)=d_b(x,y)$,
where $d_b$ is the metric on $\re(\Delta^m_{b})$ (the one defined in
\eqref{eq:reDeltaUmap} where the different vertices of the simplex have
  distance $-\log b$) because distances on $\re(\Delta^m_{b})$ are shorter than distances on $\re(\Delta^m_{a})$ and since points of those two spaces are identified, the infimum sequence of eq.~\eqref{eq:simMetric} is always going to be along a path in $\re(\Delta^m_{b})$. Hence, intuitively, for paths whose points are identified, distances are always traveled along the path that has maximal strength. Thus, we can conclude $\te{colim}(\te{diagram}\eqref{eq:colimDiagSimpleSimplices}) \simeq \re(\Delta^m_{b})$.

This has the important consequence that, if we have an arbitrarily long sequence of the form 
\begin{equation}
    \begin{split}
        \bt\cdots 
        \ar[yshift=0.7em]{r}{(\te{id},i_{ca})}\& \ar[yshift=-0.85em]{l}{S(\te{id},i_{ca})} 
        \begin{pmatrix}
            m,a\\
            s
        \end{pmatrix}
          \ar[yshift=0.7em]{r}{(\te{id},i_{ab})}\& \ar[yshift=-0.85em]{l}{S(\te{id},i_{ab})} 
          \begin{pmatrix}
              m,b\\
              s'
          \end{pmatrix} 
        \et
    \end{split}
\end{equation}
in $\bo{El}(S)$, then the metric realization functor collapses all of it to the realization of only one simplex, namely $\re(\Delta^m_{b})$. This completes the proof of the first statement of the proposition.

As a next step, we show the second statement, namely that ``the metric realization of diagrams of the form \eqref{eq:degeneracyDiagram} is also isomorphic to $\te{Re}(\Delta^m_{b})$''. The degeneracy map $\sigma:m+1\to m$ identifies two points $i,j$ of $m+1$ with one point $\sigma(i)=\sigma(j)$ of $m$ and consequently $D_S(\sigma,i_{bb})$ identifies two points $x_i$ and $x_j$ in $\te{Re}(\Delta^{m+1}_{b})$ with one point $x_{\sigma(i)}=x_{\sigma(j)}$ in $\te{Re}(\Delta^{m}_{b})$. Recall that the colimit then forces identified points to be in the same equivalence class. This takes place for pairs $(\sigma,i_{aa})$ for every $a\le b$ but a simplex is realized with its maximal strength as shown above and the maximal strength of the degenerate simplex must equal the strength of the simplex by Prop. \ref{prop:simplexFaceStrength}, and hence the identification simply leads to the realization of a metric space isomorphic to $\te{Re}(\Delta^{m}_{b})$. This completes the proof of the proposition.
\hfill$\blacksquare$

\subsection{Proof of Proposition \ref{prop:faceMap}}
\label{proof:prop:faceMap}

The face map $\delta:m\to m+1$ maps $i\in m$ to $\delta(i)\in m+1$ and is injective. Consequently, by definition, $D_S(\delta,i_{bb})$ maps $x_i\in \te{Re}(\Delta^m_{b})$ to $x_{\delta(i)}\in \te{Re}(\Delta^{m+1}_{b})$ and in this way injectively includes all elements of $\te{Re}(\Delta^m_{b})$ into $\te{Re}(\Delta^{m+1}_{b})$. At the same time, the morphism $(\te{id},i_{bd})$ includes, for every $i\in m$, $x_i$ into $\te{Re}(\Delta^m_{d})$. The colimit of the diagram then identifies all points that are mapped to each other and hence, for every $i\in m$, $x_{\delta(i)}\in \te{Re}(\Delta^{m+1}_{b})$ is identified with a point in $\te{Re}(\Delta^m_{d})$, where distances between different points are shorter, namely equal to $-\log(d)$.In the case of $\bo{UM}$, the colimit defines distances to be infimum distances among path over all identified points and therefore we obtain formula \eqref{eq:facemaprealization}. 

In the case of $\bo{EPMet}$, there might happen a gluing along an edge, where $d_a(x_i,x_j)=\infty$, while $d_b(x_i,x_j)<\infty$ but even then the glued edge will have finite distance because $\inf_{y_i\in [x_i],y_j\in [x_j]}< \infty$. Hence, the conclusion also holds in $\bo{EPMet}$.
\hfill$\blacksquare$

\subsection{Proof of Proposition \ref{prop:realizationOfFuzzySimpSet}}
\label{proof:prop:realizationOfFuzzySimpSet}

In the category $\Delta$, every simplex comes equipped with face maps for all of its faces. Consequently, those face maps are mapped to morphisms of the form \eqref{eq:faceMapDiagram} in the category of elements $\bo{El}(S)$. 

Their metric realization, as proven in Proposition \ref{prop:faceMap}, then simply consists of a simplex in which the distances between the points of any face are $(-\log)$ of the maximal strength of the face. However, the same is true for any faces of that face and so forth. This procedure can be iterated until reaching $1$-faces (edges). Therefore, the metric realization of any $n$-simplex $s$ is an $n$-point space in which the distance between any two points $x_1,x_2$ is given by $(-\log)$ of the maximal strength $\xi_1(s')$ of the $1$-face $s'$, whose realization $\Delta^1_{a}$ is identfied with the two points $x_1,x_2$ of the simplex $\Delta^n_{a}$ (for every $a$ lower or equal to the maximal strength of $s$) in the colimit. 

Recall that, if a $1$-face $s'$ has a realization $\Delta^1_{a}$ whose underlying set is identified with those two points $x_1$ and $x_2$ in the colimit, then there exist face maps $\delta_{1/2}:0\to 1$ that correspond to functions $S(\delta_{1/2},i_{00})$, and those functions fulfill the properties that $S(\delta_{1},i_{00})(s')=x_1$ and $S(\delta_{2},i_{00})(s')=x_2$ or that $S(\delta_{2},i_{00})(s')=x_1$ and $S(\delta_{1},i_{00})(s')=x_2$. (The symmetry here holds because directions are not remembered in the realization of the metric space). Furthermore, if there are multiple edges that fulfill those properties, 
then the distances defined in \eqref{eq:simMetric2} and \eqref{eq:simMetricEPMet2} according to Proposition \ref{prop:colimitUM} ensures that the one with smallest length is chosen for the realization. (As remarked before in the proof of Proposition \ref{prop:faceMap}, this also holds in $\bo{EPMet}$, where the distance \eqref{eq:simMetricEPMet2} is employed, because we are here only concerned with local identifications coming from face maps (because we only want to show that $1$-simplices suffice), the gluing of which does not require to traverse a path.) Hence, for the realization of a set of edges, whose realizations (or more precisely, their underlying sets) are all identified with the same two points, the distance between those points is given by \eqref{eq:nSimplexRealizationExplicit}. 

However, when multiple $n$-simplices (and in particular, multiple $1$-simplices) are realized all at the same time, then it might be the case that there exists a path from one point $x$ to another point $y$ that is shorter than any direct edge distance from $x$ to $y$. (Here, a path from $x$ to $y$ means that there exist edges $s_1,\cdots,s_{n-1}$ such that $x=x_1:=S(\delta_{1},i_{00})(s_1)$ and $x_i=S(\delta_{1},i_{00})(s_i)=S(\delta_{2},i_{00})(s_{i-1})$ for all $i\in\{1,\cdots,n-1\}$ and $y=x_n:=S(\delta_{2},i_{00})(s_{n-1})$. Being shorter then means that $d(x,y)$ defined in \eqref{eq:realizationExplicit} returns a lower value than $d_{\te{min}}(x,y)$.) 
Since the distance of two points of a realization of a $1$-face is always lower or equal to the distance obtained by  realizing a higher simplex of which those points are a part, it is sufficient to check, in the case of $\bo{M}=\bo{EPMet}$, if the distance of all realized $1$-faces of two points, that are identified, is infinity, and, if not or if one considers $\bo{M}=\bo{UM}$, to    
measure lengths along paths of $1$-edges when being interested in the infimum distance over all paths between the two points.
As a consequence, the infimum distance defined in \eqref{eq:simMetric} ensures that the final result of the realization has the distance defined in \eqref{eq:realizationExplicit}.
\hfill$\blacksquare$

\subsection{Proof of Proposition \ref{prop:SingFun}}
\label{proof:prop:SingFun}

Recall that the functor $ \te{Sing}:\bo{M}\to\bo{sFuz} $ (where $\bo{M}$ is either $\bo{UM}$ or $\bo{EPMet}$) is defined by 
\begin{equation*}
    \begin{split}
        \te{Sing}(X,d)(m,a) = \te{Hom}_{\bo{M}}(\te{Re}^{\te{skeleton}}(\Delta_{a}^n),(X,d)).
    \end{split}
\end{equation*}
where
\begin{equation*}
    \begin{split}
        \te{Re}^{\te{skeleton}}:~&\bo{y}(\Delta\times \bo{I})\to \bo{M},
        \quad \Delta_{a}^n \mapsto (\{x_0,\cdots,x_n\},d_a),\\
        &\te{ where }\quad d_a(x_i,x_j) := \begin{cases}
            -\log(a),&\te{if }i\ne j,\\
            0,&\te{else.}
        \end{cases}
    \end{split}
    \label{eq:reDeltaUmap}
\end{equation*}
To improve the readability of this proof we adopt the notation $ S_X $ as a shorthand for $ \te{Sing}(X,d) $. Similarly, given a morphism of uber metric spaces $ f : (X,d) \to (Y,e) $ we write $ S_f $ instead of $ \te{Sing}(f) $.

Now a morphism $ \sigma : \te{Re}^{\te{skeleton}}(\Delta_{a}^n) \to (X,d) $ in $\bo{M}$, i.e.~a non-expansive map, is the same thing as a tuple $ (r_0, \dots, r_n) \in X\times \cdots \times X $ such that $ d(r_i, r_j) \le d_a(x_i, x_j) = - \log(a) $ for all $ i \neq j $. With this in mind, define the set
\begin{align*}
    [X,d]^n_a := \{ (r_0, \dots, r_n) \in X\times \cdots \times X \, | \, d(r_i, r_j) \le - \log(a) \ \ \forall i \neq j \}.
\end{align*}
which, by construction, is isomorphic to $ S_X(n,a) $. For a fixed uber metric space $ (X,d) $ we can turn this into a fuzzy simplicial set by defining its action on a morphism $ (\phi, i_{ab}) : (n,a) \to (m, b) $ in $ \Delta \times \bo{I} $ to be the map $ [X,d]^m_b \to [X,d]^n_a; \ (r_0, \dots, r_m) \mapsto (r_{\phi(0)}, \dots, r_{\phi(n)}) $, making the following diagram commute,
\begin{equation} \label{diag:SingBracketIsoSquare}
    \begin{tikzcd}
        {S_X(m,b)} \arrow[d, "{S_X(\phi, i_{ab})}"'] \arrow[r, "\cong"] & {[X,d]^m_b} \arrow[d, "{[X,d](\phi, i_{ab})}"] \\
        {S_X(n,a)} \arrow[r, "\cong"'] & {[X,d]^n_a}          
    \end{tikzcd}
\end{equation}
giving rise to a natural isomorphism $ \te{Sing}(X,d) \cong \{X,d\} $. Given a morphism of uber or extended-pseudo metric spaces $ f : (X,d) \to (Y,e) $, the commutativity of the following diagram
\begin{equation} \label{diag:SingBracketNatIso}
    \begin{tikzcd}
        {S_X(n,a)} \arrow[d, "{S_f(n,a)}"'] \arrow[r, "\cong"] & {[X,d]^n_a} \arrow[d, "f\times \cdots \times f"] \\
        {S_Y(n,a)} \arrow[r, "\cong"']                  & {[Y,e]^n_a}                                     
    \end{tikzcd}
\end{equation}
for all $ (n,a) \in \Delta \times \bo{I} $ shows that we can in fact make the construction $ [X,d] $ functorial in $ (X,d)$, giving rise to a functor $ [-] : \bo{M} \to \bo{sFuz} $, naturally isomorphic to $ \te{Sing} $.

Now, the strength of an $n$-simplex $ s $ of some fuzzy simplicial set $ S $ is by definition the supremum of all $ a' \in \bo{I} $ such that $ s \in S(n,a') $, cf.~Def.~\ref{def:classicalFuzzySets} and eq.~\eqref{eq:maxStrength}. Since the condition $ d(r_i,r_j) \le - \log(a) $ is equivalent to $ e^{-d(r_i, r_j)} \ge a $, it follows that the strength of an $n$-simplex of $ \te{Sing}$, seen as a tuple $ (r_0, \dots, r_n) $ is 
$ \min_{0 \le i,j \le n}\exp(-d(r_{i_j},r_{i_k})) $, 

which indeed is greater than or equal to $ a $. The last statement follows from the fact that the conditions $ d(r_i,r_j) \le - \log(a) $ are equivalent to $ \min_{0 \le i,j \le n}\exp(-d(r_{i_j},r_{i_k})) \ge a $.
\hfill$\blacksquare$

\subsection{Proof of Proposition \ref{prop:UMisInSFuz2}}
\label{proof:prop:UMisInSFuz2}

Let $ f, g : (X,d) \to (Y,e) $ be two maps of uber or extended-pseudo metric spaces such that $ \te{Sing}(f) = \te{Sing}(g) $. In particular $ \te{Sing}(f)^0_a = \te{Sing}(g)^0_a $. Recalling that $ [X,d]^0_a = X $ and similarly for $ [Y,e] $,  it follows from the commutative square \eqref{diag:SingBracketNatIso} for $ n=0 $ that $ f = g :  X = [X,d]^0_a \to [Y,e]^0_a = Y $. In particular $ f = g : (X,d) \to (Y,e) $ as maps of uber or extended-pseudo metric spaces. This shows that $ \te{Sing} $ is faithful.

As for the surjectivity, assume now that we are given an arbitrary map of fuzzy simplicial sets $ \phi : \te{Sing}(X,d) \to \te{Sing}(Y,e) $. This is equivalent to giving a map $ [X,d] \to [Y,e] $ which by abuse of notation we denote again by $\phi $. Evaluation at $ (0,a) $ yields a map $ \phi^0_a : [X,d]^0_a = X \to Y = [Y,e]^0_a $. Notice that this map is independent of the choice of $a$. This is because in the square \eqref{diag:SingBracketIsoSquare}, if $ (\te{id}_0,i_{ab}) : (0,a) \to (0,b) $, then the left vertical map $ \phi_a^0 $ and thereby also the right vertical map become the identity in the following diagram
\begin{equation*}
    \begin{tikzcd}
        X \arrow[r, equal] & {[X,d]^0_a} \arrow[rr, "{[X,d](\te{id}_0,i_{ab})}", equal] \arrow[d, "\phi^0_a"'] && {[X,d]^0_{b}} \arrow[d, "\phi^0_{b}"] \arrow[r, equal] & X \\
        Y \arrow[r, equal] & {[Y,e]^0_{a}} \arrow[rr, "{[Y,e](\te{id}_0,i_{ab})}"', equal]                    && {[Y,e]^0_{b}} \arrow[r, equal]                          & Y
    \end{tikzcd}
\end{equation*}
the commutativity of which is guaranteed by the fact that $ \phi $ is a simplicial map. We take $ f := \phi^0_a : X \to Y $ to be our candidate for the preimage of $ \phi $. Consider the diagram
\begin{equation*}
    \begin{tikzcd}
        X & {[X,d]^0_a} \arrow[d, "f"'] \arrow[l, equal] && {[X,d]^n_a} \arrow[ll, "{[X,d](g_k, \te{id}_a)}"'] \arrow[d, "\phi^n_a"] & {(r_0,\dots r_n)} \arrow[d, maps to] \\
        Y & {[Y,e]^0_a} \arrow[l, equal]                 && {[Y,e]^n_a} \arrow[ll, "{[Y,e](g_k,\te{id}_a)}"]                         & {(f(r_0),\dots, f(r_n))}             
    \end{tikzcd}
\end{equation*}
which we now explain. Here $ g_k $ is an iterated face map, picking out the $k$-th 0-simplex. In this case $ [X,d](g_k, \te{id}_a) $ is the $k$-th projection. The commutativity of the diagram shows that $ \phi_a^n $ must act as depicted on the right of the diagram, i.e.~$ \phi^n_a = f\times \dots \times f $. 

It remains to show that $ f $ is non-expansive. For, if this is the case, then $ f $ is indeed in the preimage of $ \phi $, i.e. $ [f] = \phi $. To see that $ f $ is non-expansive, consider a pair $ (r_0, r_1) \in [X,d]^1_a $ as well as its image $ \phi^1_a(r_0,r_1)=(f(r_0), f(r_1)) $. Then by an argument similar to that in \eqref{prop:simplexFaceStrength}, $ \te{str}(r_0,r_1) \le \te{str}(f(r_0), f(r_1)) $. Now, by definition $ (r_0,r_1) \in [X,d]^1_a $ means $ d(r_0,r_1) \le -\log(a) $, which may be interpreted as saying that $(r_0,r_1) $ has strength at least $ a $. As we allow $ a $ to vary, this becomes an equality when we reach the actual strength, that is $ d(r_0,r_1) = -\log(\te{str}(r_0,r_1)) $. Another way of seeing this, is by recalling that the strength of an $n$-simplex is by definition the largest $ a \in \bo{I} $ such that $ [X,d]^n_a $ contains this $n$-simplex, cf.~Def.~\ref{def:classicalFuzzySets} and eq.~\eqref{eq:maxStrength}. Using the corresponding equality for $ (f(r_0),f(r_1)) \in [Y,e]^1_a $ and putting everything together we get
\begin{align}
    d(r_0,r_1) = -\log(\te{str}(r_0,r_1)) \ge - \log(\te{str}(f(r_0), f(r_1)) ) = e(f(r_0), f(r_1))    
\end{align}
which completes the proof.
\hfill$\blacksquare$

\subsection{Proof of Corollary \ref{prop:UMisInSfuz}}
\label{proof:prop:UMisInSfuz}

This follows directly by applying Proposition \ref{prop:LeftAdjointFullyFaithfulUnitIso} to Proposition \ref{prop:UMisInSFuz2}.
\hfill$\blacksquare$

\subsection{Proof of Proposition \ref{prop:SingSimpSetProperties}}
\label{proof:prop:SingSimpSetProperties}

The fact that the first two properties hold for $ \te{Sing} $ can be read off directly from the isomorphic functor
\begin{align}
    [X,d]^n_a := \{ (r_0, \dots, r_n) \in X^{n+1} \, | \, d(r_i, r_j) \le - \log(a) \ \ \forall i \neq j \}.
  \end{align}
defined in Proposition \ref{prop:SingFun}. 
To see that the third property holds, notice first that by construction $ p(r,r') $ is the strength of $ (r,r') $ (see Proposition \ref{prop:isoClassicalAndSimplicialFuzzySets}). Combining this with Proposition \ref{prop:SingFun} yields $ p(r,r') = \exp(-d(r,r')) $ and all the identities for $ p(r,r') $ then follow from those for the metric. As for the fourth property, the an extract of the relevant diagram can be visualized as follows
\begin{center}
    \begin{tikzcd}[row sep=large]
        & {[X,d]^n} \arrow[ld] \arrow[d] \arrow[rd] & & & {(r_0,\dots,r_n)} \arrow[ld, maps to] \arrow[d, maps to] \arrow[rd, maps to] & \\
        {\overset{(i,j)}{[X,d]^1_a}} \arrow[d] \arrow[rd] & {\overset{(i,k)}{[X,d]^1_a}} \arrow[ld] \arrow[rd] & {\overset{(j,k)}{[X,d]^1_a}} \arrow[ld] \arrow[d] & {(r_i,r_j)} \arrow[d, maps to] \arrow[rd, maps to] & {(r_i,r_k)} \arrow[ld, maps to] \arrow[rd, maps to] & {(r_j,r_k)} \arrow[ld, maps to] \arrow[d, maps to] \\
        X & X & X & r_i & r_j & r_k
    \end{tikzcd}
\end{center}
To show that this is a limiting cone, suppose we are given another cone with apex $ Z $ and maps $ \zeta_{ij} : Z \to [X,d]^1_a $ and $ \zeta_i : Z \to X $. Assume that this cone factorizes through the canonical cone depicted above via a map $ q : Z \to [X,d]^n_a $. Then, this map is necessarily unique. For, given some $ z \in Z $, tracing its image $ q(z) = (q(z)_0, \dots, q(z)_n) \in [X,d]^n_a $ all the way down to the bottom row forces the components to be $ q(z)_i = \zeta_i(z) $. As for the existence, we can take this as a definition and set $ q(z) = (\zeta_0(z), \dots, \zeta_n(0)) $.

Assume now that we are given a fuzzy simplicial set $ S $ satisfying the listed properties. We construct an uber metric space $ (X,d) $ by taking its underlying set to be $ X := S(0,0) $. The metric on the other hand is given by
\begin{align}
    X \times X \to \mathbb{R}_{\ge 0} \cup \{\infty\}; \quad d(r,r') := - \log\bigl(p(r,r')\bigr).
\end{align}
Here, the fact that $ \bigl( S(\delta_1,i_{00}), S(\delta_0,i_{00}) \bigr) : S(1,0) \to S(0,0)^2 $ is an isomorphism guarantees that $ p(r,r') $ has a well-defined value in $ \bo{I} $, ensuring the well-definedness of the metric $ d $ as a function. The properties satisfied by $ p(r,r') $ immediately translate into making $ (X,d) $ an uber or extended-pseudo metric space.

It remains to show that we recover $ S $ by computing $ \te{Sing}(X,d) $. Alternatively we may compare $ S $ with $ [X,d] $. By definition $ [X,d]^0_a = X = S(0,0) = S(0,a) $ for all $ a \in \bo{I} $. Regarding the 1-simplices, we begin with the calculation:
\begin{align}
    \begin{split}
        [X,d]^1_a &= \{ (r, r') \in X\times X\, | \, d(r, r') \le - \log(a) \} \\
        &= \Big\{ (r,r') \in S(0,0)^2 \, \Big| \, p(r,r') \ge a \Big\} \\
        &= \Big\{ (r,r') \in S(0,0)^2 \, \Big| \,\sup\{ b \in \bo{I} \, | \, \langle r,r' \rangle \in S(1,b) \} \ge a \Big\} \\
        &\cong \Big\{ s \in S(1,0) \ \Big| \ \sup\{ b \in \bo{I} \, | \, s \in S(1,b) \} \ge a \Big\}.
    \end{split}
    \label{eq:XdIso}
\end{align}
where in the last line we used the fact that $ \langle -, - \rangle : S(0,0)^2 \to S(1,0) $ is invertible. 
Next, note that $\xi_1(s):=\sup\{ b \in \bo{I} \, | \, s \in S(1,b) \}$ is just the strength of $s$, as defined in eq.~\eqref{eq:maxStrength}. Therefore, considering $S(1,-)$ as a fuzzy set, the last line of \eqref{eq:XdIso} is precisely equal to the preimage $\xi_1^{-1}([a,1])$, which in turn is isomorphic to $S(1,a)$ by eq.~\eqref{eq:SaProperty}. As a consequence, we have
\begin{equation}
    \begin{split}
        [X,d]^1_a \cong \Big\{ s \in S(1,0) \ \Big| \ \sup\{ b \in \bo{I} \, | \, s \in S(1,b) \} \ge a \Big\} \cong S(1,a).
    \end{split}
\end{equation}

Both isomorphisms are natural in $ a $. 
The naturality of the first isomorphism holds by construction, while that of the other one follows from Proposition \ref{prop:isoClassicalAndSimplicialFuzzySets}.

A similar argument shows that the isomorphisms are also compatible with the face and degeneracy maps between the 0-simplices and the 1-simplices. This shows that up to level 1, $ S $ and $ [X,d] $ agree as fuzzy simplicial sets. More precisely, the 1-truncations of $ S $ and $ [X,d] $ are naturally isomorphic.

Since precomposition with $F_n$ restricts the diagrams $ S(F_n(-),a)$ and $ [X,d]^{F_n(-)}_a : \mathcal D_n^\te{op} \to \bo{Sets} $ such that they involve only 0-simplices and 1-simplices, the isomorphisms constructed so far assemble into a natural isomorphism $ S(F_n(-),a) \to [X,d]^{F_n(-)}_a $. Furthermore, since the application of $\te{lim}_{\mathcal{D}_n}$ is functorial (and hence preserves isomorphisms), this induces the bottom horizontal map in the diagram
\begin{center}
    \begin{tikzcd}
        {S(n,a)} \arrow[r, "\cong"] \arrow[d, "\cong"']       & {[X,d]^n_a} \arrow[d, "\cong"]        \\
        {\lim_{\mathcal D_n} S(F_n(-),a)} \arrow[r, "\cong"'] & {\lim_{\mathcal D_n}[X,d]^{F_n(-)}_a}
    \end{tikzcd}
\end{center}
The left vertical isomorphism holds by assumption (property 4.~of Proposition \ref{prop:SingSimpSetProperties}) and the right vertical isomorphism holds because $[X,d]$ is naturally isomorphic to $\te{Sing}$, which fulfills the properties of Proposition \ref{prop:SingSimpSetProperties}, too, as mentioned at the beginning of the proof. The top horizontal map is then defined so as to make the diagram commute. 

With these isomorphisms we now have a candidate for the natural isomorphism $ S \cong [X,d] $. It only remains to check whether the given family of maps is indeed natural. To see this, consider the morphisms $ (f,i_{ab}) : (n,a) \to (m,b) $ in $ \Delta \times \bo{I} $ and the diagram
\begin{center}
    \begin{tikzcd}
        {S(m,b)} \arrow[rd, "{S(f,i_{ab})}"] \arrow[rrr, "\cong"] \arrow[ddd, "\cong"'] & & & {[X,d]^m_b} \arrow[ld, "{[X,d](f,i_{ab})}"'] \arrow[ddd, "\cong"] \\
        & {S(n,a)} \arrow[r, "\cong"] \arrow[d, "\cong"'] & {[X,d]^n_a} \arrow[d, "\cong"] & \\
        & {\lim_{\mathcal D_n} S(F_n(-),a)} \arrow[r, "\cong"'] & {\lim_{\mathcal D_n}[X,d]^{F_n(-)}_a} & \\
        {\lim_{\mathcal D_m}S(F_m(-),b)} \arrow[ru] \arrow[rrr, "\cong "'] & & & {\lim_{\mathcal D_m}[X,d]^m_b} \arrow[lu]
    \end{tikzcd}
\end{center}
which we now explain. The outer and inner square are instances of the above commutative square. The left and right trapezoids commute by the naturality of the vertical maps (as explained above the proposition). The bottom trapezoid is induced by a square of natural transformations on the respective functors which commutes thanks to the fact that the candidate map $ S \cong [X,d] $ has already been verified to be natural up to dimension $ 1 $. Finally, a diagram chase argument now shows that the above trapezoid must also commute, which shows that the family of maps $ S(n,a) \to [X,d]^n_a $ is indeed natural and completes the proof.
\hfill$\blacksquare$

\subsection{Proof of Proposition \ref{prop:simplifySpivaksIdea}}
\label{proof:prop:simplifySpivaksIdea}

By Proposition \ref{prop:SingFun}, $\te{Sing}(M,d)(n,a)$ is equal to ordered $n+1$-tuples $[x_{i_0},\cdots,x_{i_n}]$ within $M$ with strength $\min_{j,k\in \{0,\cdots,n\}}\exp(-d(x_{i_j},x_{i_k}))$ at least equal to $a$. For $n=0$ and $n=1$, this is exactly equal to $G(M)$, defined in \cite{Spivak09}. We thus obtain $G(M)=\te{tr}_1(\te{Sing}(M))$.
Therefore, 
\begin{equation}
    \begin{split}
        \te{cosk}_1(G(M))(n,a) &= \te{Hom}_{\bo{1Fuz}}(\te{tr}_1(\Delta^n_{a}),G(M))  \\
        &= \te{Hom}_{\bo{1Fuz}}(\te{tr}_1(\Delta^n_{a}),\te{tr}_1(\te{Sing}(M))) \\
        &\overset{\te{(Proposition \ref{prop:leftRightAdjointOfTruncation})}}{\simeq} \te{Hom}_{\bo{sFuz}}(\te{sk}_1(\te{tr}_1(\Delta^n_{a})),\te{Sing}(M))\\
        &\overset{(\te{Re} ~\dashv ~ \te{Sing})}{\simeq} \te{Hom}_{\bo{UM}}(\te{Re}(\te{sk}_1(\te{tr}_1(\Delta^n_{a}))),M)\\
        &\overset{\te{(Proposition \ref{prop:realizationOfFuzzySimpSet})}}{\simeq} \te{Hom}_{\bo{UM}}(\te{Re}_{c1}\circ C_1\circ \te{tr}_1(\te{sk}_1(\te{tr}_1(\Delta^n_{a}))),M)\\
        &\overset{\eqref{eq:idTr}}{\simeq} \te{Hom}_{\bo{UM}}(\te{Re}_{c1}\circ C_1 \circ \te{tr}_1(\Delta^n_{a}),M)\\
        &\overset{\te{(Proposition \ref{prop:realizationOfFuzzySimpSet})}}{\simeq} \te{Hom}_{\bo{UM}}(\te{Re}(\Delta^n_{a}),M)\overset{\te{(def)}}{=} \te{Sing}(M)(n,a).
    \end{split}
\end{equation}
\hfill$\blacksquare$

\section{Proofs of section \ref{sec:TheMergeFunctors}}

\subsection{Proof of Proposition \ref{prop:tConormAnnihilatingElement}}
\label{proof:prop:tConormAnnihilatingElement}
We know that $T(a,0)=a$ for all $a$ by property $3$. In particular, this holds for $a=1$, so we obtain $T(1,0)=1$. However, we also have that $(a,b)\le (a',b')$ implies $T(a,b)\le T(a',b')$ by property $4$ (functoriality). And if we replace $(a,b)$ by $(1,0)$, which fulfills $(1,0)\le (1,c)$ for any $c$, then functoriality implies $T(1,c)\ge T(1,0)=1$. But since $T$ is upper-bounded by $1$, we must have $T(1,c)=1$ for all $c$. Finally, by property $2$ (symmetry), we obtain $T(c,1)=T(1,c)=1$ for all $c$ as well. 
\hfill$\blacksquare$

\subsection{Proof of Proposition \ref{prop:mergeIsAfunctor}}
\label{proof:prop:mergeIsAfunctor}

Functoriality of $\te{merge}_{\bo{cFuz}}$ follows from functoriality of the t-conorm (property 4 of Definition \ref{def:tConorm}). We have to show that $\te{merge}_{\bo{cFuz}}(f_1,f_2)$ is expansive, i.e.~that
\begin{equation}
    \begin{split}
        \xi(x) &\le \lambda(\te{merge}_{\bo{cFuz}}(f,f)(x)),\text{ where }\\
        \lambda(a) = T(\lambda_1(a),&\lambda_2(a)),\qquad \xi(a) = T(\xi_1(a),\xi_1(a))
    \end{split}
\end{equation}
By functoriality of $T$, we have $a\le a'$ and $b\le b'$ implies $T(a,b)\le T(a',b')$. Furthermore, morphisms in $\bo{cFuz}$ are expansive, i.e.~$\lambda_i(f(x))\ge\xi_i(x)$ for $i\in\{1,2\}$. Combining those two properties, we obtain 
\begin{equation}
    \begin{split}
        \xi(x) &= T(\xi_1(x),\xi_2(x))\le T(\lambda_1(f(x)),\lambda_2(f(x))) \\
        &= \lambda(f(x)) = \lambda(\te{merge}_{\bo{cFuz}}(f,f)(x)).
    \end{split}
\end{equation}
This shows that $\te{merge}_{\bo{cFuz}}(f,f)$ is indeed a well-defined morphism. 
What remains to show is that $\te{merge}_{\bo{cFuz}}$ preserves identities and compositions but that follows from the fact that they are implied by preservation on the underlying sets.
\hfill$\blacksquare$

\subsection{Proof of Proposition \ref{prop:associativity}}
\label{proof:prop:associativity}
Since both arguments of $\te{merge}_{\bo{Fuz}}$ are fuzzy sets with the same underlying classical set, we only need to show that symmetry and associativity hold for the strength of the merged elements. And this follows from symmetry and associativity of the t-conorm.
\hfill$\blacksquare$

\subsection{Proof of Proposition \ref{prop:mergeAllInCfuzAndApplyInverseOnce}}
\label{proof:prop:mergeAllInCfuzAndApplyInverseOnce}

Suppose that $s=\te{merge}_{\bo{Fuz}}(s_1,s_2)$ and suppose we want to merge it with $s_3$, where $s_1(0)=s_2(0)=s_3(0)=:X$. Then 
\begin{equation}
    \begin{split}
        \te{merge}_{\bo{Fuz}}(s,s_3)&= 
        C^{-1}(\te{merge}_{\bo{cFuz}}(C(s),C(s_3)))\\
        &= 
        C^{-1}(\te{merge}_{\bo{cFuz}}(\te{merge}_{\bo{cFuz}}(s_1,s_2),C(s_3)))\\
        &= C^{-1}(\te{merge}_{\bo{cFuz}}((X,T(\xi_1,\xi_2)),(X,\xi_3)))\\
        &= C^{-1}((X,T(T(\xi_1,\xi_2),\xi_3)))\\
        &= C^{-1}(\te{merge}_{\bo{cFuz}}(\{C(s_1),C(s_2),C(s_3)\})).
    \end{split}
\end{equation}
This argument can be iterated to yield the desired result.
\hfill$\blacksquare$

\subsection{Proof of Proposition \ref{prop:fuzzySimSetMergeWelldefined}}
\label{proof:prop:fuzzySimSetMergeWelldefined}

First of all, by definition, we can be sure that $\te{merge}_{\bo{sFuz}}(S_1,S_2)(n,-)$ is again a fuzzy set for every $n$. 

However, we have to show that $\te{merge}_{\bo{sFuz}}(S_1,S_2)(n,-)$ is
consistent with the compatibility conditions of a fuzzy simplicial set
across different $n$. Since all morphisms within $\Delta$ can be
decomposed into face and degeneracy maps (as explained after
Def. \ref{def:simplicialIndexingCategory}), it is sufficient to show
that compatibility conditions imposed on those are fulfilled. By
Proposition \ref{prop:simplexFaceStrength}, ``all degeneracies of a simplex
have the same strength as the simplex'', which of course does not change
upon application of the t-conorm, so there is nothing to prove for
degeneracy maps. By Proposition \ref{prop:simplexFaceStrength}, however,
``the strength of a simplex is at most equal to the minimum of the
strengths of its faces''. Since this statement is about maximal strengths,
it is simpler to check it by invoking the equivalence $C_s$ defined in
Def. \ref{def:csFuz}. For the corresponding classical fuzzy simplicial complex, the condition stated in Proposition \ref{prop:simplexFaceStrength} then can be expressed by Proposition \ref{prop:cFuzCompat} with equation \eqref{eq:compatCondCsFuz}.
And hence we have to check that 
\begin{equation}
    \begin{split}
        T(\xi_1(s),\xi_2(s))\le T(\xi_1(f_s),\xi_2(f_s))
    \end{split}
    \label{eq:condFuzSimp}
\end{equation}
for all $s\in C(S_1(n))=C(S_2(n))$ and all faces $f_s=C_s(S)(\delta)\in C(S_1(n-1))=C(S_2(n-1))$ of $s$. The condition \eqref{eq:condFuzSimp} is indeed satisfied because equation \eqref{eq:compatCondCsFuz} holds for both $C(S_1(n))=(S_1(n)(0),\xi_n^1)$ and $C(S_2(n))=(S_2(n)(0),\xi_n^2)$ 
and yields the desired result by inserting this into condition 4 (functoriality) of $T$ (cf.~Definition \ref{def:tConorm}).
\hfill$\blacksquare$

\subsection{Proof of Proposition \ref{prop:canMerge}}
\label{proof:prop:canMerge}

Applying the $\te{Sing}$-functor to each of them, we obtain a finite set $\{S_i:=\te{Sing}(X,d_i)\}_{i\in I}$. 
By Proposition \ref{prop:SingFun}, $S_i(0,0)\simeq X$ for all $i\in I$ and $S_i(n,0)=S_j(n,0)$ for all $i,j\in I$. As a consequence, condition \eqref{eq:conditionOnMergeSet} is fulfilled, and \eqref{eq:mergeSetFuzzySim} can be applied to merge them.
\hfill$\blacksquare$

\subsection{Proof of Proposition \ref{prop:mergeUM}}
\label{proof:prop:mergeUM}

By Proposition \ref{prop:canMerge}, the fuzzy simplicial sets $S_i:=\te{Sing}(X,d_i)$ can be merged using Definition \ref{def:simplicialfuzzyPowerSetMerge}.
Invoking this definition means that we recursively apply, for each $n\in\mathbb{N}$, equation \eqref{eq:fuzzySimpSetDef} in Definition \ref{def:mergeSfuz}, each of which, in turn, is defined via equation \eqref{eq:mergeFuz} in Definition \ref{def:mergeFuz}. Applying all this to the case $n$ for two fuzzy simplicial sets $S_1,S_2$, we obtain
\begin{equation}
    \begin{aligned}
        \te{merge}_{\bo{sFuz}}(S_1,S_2)(n,a)&\ov{\ref{eq:fuzzySimpSetDef}}\te{merge}_{\bo{Fuz}}(S_1(n),S_2(n))(a)\\
        &\ov{\ref{eq:mergeFuz}} C^{-1}(\te{merge}_{\bo{cFuz}}(C(S_1(n)),C(S_2(n))))(a)\\
        &\ov{Proposition \ref{prop:isoClassicalAndSimplicialFuzzySets}} (\xi_n^{(1,2)})^{-1}([a,1]),
    \end{aligned}
\end{equation}
where
\begin{equation}
        \xi_n^{(1,2)}(s:=[x_1,\cdots,x_n]) \ov{\ref{eq:classicalTconormMerging}} T(\xi^1_n(s),\xi^2_n(s)),
\end{equation}
where
\begin{equation}
    \xi_n^{i\in\{1,2\}}([x_0,\cdots,x_n]) \ov{Proposition \ref{prop:SingFun}}\te{min}_{j,k\in\{0,\cdots,n\}}\exp(-d_i(x_j,x_k)).
\end{equation}
In particular, for $n=0$, we obtain $\xi_0^{i}([x])=1~\forall x\in X,~i\in \{1,\cdots,N\}$, and, by Proposition \ref{prop:tConormAnnihilatingElement}, we have $T(1,1)=1$. Hence, merging all vertex sets simply returns a vertex set $X$ where all elements again have maximal strength $1$. 

For $n=1$, we obtain $\xi_1^{(1,2)}([x_1,x_2])=T(e^{-d_1(x_1,x_2)},e^{-d_2(x_1,x_2)})$.
Next, we invoke Proposition \ref{prop:mergeAllInCfuzAndApplyInverseOnce}, which tells us that we can first map all $S_i(1)$ to $C(S_i(1))$, perform $\te{merge}_{\bo{cFuz}}$ there and then apply the inverse $C^{-1}$ only once. We obtain
\begin{equation}
    \begin{split}
        \te{merge}_{\bo{Fuz}}\{S_i(1)\}&=
        C^{-1}(\te{merge}_{\bo{cFuz}}\{C(S_i(1))\}) = C^{-1}(X\times X,T_1),
    \end{split}
    \label{eq:eta2T1}
\end{equation}
where $T_1$ is given by equation \eqref{eq:T1}. Furthermore, by symmetry of the metrics $d_i$, $T_1$ is symmetric. 

Knowing the explicit form of $\te{merge}_{\bo{Fuz}}\{S_i(1)\}=\te{merge}_{\bo{sFuz}}\{S_i\}(1)$ and $\te{merge}_{\bo{sFuz}}\{S_i\}(0)$, we can finally invoke Proposition \ref{prop:realizationOfFuzzySimpSet} to conclude that the final result is a metric space $(X,d)$, where $X=S_i(0,0)$ is the underlying set shared by all initial metric spaces $(X,d_i)$ and $d$ is the metric defined by equation \eqref{eq:realizationExplicit} in the case of $\top_{\bo{UM}}$ and equation \eqref{eq:realizationExplicitEPMet} in the case of $\top_{\bo{EPMet}}$, while
\begin{equation}
    \begin{split}
        d_{\te{min}}(x,y)&=\min\{-\log(\eta_1(s))~|~\exists s~:~s=[x,y]\te{ or }s=[y,x]~\}\\
        &\ov{\ref{eq:eta2T1}} \min\{-\log(T_1(s))~|~\exists s~:~s=[x,y]\te{ or }s=[y,x]~\}\\
        &\ov{symmetry of $T_1$}-\log(T_1(x,y)).
    \end{split}
\end{equation}
\hfill$\blacksquare$

\subsection{Proof of Corollary \ref{cor:specialCase}}
\label{proof:cor:specialCase}

For the merge of two classical fuzzy edges, we obtain the strength
\begin{equation}
    \begin{split}
        &\quad-\log(\max(\exp(-d_1(x_i,x_j)),\exp(-d_2(x_i,x_j))))\\
        &=-\log(\exp(\max(-d_1(x_i,x_j),-d_2(x_i,x_j))))\\
        &=\min(d_1(x_i,x_j),d_2(x_i,x_j)).
    \end{split}
\end{equation}
Iterating this argument yields the result.
\hfill$\blacksquare$

\subsection{Proof of Corollary \ref{cor:knnMetrics}}
\label{proof:cor:knnMetrics}

Note that, for all $k$ with $k\ne i$ and $k\ne j$, we have $d_k(x_i,x_j)=\infty$ by equation \eqref{eq:localDists}, which implies $e^{-d_k(x_i,x_j)}=0$. At the same time $T(a,0)=a$ by Definition \ref{def:tConorm}. As a result,
\begin{equation}
\begin{split}
    T_1(x_i,x_{j})&=T(e^{-d_1(x_i,x_{j})},T_2(x_i,x_{j})) = T(0,T_2(x_i,x_{j}))\\
    & = T_2(x_i,x_{j}) = \cdots = T(e^{-d_i(x_i,x_j)},e^{-d_j(x_i,x_j)}).
\end{split}
\end{equation}
Finally, by symmetry of the metric, $d_j(x_i,x_j)=d_j(x_j,x_i)$.
\hfill$\blacksquare$

\end{document}